\documentclass{amsart}
\usepackage{amsmath}
  \usepackage{paralist}
  \usepackage{amssymb}
 \usepackage{amsthm}
 \usepackage{esint}
 \usepackage[colorlinks=true]{hyperref}
\hypersetup{urlcolor=blue, citecolor=red}

\newtheorem{theorem}{Theorem}[section]
\newtheorem{corollary}[theorem]{Corollary}
\newtheorem{lemma}[theorem]{Lemma}
\newtheorem{proposition}[theorem]{Proposition}
\theoremstyle{definition}

\newtheorem{remark}[theorem]{Remark}

\numberwithin{equation}{section}

\title[Asymptotic profiles for Choquard equations]
      {Asymptotic profiles for Choquard equations with general critical nonlinearities}

\author[ Xiaonan Liu, Shiwang Ma, Yachen Wang]{}



\email{liuxiaonan20131110@163.com (X.N. Liu), shiwangm@nankai.edu.cn (S.W. Ma)}
\email{ yachen\_wang@mail.nankai.edu.cn (Y.C. Wang).}

\thanks{This work is partly supported by the National Natural Science Foundation of China (No. 12301128).\\
  {$^{\rm *}$ Corresponding authors: shiwangm@nankai.edu.cn, yachen\_wang@mail.nankai.edu.cn}
}

\begin{document}
\maketitle

\centerline{\scshape Xiaonan Liu$^{\rm a}$,\;  Shiwang Ma $^{\rm b \rm *}$, Yachen Wang $^{\rm b \rm *}$}
\medskip
{\footnotesize
 \centerline{$^{\rm a}$ School of Science}
 \centerline{Beijing University of Posts and Telecommunications, Beijing 100876, China}
  \centerline{$^{\rm b}$ School of Mathematical Sciences and LPMC}
 \centerline{Nankai University, Tianjin 300071, China}
 }

\bigskip

\begin{abstract}
In this paper, we study asymptotic behavior of positive ground state solutions for the nonlinear Choquard equation:
\begin{equation}\label{0.1}
-\Delta u+\varepsilon u=\big(I_{\alpha}\ast F(u)\big)F'(u),\quad
u\in H^1(\mathbb R^N),
\end{equation}
where $F(u)=|u|^{\frac{N+\alpha}{N-2}}+G(u)$, $N\geq3$ is an integer,  $I_{\alpha}$ is the Riesz potential of order $\alpha\in(0,N)$, and $\varepsilon>0$ is a parameter.  Under some mild subcritical growth assumptions on $G(u)$, we show that
as $\varepsilon \to \infty$, the ground state solutions of \eqref{0.1}, after a suitable rescaling, converge to a particular solution of the critical Choquard equation $-\Delta u=\frac{N+\alpha}{N-2}(I_{\alpha}*|u|^{\frac{N+\alpha}{N-2}})|u|^{\frac{N+\alpha}{N-2}-2}u$. We establish a novel sharp asymptotic characterisation of such a rescaling, which depends in a non-trivial way on  the asymptotic behavior of $G(u)$ at infinity and the space dimension $N=3$, $N=4$ or $N\geq5$.
\end{abstract}

\section{Introduction and Main Results}

In this paper, we consider the following Choquard equation with general nonlinearity
\begin{equation}\label{1.1}
-\Delta u+\varepsilon u=(I_{\alpha}\ast F(u))F'(u),\quad
u\in H^1(\mathbb R^N),
\end{equation}
where $N\geq3$ is an integer, $\varepsilon>0$ is a parameter, $I_{\alpha}$ is the Riesz potential of order $\alpha\in(0,N)$ and  is defined for every
$x\in \mathbb R^N\backslash \{0\}$ by
$$I_{\alpha}(x)=\frac{A_{\alpha}(N)}{|x|^{N-\alpha}}\ \ \mbox{and }\ \ A_{\alpha}(N)=\frac{\Gamma(\frac{N-\alpha}{2})}{\Gamma(\frac{\alpha}{2})\pi^{\frac{N}{2}}2^{\alpha}},$$ where $\Gamma$ denotes the Gamma function. We always  assume that  $F(u)=|u|^{\frac{N+\alpha}{N-2}}+G(u)$ with $G\in C^1(\mathbb R,\mathbb R)$ being a subcritical nonlinearity which is specified later. Particularly, $G(u)$ includes the following combined  powers nonlinearity of the form
\begin{equation}\label{1.100}
\sum_{i=1}^kc_i|u|^{q_i}, \quad c_i>0, \quad q_i\in \left(\frac{N+\alpha}{N}, \frac{N+\alpha}{N-2}\right).
\end{equation}

A prototype of such equations comes from the research of standing-wave solutions of the nonlinear Schr\"{o}dinger equation with attractive combined nonlinearities
\begin{equation}\label{1.2}
\begin{aligned}
i\psi_{t}-\Delta \psi=&(I_{\alpha}*(|\psi|^{\frac{N+\alpha}{N-2}}+|\psi|^{q_i}+|\psi|^{q_j}))\\
\cdot &(\frac{N+\alpha}{N-2}|\psi|^{\frac{N+\alpha}{N-2}-2}\psi+q_i|\psi|^{q_i-2}\psi+q_j|\psi|^{q_j-2}\psi),
\end{aligned}
\end{equation}
in $\mathbb R^N\times\mathbb R$.
One makes the ansatz $\psi(t,x)=e^{-i\varepsilon t}u(x)$ in  \eqref{1.2}, where $u: \mathbb R^N\rightarrow \mathbb C$, then  \eqref{1.2} reduces to the equation \eqref{1.1} with
$G(u)=|u|^{q_i}+|u|^{q_j}$. A theory about NLS with combined power nonlinearities was first developed by T. Tao, M. Visan and X. Zhang \cite{TVZ2007}
and  then received much attention during the last decades (cf. \cite{AIIKN2019, AIKN, LM2019, LMZ2019, LM2022}). For the background of the problem \eqref{1.1} we refer the readers to \cite{MS2017} and the references therein.

A solution $u_{\varepsilon}\in H^1(\mathbb R^N)$ of \eqref{1.1} is a critical point of the corresponding Action functional defined by
$$
I_{\varepsilon}(u)=\frac{1}{2}\int_{\mathbb R^N}|\nabla u|^2+\frac{\varepsilon}{2}\int_{\mathbb R^N}|u|^2
-\frac{1}{2}\int_{\mathbb R^N}(I_{\alpha}*(|u|^{\frac{N+\alpha}{N-2}}+G(u)))
(|u|^{\frac{N+\alpha}{N-2}}+G(u)).
$$
\smallskip

In \cite{MM2014}, Moroz and  Muratov first  study the asymptotic properties of ground states for  a class of scalar field equations   with a defocusing large exponent $p$ and a focusing smaller exponent $q$. More precisely, the following equation
\begin{equation}\label{1.5a}
-\Delta u+\varepsilon u=|u|^{p-2}u-|u|^{q-2}u, \quad {\rm in} \ \mathbb R^N,
\end{equation}
is discussed, where $N\ge 3$, $q>p>2$.  Later, in \cite{Lewin-1},  M. Lewin and S. Nodari prove a general result about the uniqueness and non-degeneracy of positive radial solutions to the above equation, and then the non-degeneracy of the unique solution $u_\varepsilon$ allows them to derive its behavior in the two limits $\varepsilon\to 0$ and $\varepsilon\to \varepsilon_*$, where $\varepsilon_*$ is a  threshold of existence.  Amongst other things, a precise asymptotic  expression of $M(\varepsilon)=\|u_\varepsilon\|_2^2$ is obtained. This gives the uniqueness of energy minimizers at fixed mass in certain regimes.

  In \cite{LM2022}, Z. Liu and Moroz extend the results in \cite{MM2014} to a class of Choquard type equation
\begin{equation}\label{1.5b}
-\Delta u+\varepsilon u=(I_\alpha \ast |u|^{p})|u|^{p-2}u- |u|^{q-2}u,
\quad {\rm in} \ \mathbb R^N,
\end{equation}
where $N\ge 3$ is an integer. Under some assumptions on the exponents $p$ and $q$, the limit profiles of the ground states are discussed in  the two cases $\varepsilon\to 0$ and $\varepsilon\to \infty$. But the precisely
asymptotic behaviors of ground states remain open.

\smallskip

In \cite{AIIKN2019, AIKN}, T.   Akahori et al. consider the following Schr\"odinger equation with two focusing  exponents $p$ and $q$:
\begin{equation}\label{1.5bbb}
-\Delta u+\varepsilon u=|u|^{p-2}u+|u|^{q-2}u, \quad {\rm in} \ \mathbb R^N,
\end{equation}
where $N\ge 3$, $2<q<p\le 2^*$ and $\varepsilon>0$ is a parameter.  It is well-known \cite{ASM2012, LLT2017, ZZ2012} that when $p=2^*$ and
\begin{equation}\label{1.40}
 q\in \begin{cases}
\begin{split}
&(2, 2^*),\ &if\ N\ge 4,\\
& (4,6),\ &if\ N=3,
\end{split}
\end{cases}
\end{equation}
 the equation \eqref{1.5bbb} admits a ground state solution $u_\varepsilon$ for any $\varepsilon>0$. The authors in \cite{AIKN}   proved that for any $q\in (2,2^*)$ and  small $\varepsilon>0$, the ground state is unique and as $\varepsilon\to 0$,  the unique ground state $u_\varepsilon$ tends to the unique positive solution of the equation $-\Delta u+u=u^{q-1}$. After a suitable rescaling, authors in \cite{AIIKN2019}  established a uniform decay estimate for the ground state $u_\varepsilon$, and then for any $q$ verifying \eqref{1.40},  they  proved  the uniqueness and nondegeneracy in $H^1_r(\mathbb R^N)$ of ground states $u_\varepsilon$   for $N\ge 5$ and large $\varepsilon>0$,  and show that for $N\ge 3$, as $\varepsilon\to\infty$, $u_\varepsilon$ tends to a particular solution of the critical Emden--Fowler equation. Later, for any  $q$ satisfying \eqref{1.40}, T. Akahori, M. Murata \cite{AM-1,AM-2}  established the uniqueness and  nondegeneracy in $H^1_r(\mathbb R^N)$ of ground states $u_\varepsilon$  in the cases $N=3,4$.
On the other hand, T. Akahori et al.\cite{ Akahori-4} and  J. Wei and Y. Wu \cite{Wei-2} proved that for $N=3, q\in (2,4]$, the equation \eqref{1.5bbb} has no ground state solutions for large $\varepsilon>0$.  By using a global bifurcation argument, Jeanjean, Zhang and Zhong \cite{Jeanjean-4} also study the asymptotic behaviour of solutions as $\varepsilon\to 0$ and $\varepsilon\to \infty$  for the equation  \eqref{1.5bbb} with a general subcritical nonlinearity and discuss the connection  to the existence, non-existence and multiplicity of  prescribed mass positive solutions to \eqref{1.5bbb} with  the  associated $L^2$ constraint condition $\int_{\mathbb R^N}|u|^2=a^2$. But the precise asymptotic profiles of positive solutions are not addressed there.  For  other related papers, we refer the reader to \cite{MM2023} and the references therein.

Li, Ma and Zhang \cite{LMZ2019}  and Li and Ma \cite{LM2019} considered the following nonlinear Choquard equation
\begin{equation}\label{1.4}
-\Delta u+ \varepsilon u=|u|^{p-2}u+(I_\alpha\ast |u|^q)|u|^{q-2}u \ \ \mbox{in}\ \ \mathbb R^N,
\end{equation}
where $N\geq3$ is an integer, $p\in \left(2,\frac{2N}{N-2}\right]$, $q\in \left[\frac{N+\alpha}{N},\frac{N+\alpha}{N-2}\right]$,  $I_{\alpha}$ is the Riesz potential with $\alpha\in(0,N)$ and $\varepsilon>0$ is a parameter, they established the existence, regularity and symmetry of the ground state solutions of \eqref{1.4}. Especially, for any $\varepsilon>0$, the equation \eqref{1.4} admits a ground state solution if $p=\frac{2N}{N-2}$ and
\begin{equation}\label{1.41}
q\in \begin{cases}
\begin{split}
&\left(1+\frac{\alpha}{N-2}, \frac{N+\alpha}{N-2}\right),\ &if\ N\ge 4,\\
&\ \ \ (2+\alpha, 3+\alpha),\ &if\ N=3.
\end{split}
\end{cases}
\end{equation}
It has been shown by the second author \cite{M-1} that if the above relation is not true and $q$ is a subcritical exponent, then  the equation \eqref{1.4} has no ground state solution for large $\varepsilon>0$, and admits two positive solutions for $\varepsilon>0$ sufficiently small under  some additional technique conditions.  On the other hand, if \eqref{1.41} holds, the second author and Moroz \cite{MM}  establish the precisely asymptotic behaviors of positive ground states of the equation \eqref{1.4} as $\varepsilon\to 0$ and $\varepsilon\to \infty$ and discuss the connection  to the existence and multiplicity of  prescribed mass positive solutions to \eqref{1.4} with  the  associated $L^2$ constraint condition $\int_{\mathbb R^N}|u|^2=a^2$.  More recently,  by employing a new technique, the second author and Moroz \cite{MM2024} extend some results in \cite{MM} to the following critical Choquard equation   with a general local perturbation:
\begin{equation}\label{1.51}
-\Delta u+ \varepsilon u=g(u)+(I_\alpha\ast |u|^p)|u|^{p-2}u \ \ \mbox{in}\ \ \mathbb R^N,
\end{equation}
where $N\ge 3$, $\varepsilon>0$ is a parameter, $I_{\alpha}$ is the Riesz potential of order $\alpha\in(0,N)$. Among other things, if  $p=\frac{N+\alpha}{N-2}$ and $\varepsilon\to \infty$, the precise asymptotic behavior of ground state solutions of the equation \eqref{1.51} is  obtained under the following mild  assumptions:

 \smallskip
{\bf (M1)} \ $g\in C([0,\infty), [0,\infty) )$ satisfies $g(s)=o(s)$ as $s\to 0$;

\smallskip
{\bf (M2)} \ there exists $q\in (2,2^*)$ verifying \eqref{1.40} such that $\lim_{s\to \infty}\frac{g(s)}{s^{q-1}}=A>0.$
\smallskip

\smallskip

Moroz and van Schaftingen \cite{Moroz-2, MS2015} have established a general theory for  the  general Choquard equation \eqref{1.1} with $F(u)$ bing a subcritical nonlinearity in the sense of Hardy-Littlewood-Sobolev inequality. More precisely, under the  assumptions:
$$
|sF'(s)|\lesssim |s|^{\frac{N+\alpha}{N}}+|s|^{\frac{N+\alpha}{N-2}},\quad
\lim_{s\to 0}\frac{F(s)}{|s|^{\frac{N+\alpha}{N}}}=\lim_{s\to\infty}\frac{F(s)}{|s|^{\frac{N+\alpha}{N-2}}}=0;
$$
and there exists $s_0\in \mathbb R\setminus \{0\}$ such that $F(s_0) \not= 0$, Moroz and van Schaftingen \cite{MS2015} proved the existence of a ground state solution of \eqref{1.1}.
Usually,$ \frac{N+\alpha}{N}$ is called the lower critical exponent and $\frac{N+\alpha}{N-2}$ is the upper critical exponent for the Choquard equation.

When $F(u)=|u|^{\frac{N+\alpha}{N-2}}$, then it follows  from the associated Poho\v zaev identity that \eqref{1.1} has no nontrivial solution in $H^1(\mathbb R^N)$ \cite{Moroz-2}. Nonetheless, there exist some  results about the existence and properties of ground states for the Choquard equations with critical nonlinearities \cite{LLT2020, LW-1}.
It turns out that the existence of ground state solutions of the equation \eqref{1.1} depends heavily on  the asymptotic behavior of $G(u)$ at infinity. For example, Li and Wang \cite{LW-1}  study  the existence and symmetry of ground state solutions for \eqref{1.1} with  $G(u)=|u|^q$ and $q\in [\frac{N+\alpha}{N}, \frac{N+\alpha}{N-2}).$ Among other things, Li and Wang \cite{LW-1} proved that for any $\varepsilon>0$,  the equation \eqref{1.1} admits at least one  ground state solution $u\in H^1(\mathbb R^N)$ if
$F(u)=|u|^{\frac{N+\alpha}{N-2}}+|u|^q$ and
\begin{equation}\label{1.42}
q\in \begin{cases}
\begin{split}
&\left(\frac{N+\alpha}{N}, \frac{N+\alpha}{N-2}\right),\ &if\ N\ge 4,\\
&\ \ \ (1+\alpha, 3+\alpha),\ &if\ N=3.
\end{split}
\end{cases}
\end{equation}
Clearly, this result is different in nature from that of the equation  \eqref{1.4}.  For other related results, we refer the readers to the references cited in \cite{LLT2020, LW-1}.

\smallskip

 However, so far nothing has been done for the asymptotic behavior of ground state solutions of  general Choquard equation \eqref{1.1}.
Inspired by the results in \cite{MM, MM2024}, in the present paper, we consider the asymptotic behavior of positive ground state solutions for the nonlinear Choquard equation \eqref{1.1}. To this end,  we make  the following  assumptions:




\smallskip
\begin{itemize}
\item[\bf{(H1)}] $G\in C^1(\mathbb R, \mathbb R)$ and $G(s)>0$ for $s>0$;
\smallskip
\item[\bf{(H2)}] there exist  $\frac{N+\alpha}{N}<q_1\le q_2<\frac{N+\alpha}{N-2}$
such that $q_1G(s)\leq G'(s)s \leq q_2G(s)$ and
\begin{equation}\label{a1}
\lim_{s\to 0}\frac{G(s)}{|s|^{q_1}}=a,\ \ \lim_{s\to \infty}\frac{G(s)}{|s|^{q_2}}=b,
\end{equation}
where $a$ and $b$ are positive constants.
\end{itemize}

\smallskip

Clearly, the function $G(u)=\sum_{i=1}^kc_i|u|^{q_i}$ given in \eqref{1.100} satisfies (H1) and (H2). In what follows, we denote $g(u):=G'(u)$, then it immediately follows  from (H1) and (H2) that,
there exist constants $0<c\leq q_1\leq q_2\leq C$ such that
\begin{equation}\label{1.5}
c(s^{q_1}+s^{q_2})\leq q_1G(s)\leq g(s)s \leq q_2G(s)\leq C(s^{q_1}+s^{q_2}).
\end{equation}
It is well-known \cite{A1989,T1976} that the  best Sobolev constant $S$ can be achieved  by the Talenti function $U_1$ given by
$$U_1(x):=[N(N-2)]^{\frac{N-2}{4}}\left(\frac{1}{1+|x|^2}\right)^{\frac{N-2}{2}}.$$
We denote $C_*(N,\alpha)=\left(\frac{1}{2\sqrt{\pi}}\right)^{\alpha}
\frac{\Gamma\left(\frac{N-\alpha}{2}\right)}{\Gamma\left(\frac{N+\alpha}{2}\right)}\left\{\frac{\Gamma(N)}{\Gamma\left(\frac{N}{2}\right)}\right\}^{\alpha}.$
Then it is known  \cite{L1983} that the function $V_1(x)=S^{-\frac{(N-2)\alpha}{4(2+\alpha)}}C_*(N,\alpha)^{-\frac{N-2}{2(2+\alpha)}}U_1(x)$
is a minimizer of the problem
\begin{equation}\label{1}
S_{\alpha}:=\inf_{u\in \mathcal{D}^{1,2}(\mathbb R^N)\setminus \{0\}}\frac{\int_{\mathbb R^N}|\nabla u|^2dx}
{\left(\int_{\mathbb R^N}(I_{\alpha}*|u|^{\frac{N+\alpha}{N-2}})|u|^{\frac{N+\alpha}{N-2}}\right)^{\frac{N-2}{N+\alpha}}}.
\end{equation}
It is also easy  to see that the family of functions
$$
W_1(x)=\left(\frac{N-2}{N+\alpha}\right)^{\frac{N-2}{2(2+\alpha)}}V_1(x), \quad W_{\rho}(x)=\rho^{-\frac{N-2}{2}}W_1(x/\rho), \quad \rho>0
$$
solve the following critical Choquard equation
\begin{equation}\label{1.10}
-\Delta w=\frac{N+\alpha}{N-2}(I_{\alpha}*|w|^{\frac{N+\alpha}{N-2}})|w|^{\frac{N+\alpha}{N-2}-2}w.
\end{equation}

\smallskip


 In this paper,  we are interested in  the limit asymptotic
profile of the ground-states  $u_{\varepsilon}$ of the problem \eqref{1.1}, and the asymptotic behavior of different norms of $u_{\varepsilon}$, as
$\varepsilon \to \infty$. Our main result is as follows:
\smallskip

\begin{theorem}\label{th1.3}
 Assume that (H1) and (H2) hold.  Suppose further that
$$
q_1>2,\quad  if\ N=4,\  3,
\leqno{\bf (H3)}
$$
and
$$
q_2>\begin{cases}
\begin{split}
&2,\ &if\ N=4,\\
&\max\{2, 1+\alpha\},\ &if\ N=3,\\
\end{split}
\end{cases}
\leqno{\bf (H4)}
$$
then \eqref{1.1} admits a positive ground state $u_\varepsilon\in H^1(\mathbb R^N)$ which is radially symmetric and radially nonincreasing, and
for  large  $\varepsilon>0$, there exists $\xi_{\varepsilon} \in (0,\infty)$ satisfying
$$
\xi_{\varepsilon}
\thicksim
\begin{cases}
\begin{split}
&\varepsilon^{-\frac{2}{4+(N-2)q_2-(N+\alpha)}}\ &if \ N\ge 5,\\
&\left(\varepsilon \ln \varepsilon\right)^{-\frac{2}{2q_2-\alpha}},\ &if\ N=4,\\
&\varepsilon^{-\frac{1}{q_2-1-\alpha}},\ &if\ N=3,\\
\end{split}
\end{cases}
$$
such that the rescaled family of ground states $w_{\varepsilon}(x)=\xi_{\varepsilon}^{\frac{N-2}{2}}u_{\varepsilon}(\xi_{\varepsilon}x)$ satisfies
$$\|\nabla w_{\varepsilon}\|_2^2
\thicksim\int_{\mathbb R^N}(I_{\alpha}*|w_{\varepsilon}|^{\frac{N+\alpha}{N-2}})|w_{\varepsilon}|^{\frac{N+\alpha}{N-2}}\thicksim 1, \quad
\|w_{\varepsilon}\|_2^2
\thicksim
\begin{cases}
\begin{split}
&1,\ &if\ N\ge 5,\\
&\ln \varepsilon,\ &if\ N=4,\\
&\varepsilon^{\frac{3+\alpha-q_2}{2(q_2-1-\alpha)}},\ &if\ N=3,
\end{split}
\end{cases}
$$
 and as $\varepsilon \to \infty$, for $N\geq5$, $w_{\varepsilon}$ converges to $W_{\rho_0}$ in $H^1(\mathbb R^N)$, where
$$
\rho_0=\left(\frac{b[(N+\alpha)-(N-2)q_2]}{2}
\frac{\int_{\mathbb R^N}(I_{\alpha}*|W_1|^{\frac{N+\alpha}{N-2}})|W_1|^{q_2}}{\int_{\mathbb R^N}|W_1|^2}\right)^{\frac{2}{4+(N-2)q_2-(N+\alpha)}},
$$
for  $N=3,4$,  $w_{\varepsilon}$ converges to $W_{1}$ in $D^{1,2}(\mathbb R^N)$ and $L^{s}(\mathbb R^N)$ for any $s>\frac{N}{N-2}$.
Furthermore, the least energy $m_{\varepsilon}$
of the ground state satisfies
$$
\frac{2+\alpha}{2(N-2)}\left(\frac{N-2}{N+\alpha}S_{\alpha}\right)^{\frac{N+\alpha}{2+\alpha}}-m_{\varepsilon}
\thicksim
\begin{cases}
\begin{split}
&\varepsilon^{-\frac{(N+\alpha)-(N-2)q_2}{4+(N-2)q_2-(N+\alpha)}},\ &if\ N\ge 5,\\
&\left(\varepsilon\ln \varepsilon\right)^{-\frac{4+\alpha-2q_2}{2q_2-\alpha}},\ &if\ N=4,\\
&\varepsilon^{-\frac{3+\alpha-q_2}{2(q_2-1-\alpha)}},\ &if\ N=3.\\
\end{split}
\end{cases}
$$
\end{theorem}

\smallskip

We remark that $q_i>2 \ (i=1,2)$ in (H3) and (H4) are technical  conditions.  While, in the case $N=3$, the condition $q_2>1+\alpha$  in (H4) is essential in our main result. The situation $N=3$ and $q_2\le 1+\alpha$  will be treated in a forthcoming paper.
\smallskip

As a direct consequence of Theorem 1.1, we have the following result.

\begin{corollary}
 Assume that (H1), (H2), (H3) and (H4) hold.  If   $u_\varepsilon\in H^1(\mathbb R^N)$ is the ground state solution of \eqref{1.1}, then as $\varepsilon\to \infty$, there hold
$$
u_{\varepsilon}(0)
\thicksim
\begin{cases}
\begin{split}
&\varepsilon^{\frac{N-2}{4+(N-2)q_2-(N+\alpha)}},\ \ &if\ \ N\geq5,\\
&\left(\varepsilon \ln \varepsilon\right)^{\frac{2}{2q_2-\alpha}},\ \ &if\ \ N=4,\\
&\varepsilon^{\frac{1}{2(q_2-1-\alpha)}},\ \ &if\ \ N=3,
\end{split}
\end{cases}
$$
$$
\|u_{\varepsilon}\|_2^2
\thicksim
\begin{cases}
\begin{split}
&\varepsilon^{-\frac{4}{4+(N-2)q_2-(N+\alpha)}},\ \ &if\ \ N\geq5,\\
&\varepsilon^{-\frac{4}{2q_2-\alpha}}\left(\ln \varepsilon\right)^{-\frac{4+\alpha-2q_2}{2q_2-\alpha}},\ \ &if\ \ N=4,\\
&\varepsilon^{-\frac{q_2+1-\alpha}{2(q_2-1-\alpha)}},\ \ &if\ \ N=3,
\end{split}
\end{cases}
$$

$$
\|\nabla u_{\varepsilon}\|_2^2
=\left(\frac{N-2}{N+\alpha}\right)^{\frac{N-2}{2+\alpha}}
S_{\alpha}^{\frac{N+\alpha}{2+\alpha}}+
\begin{cases}
\begin{split}
&O(\varepsilon^{-\frac{(N+\alpha)-(N-2)q_2}{4+(N-2)q_2-(N+\alpha)}}),\ \ &if\ \ N\geq5,\\
&O((\varepsilon\ln \varepsilon)^{-\frac{4+\alpha-2q_2}{2q_2-\alpha}}),\ \ &if\ \ N=4,\\
&O(\varepsilon^{-\frac{3+\alpha-q_2}{2(q_2-1-\alpha)}}),\ \ &if\ \ N=3,
\end{split}
\end{cases}
$$
$$
\begin{aligned}
\int_{\mathbb R^N}(I_{\alpha}*|u_{\varepsilon}|^{\frac{N+\alpha}{N-2}})|u_{\varepsilon}|^{\frac{N+\alpha}{N-2}}
&=\left(\frac{N-2}{N+\alpha}S_{\alpha}\right)^{\frac{N+\alpha}{2+\alpha}}\\
&\quad\quad\quad+
\begin{cases}
\begin{split}
&O(\varepsilon^{-\frac{(N+\alpha)-(N-2)q_2}{4+(N-2)q_2-(N+\alpha)}}),\ \ &if\ \ N\geq5,\\
&O((\varepsilon\ln \varepsilon)^{-\frac{4+\alpha-2q_2}{2q_2-\alpha}}),\ \ &if\ \ N=4,\\
&O(\varepsilon^{-\frac{3+\alpha-q_2}{2(q_2-1-\alpha)}}),\ \ &if\ \ N=3.
\end{split}
\end{cases}
\end{aligned}
$$
\end{corollary}


\smallskip

\noindent
\textbf{Organization of the paper}. In Section 2, we give  some preliminary results which are needed in the proof of our main result.  Sections 3--4 are devoted to the proof of Theorem 1.1. Finally, In the Appendix, we  discuss the existence  and their optimal uniform decay estimates of  ground states and  outline the proofs of  Proposition \ref{Pro 4.1} and  the existence of ground state solution to the equation \eqref{1.1} under our assumptions.

\smallskip

\noindent
\textbf{Basic notations}. Throughout this paper, we assume $N\geq
3$. $L^p(\mathbb{R}^N)$ with $1\leq p<\infty$ denotes the Lebesgue space
with the norms
$\|u\|_p=\left(\int_{\mathbb{R}^N}|u|^p\right)^{1/p}$.
 $ H^1(\mathbb{R}^N)$ is the usual Sobolev space with norm
$\|u\|_{H^1(\mathbb{R}^N)}=\left(\int_{\mathbb{R}^N}|\nabla
u|^2+|u|^2\right)^{1/2}$. $ D^{1,2}(\mathbb{R}^N)=\{u\in
L^{2^*}(\mathbb{R}^N): |\nabla u|\in
L^2(\mathbb{R}^N)\}$. $H_r^1(\mathbb{R}^N)=\{u\in H^1(\mathbb{R}^N):
u\  \mathrm{is\ radially \ symmetric}\}$.  $B_r$ denotes the ball in $\mathbb R^N$ with radius $r>0$ and centered at the origin,  $|B_r|$ and $B_r^c$ denote its Lebesgue measure  and its complement in $\mathbb R^N$, respectively.  As usual, $C$, $c$, etc., denote generic positive constants.
 For any  large $\epsilon>0$ and two nonnegative functions $f(\epsilon)$ and  $g(\epsilon)$, we write

(1)  $f(\epsilon)\lesssim g(\epsilon)$ or $g(\epsilon)\gtrsim f(\epsilon)$ if there exists a positive constant $C$ independent of $\epsilon$ such that $f(\epsilon)\le Cg(\epsilon)$.

(2) $f(\epsilon)\sim g(\epsilon)$ if $f(\epsilon)\lesssim g(\epsilon)$ and $f(\epsilon)\gtrsim g(\epsilon)$.\\
If $|f(\epsilon)|\lesssim |g(\epsilon)|$, we write $f(\epsilon)=O((g(\epsilon))$. Finally, if $\lim_{\epsilon\to \infty} f(\epsilon)/g(\epsilon)=0$, then we write $f(\epsilon)=o(g(\epsilon))$ as $\epsilon\to \infty$.

\section{Preliminary}
In this section, we mainly give some preliminaries that will be used in this  paper.
To describe the  asymptotic behavior of $u_{\varepsilon}$, we use a canonical rescaling:
\begin{equation}\label{1.6}
v(x)=\varepsilon^{-\frac{N-2}{4}}u(\varepsilon^{-\frac{1}{2}} x),
\end{equation}
then \eqref{1.1} transforms into the equation
\begin{equation}\label{1.7}
\begin{aligned}
-\Delta v+ v=&\left(I_{\alpha}*(|v|^{\frac{N+\alpha}{N-2}}+\varepsilon^{-\frac{N+\alpha}{4}}G(\varepsilon^{\frac{N-2}{4}}v))\right)\\
\cdot &\left(\frac{N+\alpha}{N-2}|v|^{\frac{N+\alpha}{N-2}-2}v+\varepsilon^{-\frac{N+\alpha}{4}+\frac{N-2}{4}}g(\varepsilon^{\frac{N-2}{4}}v)\right).\\
\end{aligned}
\end{equation}

The formal limit equation for \eqref{1.7} as $\varepsilon\to\infty$ is given by
\begin{equation}\label{1.7a}
-\Delta v+ v=\frac{N+\alpha}{N-2}(I_{\alpha}*|v|^{\frac{N+\alpha}{N-2}})|v|^{\frac{N+\alpha}{N-2}-2}v.
\end{equation}
Recall that the critical Choquard equation \eqref{1.7a} has no nontrivial solutions in $H^1(\mathbb R^N)$, which  follows from the associated  Poho\v{z}aev's identity. We denote the Nehari manifolds for \eqref{1.7} and
\eqref{1.7a} as follows:
$$\mathcal{N}_{\varepsilon}^*:=\{v \in H^1(\mathbb R^N)\backslash\{0\}\mid N_{\varepsilon}^*(v)=0\},$$
and
$$\mathcal{N}_{\infty}^*:=\{v \in H^1(\mathbb R^N)\backslash\{0\}\mid
\int_{\mathbb R^N}|\nabla v|^2+\int_{\mathbb R^N}|v|^2=\frac{N+\alpha}{N-2}(I_{\alpha}*|v|^{\frac{N+\alpha}{N-2}})|v|^{\frac{N+\alpha}{N-2}}\},$$
respectively, where
$$\begin{aligned}
N_{\varepsilon}^*(v)=\int_{\mathbb R^N}|\nabla v|^2+\int_{\mathbb R^N}|v|^2
-\int_{\mathbb R^N}&(I_{\alpha}*(|v|^{\frac{N+\alpha}{N-2}}
+\varepsilon^{-\frac{N+\alpha}{4}}G(\varepsilon^{\frac{N-2}{4}}v)))\\
\cdot&(\frac{N+\alpha}{N-2}|v|^{\frac{N+\alpha}{N-2}}+\varepsilon^{-\frac{N+\alpha}{4}+\frac{N-2}{4}}g(\varepsilon^{\frac{N-2}{4}}v)v).\\
\end{aligned}$$
The corresponding energy functional of \eqref{1.7} is given by
$$\begin{aligned}
I_{\varepsilon}^*(v)=\frac{1}{2}\int_{\mathbb R^N}|\nabla v|^2+|v|^2
-\frac{1}{2}\int_{\mathbb R^N}&(I_{\alpha}*(|v|^{\frac{N+\alpha}{N-2}}+\varepsilon^{-\frac{N+\alpha}{4}}G(\varepsilon^{\frac{N-2}{4}}v)))\\
\cdot&(|v|^{\frac{N+\alpha}{N-2}}+\varepsilon^{-\frac{N+\alpha}{4}}G(\varepsilon^{\frac{N-2}{4}}v)),\\
\end{aligned}$$
and the limiting energy functional $I_\infty^*:H^1(\mathbb R^N)\to \mathbb R$ is
$$
I_{\infty}^*(v)=\frac{1}{2}\int_{\mathbb R^N}|\nabla v|^2+|v|^2
-\frac{1}{2}\int_{\mathbb R^N}(I_{\alpha}*|v|^{\frac{N+\alpha}{N-2}})|v|^{\frac{N+\alpha}{N-2}}.$$
It is easy to see that
$$m_{\varepsilon}^*:=\inf_{v\in \mathcal{N}_{\varepsilon}^*}I_{\varepsilon}^*(v),\ \ \ m_{\infty}^*:=\inf_{v \in \mathcal{N}_{\infty}^*}I_{\infty}^*(v),$$
are well defined and positive.


Next, we consider the rescaling again
\begin{equation}\label{1.8}
w(x)=\varepsilon^{-\frac{(N-2)\sigma}{4}}v(\varepsilon^{-\frac{\sigma}{2}}x),
\end{equation}
where $\sigma=\frac{(N+\alpha)-(N-2)q_2}{4+(N-2)q_2-(N+\alpha)}>0$, then equation \eqref{1.7} transforms into the equation
\begin{equation}\label{1.9}
-\Delta w+ \varepsilon^{-\sigma}w=\left(I_{\alpha}*(|w|^{\frac{N+\alpha}{N-2}} +\varepsilon_1^{-1}G(\varepsilon_2w))\right)
\left(\frac{N+\alpha}{N-2}|w|^{\frac{N+\alpha}{N-2}-2}w
+\varepsilon_1^{-1}\varepsilon_2g(\varepsilon_2w)\right),
\end{equation}
where
$$
\varepsilon_1=\varepsilon^{\frac{(N+\alpha)(1+\sigma)}{4}},\ \ \ \varepsilon_2=\varepsilon^{\frac{(N-2)(1+\sigma)}{4}}.
$$
Clearly, we have $\lim_{\varepsilon\to\infty}\varepsilon_1=\lim_{\varepsilon\to \infty}\varepsilon_2=\infty$, and  the relation
\begin{equation}\label{A}
\varepsilon^{-\sigma}=\varepsilon_1^{-1}\varepsilon_2^{q_2}.
\end{equation}
The corresponding energy functional of \eqref{1.9} is given by
$$\begin{aligned}
J_{\varepsilon}(w)=\frac{1}{2}\int_{\mathbb R^N}|\nabla w|^2+\varepsilon^{-\sigma}|w|^2-\frac{1}{2}\int_{\mathbb R^N}
&(I_{\alpha}*(|w|^{\frac{N+\alpha}{N-2}}+\varepsilon_1^{-1}G(\varepsilon_2w)))\\
\cdot&(|w|^{\frac{N+\alpha}{N-2}}+\varepsilon_1^{-1}G(\varepsilon_2w)).\\
\end{aligned}$$
Arguing as in \cite{LM2019, LMZ2019}, we know that any ground state solution of \eqref{1.9} satisfies $w_\varepsilon\in \mathcal{N}_{\varepsilon}\cap \mathcal{P}_{\varepsilon}$,   where $\mathcal{N}_{\varepsilon}$ and $\mathcal{P}_{\varepsilon}$ are the Nehari manifold and Poho\v{z}aev manifold  corresponding to \eqref{1.9} and given by
$$\mathcal{N}_{\varepsilon}:=\{w \in H^1(\mathbb R^N)\backslash\{0\}\mid N_{\varepsilon}(w)=0 \},$$
$$\mathcal{P}_{\varepsilon}:=\{w \in H^1(\mathbb R^N)\backslash\{0\}\mid P_{\varepsilon}(w)=0\},$$
 respectively, where
$$\begin{aligned}
N_{\varepsilon}(w)=\int_{\mathbb R^N}|\nabla w|^2+\varepsilon^{-\sigma}|w|^2
-\int_{\mathbb R^N}&(I_{\alpha}*(|w|^{\frac{N+\alpha}{N-2}}+\varepsilon_1^{-1}G(\varepsilon_2w)))\\
\cdot&(\frac{N+\alpha}{N-2}|w|^{\frac{N+\alpha}{N-2}}+\varepsilon_1^{-1}\varepsilon_2g(\varepsilon_2w)w),\\
\end{aligned}$$
$$\begin{aligned}
P_{\varepsilon}(w)&=\frac{N-2}{2}\int_{\mathbb R^N}|\nabla w|^2+\frac{N\varepsilon^{-\sigma}}{2}\int_{\mathbb R^N}|w|^2\\
&\ \ \ -\frac{N+\alpha}{2}\int_{\mathbb R^N}(I_{\alpha}*(|w|^{\frac{N+\alpha}{N-2}}
+\varepsilon_1^{-1}G(\varepsilon_2w)))(|w|^{\frac{N+\alpha}{N-2}}+\varepsilon_1^{-1}G(\varepsilon_2w)).\\
\end{aligned}$$
As $\varepsilon \to \infty$, the limit of the equation \eqref{1.9} is the critical Choquard equation \eqref{1.10},
 the associated energy functional  is given by
$$J_{\infty}(w)=\frac{1}{2}\int_{\mathbb R^N}|\nabla w|^2-\frac{1}{2}\int_{\mathbb R^N}(I_{\alpha}*|w|^{\frac{N+\alpha}{N-2}})|w|^{\frac{N+\alpha}{N-2}},$$
and the corresponding Nehari and Poho\v zaev manifolds are defined by
$$\mathcal{N}_{\infty}=\mathcal{P}_{\infty}:=\left\{w \in \mathcal{D}^{1,2}(\mathbb R^N)\backslash\{0\} \mid
\int_{\mathbb R^N}|\nabla w|^2=\frac{N+\alpha}{N-2}\int_{\mathbb R^N}(I_{\alpha}*|w|^{\frac{N+\alpha}{N-2}})|w|^{\frac{N+\alpha}{N-2}}\right\}.$$
It is known that
$$m_{\infty}:=\inf_{w \in \mathcal{N}_{\infty}} J_{\infty}(w)=\inf_{w \in \mathcal{P}_{\infty}} J_{\infty}(w),$$
is well-defined and positive. Similarly, arguing as in \cite{LM2019, LMZ2019}, it is  easy to know that
$$m_{\varepsilon}:=\inf_{w \in \mathcal{N}_{\varepsilon}} J_{\varepsilon}(w)=\inf_{w \in \mathcal{P}_{\varepsilon}} J_{\varepsilon}(w)$$
is well-defined and positive.

\smallskip
In order to prove our results, we give some lemmas, which are useful for the subsequent proof. First, we give the well known Hardy-Littlewood-Sobolev inequality, which can be found in \cite{LL2001}.
\begin{lemma}\label{lem 2.2}
Let $p$, $r>1$ and $0<\alpha<N$ with $1/p+(N-\alpha)/N+1/r=2$. Let $u \in L^{p}(\mathbb R^N)$ and $v \in L^{r}(\mathbb R^N)$. Then there exists a sharp constant
$C(N,\alpha,p)$, independent of $u$ and $v$, such that
$$\left|\int_{\mathbb R^N}\int_{\mathbb R^N}\frac{u(x)v(y)}{|x-y|^{N-\alpha}}\right|\leq C(N,\alpha,p)\|u\|_{p}\|v\|_{r}.$$
If $p=r=\frac{2N}{N+\alpha}$, then
$$C(N,\alpha,p)=C_{\alpha}(N)=\pi^{\frac{N-\alpha}{2}}
\frac{\Gamma(\frac{\alpha}{2})}{\Gamma(\frac{N+\alpha}{2})}\left\{\frac{\Gamma(\frac{N}{2})}{\Gamma(N)}\right\}^{-\frac{\alpha}{N}}.$$
\end{lemma}

\begin{remark}
By the Hardy-Littlewood-Sobolev inequality, for any $v\in L^{s}(\mathbb R^N)$ with
$s \in (1,\frac{N}{\alpha})$, $I_{\alpha}*v \in L^{\frac{Ns}{N-\alpha s}}(\mathbb R^N)$ and
$$\|I_{\alpha}*v\|_{\frac{Ns}{N-\alpha s}}\leq A_{\alpha}(N)C(N,\alpha,s)\|v\|_s.$$
\end{remark}

Next, we describe the following minimax characterizations for the least energy $m_{\varepsilon}$.
\begin{lemma}\label{lem 2.1}
Let
$$u_{t}(x)=
\left\{
\begin{aligned}
&u \left(\frac{x}{t}\right),
\ & if& \ t>0,\\
&0,
\ & if&\ t=0,\\
\end{aligned}
\right.$$
then
$$m_{\varepsilon}=\inf_{u\in H^{1}(\mathbb R^N)\backslash\{0\}}\sup_{t\geq 0}J_{\varepsilon}(tu)=\inf_{u\in H^{1}(\mathbb R^N)\backslash\{0\}}\sup_{t\geq 0}J_{\varepsilon}(u_t).$$
In particular, we have $m_{\varepsilon}=J_{\varepsilon}(w_{\varepsilon})=\sup_{t>0}J_{\varepsilon}(tw_{\varepsilon})=\sup_{t>0}J_{\varepsilon}((w_{\varepsilon})_t)$.
\end{lemma}
\begin{proof}
The proof is standard, so we omit it here, and we refer readers to \cite{JJ2002,LM2019}.
\end{proof}

\begin{lemma}\label{lem 1.1a}
Assume that  (H1) and (H2) hold, then the rescaled family of solutions $\{w_{\varepsilon}\}$ is bounded in $H^1(\mathbb R^N)$, and satisfies
\begin{equation}\label{abc1}
\|w_{\varepsilon}\|_2^2=\frac{b[(N+\alpha)-(N-2)q_2]}{2}\int_{\mathbb R^N}(I*|w_{\varepsilon}|^{\frac{N+\alpha}{N-2}})|w_{\varepsilon}|^{q_2}+o_{\varepsilon}(1).
\end{equation}
\end{lemma}

\begin{proof}
It is not hard to show that $m_\varepsilon=m_{\varepsilon}^*\leq m_{\infty}^*$. Moreover, according to \eqref{1.5} and \eqref{1.8}, we have
$$
\begin{aligned}
m_{\varepsilon}&=J_{\varepsilon}(w_{\varepsilon})=J_{\varepsilon}(w_{\varepsilon})-\frac{1}{2q_1}{J'_{\varepsilon}}(w_{\varepsilon})w_{\varepsilon}\\
&\geq \left(\frac{1}{2}-\frac{1}{2q_1}\right)\int_{\mathbb R^N}|\nabla w_{\varepsilon}|^2+\varepsilon^{-\sigma}|w_{\varepsilon}|^2\\
&+\left(\frac{N+\alpha}{2q_1(N-2)}-\frac{1}{2}\right)\int_{\mathbb R^N}(I_{\alpha}*(|w_{\varepsilon}|^{\frac{N+\alpha}{N-2}}
+\varepsilon_1^{-1}G(\varepsilon_2w_{\varepsilon})))|w_{\varepsilon}|^{\frac{N+\alpha}{N-2}}\\
&\geq \left(\frac{1}{2}-\frac{1}{2q_1}\right)\int_{\mathbb R^N}|\nabla w_{\varepsilon}|^2+\varepsilon^{-\sigma}|w_{\varepsilon}|^2.\\
\end{aligned}
$$
Hence,  $\{w_{\varepsilon}\}$ is bounded in $\mathcal{D}^{1,2}(\mathbb R^N)$ and $\varepsilon^{-\sigma}\|w_{\varepsilon}\|_2^2$ is bounded.
It suffices to show that $\{w_{\varepsilon}\}$ is also bounded in $L^2(\mathbb R^N)$.
Since $w_{\varepsilon} \in \mathcal{N}_{\varepsilon}\cap \mathcal{P}_{\varepsilon}$, we obtain
\begin{equation}\label{2.1}
\begin{aligned}
\varepsilon^{-\sigma}\|w_{\varepsilon}\|_2^2&=
\int_{\mathbb R^N}\int_{\mathbb R^N}\frac{(|w_{\varepsilon}|^{\frac{N+\alpha}{N-2}}+\varepsilon_1^{-1}G(\varepsilon_2w_{\varepsilon}))
(\frac{N+\alpha}{2}\varepsilon_1^{-1}G(\varepsilon_2w_{\varepsilon}))}{|x-y|^{N-\alpha}}dxdy\\
&\ \ \ -\int_{\mathbb R^N}\int_{\mathbb R^N}\frac{(|w_{\varepsilon}|^{\frac{N+\alpha}{N-2}}+\varepsilon_1^{-1}G(\varepsilon_2w_{\varepsilon}))
(\frac{N-2}{2}\varepsilon_1^{-1}\varepsilon_2g(\varepsilon_2w_{\varepsilon})w_{\varepsilon})}{|x-y|^{N-\alpha}}dxdy.\\
\end{aligned}
\end{equation}
By \eqref{1.5}, we have
$$
\begin{aligned}
&\varepsilon^{-\sigma}\|w_{\varepsilon}\|_2^2\leq
C\int_{\mathbb R^N}\int_{\mathbb R^N}
\frac{|w_{\varepsilon}|^{\frac{N+\alpha}{N-2}}(\varepsilon_1^{-1}\varepsilon_2^{q_1}|w_{\varepsilon}|^{q_1}+\varepsilon_1^{-1}\varepsilon_2^{q_2}|w_{\varepsilon}|^{q_2})}{|x-y|^{N-\alpha}}dxdy\\
+&C\int_{\mathbb R^N}\int_{\mathbb R^N}
\frac{2\varepsilon_1^{-2}\varepsilon_2^{q_1+q_2}|w_{\varepsilon}|^{q_1}|w_{\varepsilon}|^{q_2}+\varepsilon_1^{-2}\varepsilon_2^{2q_1}|w_{\varepsilon}|^{q_1}|w_{\varepsilon}|^{q_1}
+\varepsilon_1^{-2}\varepsilon_2^{2q_2}|w_{\varepsilon}|^{q_2}|w_{\varepsilon}|^{q_2}}{|x-y|^{N-\alpha}}dxdy.\\
\end{aligned}
$$
From the Hardy-Littlewood-Sobolev, the Sobolev and the interpolation inequalities, we obtain
\begin{equation}\label{2.2}
\begin{aligned}
&\int_{\mathbb R^N}\int_{\mathbb R^N}\frac{|w_{\varepsilon}(x)|^{\frac{N+\alpha}{N-2}}|w_{\varepsilon}(y)|^{q_i}}{|x-y|^{N-\alpha}}dxdy\\
&\lesssim \|w_{\varepsilon}\|_2^{\frac{2^*-\tilde{q}_i}{2^*-2}\cdot\frac{N+\alpha}{N}}
(S^{-1}\|\nabla w_{\varepsilon}\|_2^2)^{\frac{2^*}{2}\cdot(1+\frac{\tilde{q}_i-2}{2^*-2})\cdot\frac{N+\alpha}{2N}}\lesssim \|w_{\varepsilon}\|_2^{\frac{2^*-\tilde{q}_i}{2^*-2}\cdot\frac{N+\alpha}{N}},\\
\end{aligned}
\end{equation}
where $\tilde{q_i}=\frac{2Nq_i}{N+\alpha}$. Similarly, we also have
\begin{equation}\label{2.3}
\begin{aligned}
\int_{\mathbb R^N}\int_{\mathbb R^N}\frac{|w_{\varepsilon}(x)|^{q_i}|w_{\varepsilon}(y)|^{q_j}}{|x-y|^{N-\alpha}}dxdy
\lesssim\|w_{\varepsilon}\|_2^{\left(\frac{2^*-\tilde{q}_i}{2^*-2}+\frac{2^*-\tilde{q}_j}{2^*-2}\right)\cdot\frac{N+\alpha}{N}},
\end{aligned}
\end{equation}
According \eqref{2.2}, \eqref{2.3}, note that $\varepsilon^{-\sigma}\|w_\varepsilon\|_2^2=\varepsilon_1^{-1}\varepsilon_2^{q_2}\|w_\varepsilon\|_2^2$ is bounded,  we get
$$\begin{aligned}
\|w_{\varepsilon}\|_2^2&\lesssim \varepsilon_2^{q_i-q_2}\|w_{\varepsilon}\|_2^{\frac{2^*-\tilde{q}_i}{2^*-2}\cdot\frac{N+\alpha}{N}}
+\varepsilon_1^{-1}\varepsilon_2^{q_i+q_j-q_2}\|w_{\varepsilon}\|_2^{\left(\frac{2^*-\tilde{q}_i}{2^*-2}+\frac{2^*-\tilde{q}_j}{2^*-2}\right)\cdot\frac{N+\alpha}{N}}\\
&\lesssim \varepsilon_2^{q_i-q_2}\|w_{\varepsilon}\|_2^{\frac{2^*-\tilde{q}_i}{2^*-2}\cdot\frac{N+\alpha}{N}}
+\varepsilon_2^{q_i+q_j-2q_2}\|w_{\varepsilon}\|_2^{\left(\frac{2^*-\tilde{q}_i}{2^*-2}+\frac{2^*-\tilde{q}_j}{2^*-2}\right)\cdot\frac{N+\alpha}{N}-2}\\
&\lesssim \|w_{\varepsilon}\|_2^{\frac{2^*-\tilde{q}_i}{2^*-2}\cdot\frac{N+\alpha}{N}}
+\|w_{\varepsilon}\|_2^{\left(\frac{2^*-\tilde{q}_i}{2^*-2}+\frac{2^*-\tilde{q}_j}{2^*-2}\right)\cdot\frac{N+\alpha}{N}-2},\\
\end{aligned}$$
by $q_i\in \left(\frac{N+\alpha}{N}, \frac{N+\alpha}{N-2}\right)$, we get that
$$\frac{2^*-\tilde{q}_i}{2^*-2}\cdot\frac{N+\alpha}{N}<2,\quad \left(\frac{2^*-\tilde{q}_i}{2^*-2}+\frac{2^*-\tilde{q}_j}{2^*-2}\right)\cdot\frac{N+\alpha}{N}<4.
$$
Hence, $w_\varepsilon$ is bounded in $L^2(\mathbb R^N)$.

Finally, by virtue of \eqref{2.1},  \eqref{a1} and Lemma \ref{lem 1.1a} and arguing as in \cite[Lemma 4.2]{MM2024}, note that $\varepsilon^{-\sigma}=\varepsilon_1^{-1}\varepsilon_2^{q_2}$,  we have
$$
\begin{aligned}
\varepsilon^{-\sigma}\|w_{\varepsilon}\|_2^2&=
\int_{\mathbb R^N}\int_{\mathbb R^N}\frac{(|w_{\varepsilon}|^{\frac{N+\alpha}{N-2}}+b\varepsilon_1^{-1}|\varepsilon_2w_{\varepsilon}|^{q_2})
\frac{(N+\alpha)b}{2}\varepsilon_1^{-1}|\varepsilon_2w_{\varepsilon}|^{q_2}}{|x-y|^{N-\alpha}}dxdy\\
&-\int_{\mathbb R^N}\int_{\mathbb R^N}\frac{(|w_{\varepsilon}|^{\frac{N+\alpha}{N-2}}+b\varepsilon_1^{-1}|\varepsilon_2w_{\varepsilon}|^{q_2})
\frac{(N-2)b}{2}q_2\varepsilon_1^{-1}|\varepsilon_2w_{\varepsilon}|^{q_2}}{|x-y|^{N-\alpha}}dxdy+o_{\varepsilon}(1).\\
\end{aligned}
$$
Hence, we have
$$\|w_{\varepsilon}\|_2^2=\frac{b[(N+\alpha)-(N-2)q_2]}{2}\int_{\mathbb R^N}(I*|w_{\varepsilon}|^{\frac{N+\alpha}{N-2}})|w_{\varepsilon}|^{q_2}+o_{\varepsilon}(1).$$
The proof is complete. \end{proof}

Finally, we give several classical lemmas, which are all used during the process of our proof.
\begin{lemma}\label{lem 2.4}
(P.L.Lions \cite{L1984})Let $r>0$ and $2\leq q\leq 2^*$. If $(u_n)$ is bounded in $H^1(\mathbb R^N)$ and if
$$\sup_{y \in \mathbb R^N}\int_{B_r(y)}|u_n|^q \to 0,\ \ as\ \ n \to \infty,$$
then $u_n \to 0$ in $L^s({\mathbb R^N})$ for $2<s<2^*$. Moreover, if $q=2^*$, then $u_n\to 0$ in $L^{2^*}(\mathbb R^N)$.
\end{lemma}

\begin{lemma}\label{lem 2.5}
(Radial Lemma A.II, H. Berestycki and P. L. Lions \cite{BL1983}) Let $N\geq2$, then every radial function $u \in H^1(\mathbb R^N)$ is almost everywhere equal to a
function $\tilde{u}(x)$, continuous for $x\neq0$, such that
\begin{equation}\label{2.4}
|\tilde{u}(x)|\leq C_{N}|x|^{-\frac{N-1}{2}}\|u\|_{H^1(\mathbb R^N)}\ \ \ \mbox{for}\ \ \ |x|\geq\alpha_{N},
\end{equation}
where $C_N$ and $\alpha_{N}$ depend only on the dimension $N$.
\end{lemma}

\begin{lemma}\label{lem 2.6}
(Radial Lemma A.III, H. Berestycki and P. L. Lions \cite{BL1983}) Let $N\geq3$, then every radial function $u \in \mathcal{D}^{1,2}(\mathbb R^N)$ is almost everywhere equal to a
function $\tilde{u}(x)$, continuous for $x\neq0$, such that
\begin{equation}\label{2.5}
|\tilde{u}(x)|\leq C_{N}|x|^{-\frac{N-2}{2}}\|u\|_{\mathcal{D}^{1,2}(\mathbb R^N)}\ \ \ \mbox{for}\ \ \ |x|\geq 1,
\end{equation}
where $C_N$ only depends on the dimension $N$.
\end{lemma}


\begin{lemma}\label{lem 2.7}
Let $0<\alpha<N$, $0\leq f \in L^1(\mathbb R^N)$ be a radially symmetric function such that
\begin{equation}\label{2.6}
\lim_{|x|\to \infty}f(|x|)|x|^N=0.
\end{equation}
If $\alpha\leq1$, we additionally assume that $f$ is monotone non-increasing. Then as $|x|\to +\infty$, we have
\begin{equation}\label{2.8}
\int_{\mathbb R^N}\frac{f(y)dy}{|x-y|^{N-\alpha}}=\frac{\|f\|_{L^1}}{|x|^{N-\alpha}}+o\left(\frac{1}{|x|^{N-\alpha}}\right).
\end{equation}
\end{lemma}
\begin{proof}
The proof comes from Lemma 3.10 of \cite{MM}.
\end{proof}

The following Moser iteration lemma is given in Proposition B.1 of \cite{AIIKN2019}. See also \cite{GT1983, LLW2013}.
\begin{lemma}\label{lem 2.8}
Assume $N\geq3$. Let $a(x)$ and $b(x)$ be functions on $B_4$, and let $u\in H^{1}(B_4)$ be a weak solution to
\begin{equation}\label{2.9}
-\Delta u+a(x)u=b(x)u\ \ \ in\ \ \ B_4.
\end{equation}
Suppose that $a(x)$ and $u$ satisfy that
\begin{equation}\label{2.10}
a(x)\geq 0\ \ \ \mbox{for} \ \ \ a.e.\ \ \ x\in B_4,
\end{equation}
and
\begin{equation}\label{2.11}
\int_{B_4}a(x)|u(x)v(x)|dx<\infty \ \ \ \mbox{for}\ \ \ each\ \ \ v\in H_0^1(B_4).
\end{equation}
(i) \ Assume that for any $\kappa \in (0,1)$, there exists $t_{\kappa}>0$ such that
$$\|\chi_{[|b|>t_{\kappa}]}b\|_{L^{N/2}(B_4)}\leq \kappa,$$
where $[|b|>t]:=\{x\in B_4: |b(x)|>t\}$, and $\chi_{A}(x)$ denotes the characteristic function of $A\subset \mathbb R^N$. Then for any $r\in (0,\infty)$, there
exists a constant $C(N,r,t_{\kappa})$ such that
$$\||u|^{r+1}\|_{H^1(B_1)}\leq C(N,r,t_{\kappa})\|u\|_{L^{2^*}(B_4)}.$$
(ii) Let $s>\frac{N}{2}$ and assume that $b\in L^{s}(B_4)$. Then there exists a constant $C(N,s,\|b\|_{L^{s}(B_4)})$ such that
$$\|u\|_{L^{\infty}(B_1)}\leq C(N,s,\|b\|_{L^{s}(B_4)})\|u\|_{L^{2^*}(B_4)}.$$
Here, the constants $C(N,r,t_{\kappa})$ and $C(N,s,\|b\|_{L^{s}(B_4)})$ in $(i)$ and $(ii)$ remain bounded as long as $r$, $t_{\kappa}$ and
$\|b\|_{L^{s}(B_4)}$ are bounded.
\end{lemma}

\section{The proof of main result  for the cases $N\geq5$}
In this section, we prove that the asymptotic behavior of positive ground state solutions of \eqref{1.1} for the cases $N\geq5$.
For $w \in H^1(\mathbb R^N)\backslash\{0\}$, let
\begin{equation}\label{3.1}
\tau(w)=\frac{N-2}{N+\alpha}\cdot \frac{\int_{\mathbb R^N}|\nabla w|^2}{\int_{\mathbb R^N}(I_{\alpha}*|w|^{\frac{N+\alpha}{N-2}})|w|^{\frac{N+\alpha}{N-2}}},
\end{equation}
then $\tau(w)^{\frac{N-2}{2(2+\alpha)}}w \in \mathcal{N}_{\infty}$ for any $w \in H^1(\mathbb R^N)\backslash\{0\}$ and $w \in \mathcal{N}_{\infty}$ if and only if $\tau(w)=1$.

\begin{lemma}\label{lem 1.1b}
Assume that  (H1) and (H2) hold, then there holds
$$m_{\infty}-m_{\varepsilon}\sim \varepsilon^{-\sigma} \ \ as \ \ \varepsilon \ \to \ \ \infty.$$
\end{lemma}
\begin{proof}
First, we claim that there exists a constant $C>0$ such that
\begin{equation}\label{3.2}
1 \leq \tau(w_{\varepsilon})\leq 1+C\varepsilon^{-\sigma}.
\end{equation}
In fact, since $w_{\varepsilon} \in \mathcal{N}_{\varepsilon}$, then it follows from \eqref{1.5} that
$$\begin{aligned}
&\tau(w_{\varepsilon})=\frac{N-2}{N+\alpha}\cdot \frac{\int_{\mathbb R^N}|\nabla w_{\varepsilon}|^2}{\int_{\mathbb R^N}(I_{\alpha}*|w_{\varepsilon}|^{\frac{N+\alpha}{N-2}})|w_{\varepsilon}|^{\frac{N+\alpha}{N-2}}}\\
=&1+\frac{N-2}{N+\alpha}\frac{\int_{\mathbb R^N}(I_{\alpha}*\varepsilon_1^{-1}G(\varepsilon_2w_{\varepsilon}))\varepsilon_1^{-1}\varepsilon_2g(\varepsilon_2w_{\varepsilon})w_{\varepsilon}}
{\int_{\mathbb R^N}(I_{\alpha}*|w_{\varepsilon}|^{\frac{N+\alpha}{N-2}})|w_{\varepsilon}|^{\frac{N+\alpha}{N-2}}}\\
+&\frac{N-2}{N+\alpha}\frac{\int_{\mathbb R^N}(I_{\alpha}*|w_{\varepsilon}|^{\frac{N+\alpha}{N-2}})(\frac{N+\alpha}{N-2}\varepsilon_1^{-1}G(\varepsilon_2w_{\varepsilon})
+\varepsilon_1^{-1}\varepsilon_2g(\varepsilon_2w_{\varepsilon})w_{\varepsilon})-\varepsilon^{-\sigma}\int_{\mathbb R^N}|w_{\varepsilon}|^2}
{\int_{\mathbb R^N}(I_{\alpha}*|w_{\varepsilon}|^{\frac{N+\alpha}{N-2}})|w_{\varepsilon}|^{\frac{N+\alpha}{N-2}}}\\
\lesssim &1+C\varepsilon^{-\sigma}\frac{\int_{\mathbb R^N}(I_{\alpha}*|w_{\varepsilon}|^{\frac{N+\alpha}{N-2}})|w_{\varepsilon}|^{q_2}-\int_{\mathbb R^N}|w_{\varepsilon}|^2}
{\int_{\mathbb R^N}(I_{\alpha}*|w_{\varepsilon}|^{\frac{N+\alpha}{N-2}})|w_{\varepsilon}|^{\frac{N+\alpha}{N-2}}}.
\end{aligned}$$
By the Hardy-Littlewood-Sobolev inequality, we get that
$$\begin{aligned}
&\frac{\int_{\mathbb R^N}(I_{\alpha}*|w_{\varepsilon}|^{\frac{N+\alpha}{N-2}})|w_{\varepsilon}|^{q_2}-\int_{\mathbb R^N}|w_{\varepsilon}|^2}{\int_{\mathbb R^N}|w_{\varepsilon}|^{2^*}}\\
\leq &\frac{\left(\int_{\mathbb R^N}|w_{\varepsilon}|^2\right)^{\frac{2^*-\tilde{q}_2}{2^*-2}\cdot\frac{N+\alpha}{2N}}
\left(\int_{\mathbb R^N}|w_{\varepsilon}|^{2^*}\right)^{(1+\frac{\tilde{q}_2-2}{2^*-2})
\cdot\frac{N+\alpha}{2N}}-\int_{\mathbb R^N}|w_{\varepsilon}|^2}{\int_{\mathbb R^N}|w_{\varepsilon}|^{2^*}}\\
=&\left(\frac{\int_{\mathbb R^N}|w_{\varepsilon}|^{2}}{\int_{\mathbb R^N}|w_{\varepsilon}|^{2^*}}\right)^{\frac{2^*-\tilde{q}_2}{2^*-2}\cdot\frac{N+\alpha}{2N}}
(\int_{\mathbb R^N}|w_{\varepsilon}|^{2^*})^{\frac{\alpha}{N}}-\left(\frac{\int_{\mathbb R^N}|w_{\varepsilon}|^{2}}{\int_{\mathbb R^N}|w_{\varepsilon}|^{2^*}}\right)\\
\leq &C\eta_{\varepsilon}^{\frac{2^*-\tilde{q}_2}{2^*-2}\cdot\frac{N+\alpha}{2N}}-\eta_{\varepsilon}=G(\eta_{\varepsilon}),
\end{aligned}$$
where $\eta_{\varepsilon}:=\frac{\int_{\mathbb R^N}|w_{\varepsilon}|^{2}}{\int_{\mathbb R^N}|w_{\varepsilon}|^{2^*}}$.
According to $\frac{2^*-\tilde{q}_2}{2^*-2}\cdot\frac{N+\alpha}{2N}<1$, we have $$G(\bar{\eta}_{\varepsilon})=\max G(\eta_{\varepsilon})\leq C.$$
Therefore, by the boundedness of $w_{\varepsilon}$ in $\mathcal{D}^{1,2}(\mathbb R^N)$, we obtain
$$\begin{aligned}
\tau(w_{\varepsilon})&\lesssim 1+C\varepsilon^{-\sigma} G(\bar{\eta}_{\varepsilon})\frac{\int_{\mathbb R^N}|w_{\varepsilon}|^{2^*}}
{\int_{\mathbb R^N}(I_{\alpha}*|w_{\varepsilon}|^{\frac{N+\alpha}{N-2}})|w_{\varepsilon}|^{\frac{N+\alpha}{N-2}}}\\
&\lesssim 1+C\varepsilon^{-\sigma} G(\bar{\eta}_{\varepsilon})\frac{S^{-\frac{N}{N-2}}(\int_{\mathbb R^N}|\nabla w_{\varepsilon}|^{2}dx)^{\frac{N}{N-2}}}
{\int_{\mathbb R^N}(I_{\alpha}*|w_{\varepsilon}|^{\frac{N+\alpha}{N-2}})|w_{\varepsilon}|^{\frac{N+\alpha}{N-2}}}\\
&= 1+C\varepsilon^{-\sigma} G(\bar{\eta}_{\varepsilon})S^{-\frac{N}{N-2}}\tau(w_{\varepsilon})(\int_{\mathbb R^N}|\nabla w_{\varepsilon}|^{2}dx)^{\frac{2}{N-2}}\\
&\leq 1+C\varepsilon^{-\sigma}\tau(w_{\varepsilon}),
\end{aligned}$$
and hence for large $\varepsilon>0$, there holds
$$\tau(w_{\varepsilon})\leq \frac{1}{1-C\varepsilon^{-\sigma}}=\frac{1-C\varepsilon^{-\sigma}+C\varepsilon^{-\sigma}}{1-C\varepsilon^{-\sigma}}\leq 1+C\varepsilon^{-\sigma}.$$

On the other hand, by \eqref{1.5} and \eqref{2.1}, we have that
$$\int_{\mathbb R^N}(I_{\alpha}*|w_{\varepsilon}|^{\frac{N+\alpha}{N-2}})(\frac{N+\alpha}{N-2}\varepsilon_1^{-1}G(\varepsilon_2w_{\varepsilon})
+\varepsilon_1^{-1}\varepsilon_2g(\varepsilon_2w_{\varepsilon})w_{\varepsilon})-\varepsilon^{-\sigma}\int_{\mathbb R^N}|w_{\varepsilon}|^2>0.$$
Therefore $\tau(w_{\varepsilon})>1$. The claim \eqref{3.2} is proved.

Second, if $N\geq3$, then by \eqref{1.5}, Lemma \ref{lem 2.1} and the boundedness of $\{w_{\varepsilon}\}$, we find that
$$\begin{aligned}
m_{\infty}&\leq \sup_{t\geq0}J_{\infty}((w_{\varepsilon})_t)=J_{\infty}((w_{\varepsilon})_{t_{\varepsilon}})\\
&\leq \sup_{t\geq0}J_{\varepsilon}((w_{\varepsilon})_t)-\frac{1}{2}\varepsilon^{-\sigma}\int_{\mathbb R^N}|(w_{\varepsilon})_{t_{\varepsilon}}|^2dx\\
&\ \ \ +\frac{1}{2}\int_{\mathbb R^N}(I_{\alpha}*(2|(w_{\varepsilon})_{t_{\varepsilon}}|^{\frac{N+\alpha}{N-2}}+\varepsilon_1^{-1}G(\varepsilon_2(w_{\varepsilon})_{t_{\varepsilon}})))
\varepsilon_1^{-1}G(\varepsilon_2(w_{\varepsilon})_{t_{\varepsilon}})\\
&\leq m_{\varepsilon}+o(\varepsilon^{-\sigma})\\
&\ \ \ \ \ \ \ \ +\varepsilon^{-\sigma}\Big(\tau(w_{\varepsilon})^{\frac{N+\alpha}{2+\alpha}}C\int_{\mathbb R^N}(I_{\alpha}*|w_{\varepsilon}|^{\frac{N+\alpha}{N-2}})|w_{\varepsilon}|^{q_2}
-\frac{\tau(w_{\varepsilon})^{\frac{N}{2+\alpha}}}{2}\int_{\mathbb R^N}|w_{\varepsilon}|^2\Big)\\
&\leq  m_{\varepsilon}+C\varepsilon^{-\sigma},
\end{aligned}$$
where, we have used  the fact that
\begin{equation}\label{3.3}
t_{\varepsilon}=\left(\frac{N-2}{N+\alpha}\cdot \frac{\int_{\mathbb R^N}|\nabla w_{\varepsilon}|^2}{\int_{\mathbb R^N}(I_{\alpha}*|w_{\varepsilon}|^{\frac{N+\alpha}{N-2}})|w_{\varepsilon}|^{\frac{N+\alpha}{N-2}}}\right)^{\frac{1}{2+\alpha}}=\tau(w_{\varepsilon})^{\frac{1}{2+\alpha}}.
\end{equation}

Third, for $N\geq5$ and each $\rho>0$, the family $W_{\rho}:=\rho^{-\frac{N-2}{2}}W_1(x/\rho)$ are radial ground states of
$-\Delta w=\frac{N+\alpha}{N-2}(I_{\alpha}*|w|^{\frac{N+\alpha}{N-2}})|w|^{\frac{N+\alpha}{N-2}-2}w$ and satisfy that
$$\|W_{\rho}\|_2^2=\rho^2\|W_{1}\|_2^2,$$
$$\int_{\mathbb R^N}(I_{\alpha}*|W_\rho|^{\frac{N+\alpha}{N-2}})|W_\rho|^{q_2}
=\rho^{\frac{N+\alpha}{2}-\frac{(N-2)q_2}{2}}\int_{\mathbb R^N}(I_{\alpha}*|W_1|^{\frac{N+\alpha}{N-2}})|W_1|^{q_2}.$$
Let
$$\begin{aligned}
g_0(\rho)=&b\int_{\mathbb R^N}(I_{\alpha}*|W_\rho|^{\frac{N+\alpha}{N-2}})|W_\rho|^{q_2}-\frac{1}{2}\int_{\mathbb R^N}|W_\rho|^2\\
=&b\rho^{\frac{N+\alpha}{2}-\frac{(N-2)q_2}{2}}\int_{\mathbb R^N}(I_{\alpha}*|W_1|^{\frac{N+\alpha}{N-2}})|W_1|^{q_2}-\frac{1}{2}\rho^2\int_{\mathbb R^N}|W_1|^2.\\
\end{aligned}$$
Then there exists $\rho_0=\rho(q_2)\in (0,\infty)$ with
$$\rho_0=\left(\frac{b[(N+\alpha)-(N-2)q_2]}{2}
\frac{\int_{\mathbb R^N}(I_{\alpha}*|W_1|^{\frac{N+\alpha}{N-2}})|W_1|^{q_2}}{\int_{\mathbb R^N}|W_1|^2}\right)^{\frac{2}{4+(N-2)q_2-(N+\alpha)}},$$
such that
$$\begin{aligned}
g_0(\rho_0)&=\sup_{\rho>0}g_0(\rho)\\
&=\frac{b[4+(N-2)q_2-(N+\alpha)]}{4}\cdot\rho_0^{\frac{(N+\alpha)-(N-2)q_2}{2}}
\int_{\mathbb R^N}(I_{\alpha}*|W_1|^{\frac{N+\alpha}{N-2}})|W_1|^{q_2}.\\
\end{aligned}$$

Let $W_0=W_{\rho_0}$,
then by \eqref{a1}, \eqref{1.5} and Lemma \ref{lem 2.1}, there exists $t_{\varepsilon} \in (0,\infty)$ such that
$$\begin{aligned}
m_{\varepsilon}&\leq \sup_{t\geq0}J_{\varepsilon}(tW_0)=J_{\varepsilon}(t_{\varepsilon}W_0)\\
&\leq\sup_{t\geq 0}\left(\frac{t^2}{2}\int_{\mathbb R^N}|\nabla W_0|^2
-\frac{t^{\frac{2(N+\alpha)}{N-2}}}{2}\int_{\mathbb R^N}(I_{\alpha}*|W_0|^{\frac{N+\alpha}{N-2}})|W_0|^{\frac{N+\alpha}{N-2}}\right)\\
&\ -\varepsilon^{-\sigma}\left(bt_{\varepsilon}^{\frac{N+\alpha}{N-2}+q_2}
\int_{\mathbb R^N}(I_{\alpha}*|W_0|^{\frac{N+\alpha}{N-2}})|W_0|^{q_2}-\frac{t_{\varepsilon}^2}{2}\int_{\mathbb R^N}|W_0|^2\right)+o(\varepsilon^{-2\sigma})\\
&=m_{\infty}-\varepsilon^{-\sigma}\left(bt_{\varepsilon}^{\frac{N+\alpha}{N-2}+q_2}
\int_{\mathbb R^N}(I_{\alpha}*|W_0|^{\frac{N+\alpha}{N-2}})|W_0|^{q_2}-\frac{t_{\varepsilon}^2}{2}\int_{\mathbb R^N}|W_0|^2\right)+o(\varepsilon^{-2\sigma}).\\
\end{aligned}$$
According to $t_{\varepsilon}W_0 \in \mathcal{N}_{\varepsilon}$, we have
$$\begin{aligned}
t_{\varepsilon}^2&\int_{\mathbb R^N}|\nabla W_0|^2+\varepsilon^{-\sigma}|W_0|^2
=t_{\varepsilon}^{\frac{2(N+\alpha)}{N-2}}\cdot\frac{N+\alpha}{N-2}\int_{\mathbb R^N}(I_{\alpha}*|W_0|^{\frac{N+\alpha}{N-2}})|W_0|^{\frac{N+\alpha}{N-2}}\\
&+\int_{\mathbb R^N}(I_{\alpha}*|t_{\varepsilon}W_0|^{\frac{N+\alpha}{N-2}})\Big(\frac{N+\alpha}{N-2}\varepsilon_1^{-1}G(\varepsilon_2t_{\varepsilon}W_0)
+\varepsilon_1^{-1}\varepsilon_2g(\varepsilon_2t_{\varepsilon}W_0)t_{\varepsilon}W_0\Big)\\
&+\int_{\mathbb R^N}(I_{\alpha}*\varepsilon_1^{-1}G(\varepsilon_2t_{\varepsilon}W_0))\varepsilon_1^{-1}\varepsilon_2g(\varepsilon_2t_{\varepsilon}W_0)t_{\varepsilon}W_0.\\
\end{aligned}$$
If $t_{\varepsilon}\geq1$, by \eqref{1.5}, we have
$$\begin{aligned}
\int_{\mathbb R^N}|\nabla W_0|^2+\varepsilon^{-\sigma}|W_0|^2
\geq t_{\varepsilon}^{2q_2-2}\Big\{&\frac{N+\alpha}{N-2}\int_{\mathbb R^N}(I_{\alpha}*|W_0|^{\frac{N+\alpha}{N-2}})|W_0|^{\frac{N+\alpha}{N-2}}\\
&+\varepsilon^{-\sigma} c\int_{\mathbb R^N}(I_{\alpha}*|W_0|^{\frac{N+\alpha}{N-2}})|W_0|^{q_2}+o(\varepsilon^{-\sigma})\Big\}.\\
\end{aligned}$$
Hence,
$$
t_{\varepsilon}\leq\max\Bigg\{1,\left(\frac{\frac{2(N+\alpha)m_{\infty}}{2+\alpha}
+\varepsilon^{-\sigma}\int_{\mathbb R^N}|W_0|^2}{\frac{2(N+\alpha)m_{\infty}}{2+\alpha}+\varepsilon^{-\sigma} c\int_{\mathbb R^N}
(I_{\alpha}*|W_0|^{\frac{N+\alpha}{N-2}})|W_0|^{q_2}+o(\varepsilon^{-\sigma})}\right)^{\frac{1}{2q_2-2}}\Bigg\}.
$$
If $t_{\varepsilon}\leq1$, then
$$\begin{aligned}
\int_{\mathbb R^N}|\nabla W_0|^2+\varepsilon^{-\sigma}|W_0|^2
\leq t_{\varepsilon}^{2q_2-2}\Bigg\{&\frac{N+\alpha}{N-2}\int_{\mathbb R^N}(I_{\alpha}*|W_0|^{\frac{N+\alpha}{N-2}})|W_0|^{\frac{N+\alpha}{N-2}}\\
&+\varepsilon^{-\sigma} C\int_{\mathbb R^N}(I_{\alpha}*|W_0|^{\frac{N+\alpha}{N-2}})|W_0|^{q_2}+o(\varepsilon^{-\sigma})\Bigg\}.\\
\end{aligned}$$
Hence,
$$
t_{\varepsilon}\geq\min\Bigg\{1,\left(\frac{\frac{2(N+\alpha)m_{\infty}}{2+\alpha}
+\varepsilon^{-\sigma}\int_{\mathbb R^N}|W_0|^2}{\frac{2(N+\alpha)m_{\infty}}{2+\alpha}+\varepsilon^{-\sigma} C\int_{\mathbb R^N}
(I_{\alpha}*|W_0|^{\frac{N+\alpha}{N-2}})|W_0|^{q_2}+o(\varepsilon^{-\sigma})}\right)^{\frac{1}{2q_2-2}}\Bigg\}.
$$
So we conclude that $t_{\varepsilon}\to 1$ as $\varepsilon \to \infty$.
Since
$$b\int_{\mathbb R^N}(I_{\alpha}*|W_\rho|^{\frac{N+\alpha}{N-2}})|W_\rho|^{q_2}>\frac{1}{2}\int_{\mathbb R^N}|W_\rho|^2,$$
there exists a constant $C>0$ such that
$$m_{\varepsilon}\leq m_{\infty}-C\varepsilon^{-\sigma}$$ for large $\varepsilon>0$,
and the conclusion follows.
\end{proof}

\begin{lemma}\label{lem 1.1c}
Assume that  (H1) and (H2) hold, then $\|w_{\varepsilon}\|_2\sim\|w_\varepsilon\|_{\frac{2Nq_2}{N+\alpha}}\sim1$ as $\varepsilon\to \infty$.
\end{lemma}
\begin{proof}
Combining \eqref{1.5}, \eqref{2.1}, \eqref{3.2} with \eqref{3.3}, we have
$$\begin{aligned}
m_{\infty}&\leq \sup_{t\geq0}J_{\infty}((w_{\varepsilon})_t)\\
&\leq m_{\varepsilon}+\left(1-\frac{(N+\alpha)-(N-2)q_2}{4}\right)
\int_{\mathbb R^N}\int_{\mathbb R^N}\frac{{|w_{\varepsilon}}|^{\frac{N+\alpha}{N-2}}\cdot\varepsilon_1^{-1}G(\varepsilon_2w_{\varepsilon})}{|x-y|^{N-\alpha}}dxdy\\
& \quad + \frac{1}{2}\int_{\mathbb R^N}\int_{\mathbb R^N}\frac{\varepsilon_1^{-1}G(\varepsilon_2(w_{\varepsilon})_{t_{\varepsilon}})
\cdot \varepsilon_1^{-1}G(\varepsilon_2(w_{\varepsilon})_{t_{\varepsilon}})}{|x-y|^{N-\alpha}}dxdy\\
&\leq m_{\varepsilon}+\varepsilon^{-\sigma} C\left(1-\frac{(N+\alpha)-(N-2)q_2}{4}\right)\int_{\mathbb R^N}(I_{\alpha}*|w_{\varepsilon}|^{\frac{N+\alpha}{N-2}})
|w_{\varepsilon}|^{q_2}+o(\varepsilon^{-\sigma}),\\
\end{aligned}$$
and it follows that
$$\int_{\mathbb R^N}(I_{\alpha}*|w_{\varepsilon}|^{\frac{N+\alpha}{N-2}})|w_{\varepsilon}|^{q_2}
\geq \frac{4C(m_{\infty}-m_{\varepsilon})\varepsilon^{\sigma}}{4-(N+\alpha)+q_2(N-2)}\geq C>0.
$$
On the other hand, by \eqref{2.2} and the boundedness of $\{w_{\varepsilon}\}$, we have
$\int_{\mathbb R^N}(I_{\alpha}*|w_{\varepsilon}|^{\frac{N+\alpha}{N-2}})|w_{\varepsilon}|^{q_2}\leq C$.
Therefore, $\int_{\mathbb R^N}(I_{\alpha}*|w_{\varepsilon}|^{\frac{N+\alpha}{N-2}})|w_{\varepsilon}|^{q_2}\sim 1$.
Finally, it follows from \eqref{1.5}, \eqref{abc1}, Lemma  \ref{lem 2.2} and the boundedness of $w_\varepsilon$ in $H^1(\mathbb R^N)$ that
$$\|w_{\varepsilon}\|_2^2\sim1\ \ as\ \ \varepsilon\to \infty,$$
and
$$
1\sim \int_{\mathbb R^N}(I_\alpha\ast |w_\varepsilon|^{\frac{N+\alpha}{N-2}})|w_\varepsilon|^{q_2}\le C\|w_\varepsilon\|_{\frac{2Nq_2}{N+\alpha}}^{q_2}\le C<\infty.
$$
The proof is complete.
\end{proof}

\begin{lemma}\label{lem 1.1d} Assume that  $N\ge 5$, (H1) and  (H2) hold,
then there exists $\xi_{\varepsilon} \in (0,\infty)$ satisfying
$$\xi_{\varepsilon}\thicksim \varepsilon^{-\frac{(N+\alpha)-(N-2)q_2}{2[4+(N-2)q_2-(N+\alpha)]}}$$
such that the rescaled ground states
$$w_{\varepsilon}(x)=\xi_{\varepsilon}^{\frac{N-2}{2}}v_{\varepsilon}(\xi_{\varepsilon}x)$$
converge to $W_{\rho_0}$ in $H^1(\mathbb R^N)$ as
$\varepsilon\to \infty$, where $W_{\rho_0}$ is a positive ground state of the equation \eqref{1.10} with
$$\rho_0=\left(\frac{b[(N+\alpha)-(N-2)q_2]}{2}
\frac{\int_{\mathbb R^N}(I_{\alpha}*|W_1|^{\frac{N+\alpha}{N-2}})|W_1|^{q_2}}{\int_{\mathbb R^N}|W_1|^2}\right)^{\frac{2}{4+(N-2)q_2-(N+\alpha)}}.$$
\end{lemma}
\begin{proof}
The proof is similar to \cite[Lemma 5.7]{MM}, and for the readers' convenience we give the details. Note that $w_{\varepsilon}$ is a positive radially symmetric function and by Lemma \ref{lem 1.1a}, we know that $\{w_{\varepsilon}\}$ is bounded in $H^1(\mathbb R^N)$.
Then there exists $w_{\infty} \in H^1(\mathbb R^N)$ satisfying $-\Delta w=\frac{N+\alpha}{N-2}(I_{\alpha}*|w|^{\frac{N+\alpha}{N-2}})|w|^{\frac{N+\alpha}{N-2}-2}w$ such that
$$w_{\varepsilon}\rightharpoonup w_{\infty}\ \ \mbox{in}\ H^1(\mathbb R^N),\ \ \ w_{\varepsilon}\rightarrow w_{\infty}\ \ \mbox{in}\ \ L^{q}(\mathbb R^N)\ \  \mbox{for} \ \ \ \forall\ q\in(2,2^*),$$
and
$$w_{\varepsilon}\rightarrow w_{\infty}\ \ \mbox{in}\ \ L_{loc}^2(\mathbb R^N),\ \ \ w_{\varepsilon}(x)\rightarrow w_{\infty}(x)\ \ \mbox{a.e.  on} \ \ \mathbb R^N.$$
Observe that
$$\begin{aligned}
J_{\infty}(w_{\varepsilon})&=J_{\varepsilon}(w_{\varepsilon})+\Bigg(\int_{\mathbb R^N}\int_{\mathbb R^N}\frac{|w_{\varepsilon}|^{\frac{N+\alpha}{N-2}}\cdot \varepsilon_1^{-1}G(\varepsilon_2w_{\varepsilon})}{|x-y|^{N-\alpha}}dxdy
-\frac{1}{2}\varepsilon^{-\sigma}\int_{\mathbb R^N}|w_{\varepsilon}|^2dx\Bigg)\\
&\ \ \ +\frac{1}{2}\int_{\mathbb R^N}\int_{\mathbb R^N}\frac{\varepsilon_1^{-1}G(\varepsilon_2w_{\varepsilon})
\cdot \varepsilon_1^{-1}G(\varepsilon_2w_{\varepsilon})}{|x-y|^{N-\alpha}}dxdy\\
&=m_{\varepsilon}+o(1)=m_{\infty}+o(1),\\
\end{aligned}$$
and
$$\begin{aligned}
J'_{\infty}(w_{\varepsilon})w_{\varepsilon}&=J'_{\varepsilon}(w_{\varepsilon})w_{\varepsilon}+\int_{\mathbb R^N}\int_{\mathbb R^N}
\frac{\varepsilon_1^{-1}G(\varepsilon_2w_{\varepsilon})\varepsilon_1^{-1}\varepsilon_2g(\varepsilon_2w_{\varepsilon})w_{\varepsilon}}{|x-y|^{N-\alpha}}dxdy
-\varepsilon^{-\sigma}\int_{\mathbb R^N}|w_{\varepsilon}|^2dx\\
&\ \ \ +\int_{\mathbb R^N}\int_{\mathbb R^N}\frac{|w_{\varepsilon}|^{\frac{N+\alpha}{N-2}}(\frac{N+\alpha}{N-2}\varepsilon_1^{-1}G(\varepsilon_2w_{\varepsilon})
+\varepsilon_1^{-1}\varepsilon_2g(\varepsilon_2w_{\varepsilon})w_{\varepsilon})}{|x-y|^{N-\alpha}}dxdy=o(1).\\
\end{aligned}$$
By Lion's Lemma, it is standard to show that there exist $\zeta^{(j)}_{\varepsilon}\in (0,\infty)$ and $w^{(j)}\in \mathcal{D}^{1,2}(\mathbb R^N)$ ($j=1,2,\ldots,k$) such that
$$w_{\varepsilon}=w_{\infty}+\sum_{j=1}^{k}(\zeta_{\varepsilon}^{(j)})^{-\frac{N-2}{2}}w^{(j)}((\zeta_{\varepsilon}^{(j)})^{-1}x)+\tilde{w}_{\varepsilon},$$
where $\tilde{w}_{\varepsilon}\ \to\ 0$ in $L^{2^*}(\mathbb R^N)$ and $w^{(j)}$ are nontrivial solutions of the limit equation
$-\Delta w=\frac{N+\alpha}{N-2}(I_{\alpha}*|w|^{\frac{N+\alpha}{N-2}})|w|^{\frac{N+\alpha}{N-2}-2}w.$
Besides, we have
$$\limsup_{\varepsilon \to \infty}\|w_{\varepsilon}\|_{\mathcal{D}^{1,2}(\mathbb R^N)}^2\geq \|w_{\infty}\|_{\mathcal{D}^{1,2}(\mathbb R^N)}^2
+\sum_{j=1}^{k}\|w^{(j)}\|_{\mathcal{D}^{1,2}(\mathbb R^N)}^2,$$
and
$$m_{\infty}=J_{\infty}(w_{\infty})+\sum_{j=1}^{k}J_{\infty}(w^{(j)}).$$
Moreover, $J_{\infty}(w_{\infty})\geq 0$ and $J_{\infty}(w^{(j)})\geq m_{\infty}$ for all $j=1,2,\ldots,k$.

If $N\geq5$, then by Lemma \ref{lem 1.1c}, we have $\|w_\infty\|_{\frac{2Nq_2}{N+\alpha}}\not=0$. That is, $w_{\infty}\neq0$ and hence  $J_{\infty}(w_{\infty})=m_{\infty}$ and $k=0$. Thus, $w_{\varepsilon}\to w_{\infty}$ in $L^{2^*}(\mathbb R^N)$.
Since $J'_{\infty}(w_{\varepsilon})\to 0$, it follows that $w_{\varepsilon}\to w_{\infty}$ in $\mathcal{D}^{1,2}(\mathbb R^N)$. Also because $w_{\varepsilon}(x)$ is radial and radially decreasing for
every $x \in \mathbb R^N\backslash\{0\}$, we have
$$w_{\varepsilon}^2(x)\leq \frac{1}{|B_{|x|}|}\int_{B_{|x|}}|w_{\varepsilon}|^2\leq \frac{1}{|x|^{N}}\int_{\mathbb R^N}|w_{\varepsilon}|^2\leq \frac{C}{|x|^{N}},$$
then
\begin{equation}\label{3.4}
w_{\varepsilon}(x)\leq C|x|^{-\frac{N}{2}},\ \ |x|\geq1.
\end{equation}

If $\alpha>N-4$, then we have $\frac{N+\alpha}{N-2}>2$ and hence
$$|w_{\varepsilon}|^{\frac{N+\alpha}{N-2}}|x|^{N}\leq C|x|^{-\frac{N}{2}\cdot\frac{N+\alpha}{N-2}+N}=C|x|^{-\frac{N}{2}\left(\frac{N+\alpha}{N-2}-2\right)}\ \to\ 0
\ \mbox{as}\ |x|\ \to\ \infty.$$
By virtue of Lemma \ref{lem 2.7}, we obtain
\begin{equation}\label{3.5}
\left(I_{\alpha}*|w_{\varepsilon}|^{\frac{N+\alpha}{N-2}}\right)(x)=\frac{\left\| |w_{\varepsilon}|^{\frac{N+\alpha}{N-2}}\right\|_{L^1}}{|x|^{N-\alpha}}
+o\left(\frac{1}{|x|^{N-\alpha}}\right)\leq C|x|^{-N+\alpha},\ \ |x|\geq1,
\end{equation}
and then
\begin{equation}\label{3.6}
\frac{N+\alpha}{N-2}\left(I_{\alpha}*|w_{\varepsilon}|^{\frac{N+\alpha}{N-2}}\right)(x)|w_{\varepsilon}|^{\frac{N+\alpha}{N-2}-2}(x)
\leq C|x|^{-\frac{N^2-N\alpha+4\alpha}{2(N-2)}},\ \ |x|\geq \tilde{R}.\\
\end{equation}
Since
$$\begin{aligned}
&\Big(-\Delta-C|x|^{-\frac{N^2-N\alpha+4\alpha}{2(N-2)}}\Big)w_{\varepsilon}\\
\leq &\Big(-\Delta+\varepsilon^{-\sigma}-\frac{N+\alpha}{N-2}(I_{\alpha}*|w_{\varepsilon}|^{\frac{N+\alpha}{N-2}})|w_{\varepsilon}|^{\frac{N+\alpha}{N-2}-2}\\
&\ -(I_{\alpha}*|w_{\varepsilon}|^{\frac{N+\alpha}{N-2}})\varepsilon_1^{-1}\varepsilon_2g(\varepsilon_2w_{\varepsilon})w_{\varepsilon}^{-1}
-(I_{\alpha}*\varepsilon_1^{-1}G(\varepsilon_2w_{\varepsilon}))|w_{\varepsilon}|^{\frac{N+\alpha}{N-2}-2}\\
&\ -(I_{\alpha}*\varepsilon_1^{-1}G(\varepsilon_2w_{\varepsilon}))\varepsilon_1^{-1}\varepsilon_2g(\varepsilon_2w_{\varepsilon})w_{\varepsilon}^{-1}\Big)w_{\varepsilon}=0,\\
\end{aligned}$$
for large $|x|$. Also,
$$\begin{aligned}
&\left(-\Delta-C|x|^{-\frac{N^2-N\alpha+4\alpha}{2(N-2)}}\right)\frac{1}{|x|^{N-2-\kappa_0}}\\
=&\left(\varepsilon_0(N-2-\varepsilon_0)-C|x|^{-\frac{(N-4)(N-\alpha)+8}{2(N-2)}}\right)\frac{1}{|x|^{N-\kappa_0}}>0,\\
\end{aligned}$$
which is positive for $|x|$ large enough. By the maximum principle on $\mathbb R^N\backslash B_{R}$, we deduce that
\begin{equation}\label{3.7}
w_{\varepsilon}(x)\leq \frac{w_{\varepsilon}(R)R^{N-2-\kappa_0}}{|x|^{N-2-\kappa_0}}\leq \frac{CR^{N-2-\kappa_0}}{|x|^{N-2-\kappa_0}},\ \ \mbox{for} \ \ |x|\geq R.
\end{equation}
When $\kappa_0>0$ is small enough, the domination is in $L^2(\mathbb R^N)$ for $N\geq5$, and this shows,
by the dominated convergence theorem, that $w_{\varepsilon}\to w_{\infty}$ in $L^2(\mathbb R^N)$. Thus, we conclude that $w_{\varepsilon}\to w_{\infty}$ in $H^1(\mathbb R^N)$.
Moreover, by \eqref{abc1}, we obtain
$$\|w_{\infty}\|_2^2=\frac{b[(N+\alpha)-(N-2)q_2]}{2}\int_{\mathbb R^N}(I*|w_{\infty}|^{\frac{N+\alpha}{N-2}})|w_{\infty}|^{q_2},$$
and it follows that $w_{\infty}=W_{\rho_0}$ with
$$\rho_0=\left(\frac{b[(N+\alpha)-(N-2)q_2]}{2}
\frac{\int_{\mathbb R^N}(I_{\alpha}*|W_1|^{\frac{N+\alpha}{N-2}})|W_1|^{q_2}}{\int_{\mathbb R^N}|W_1|^2}\right)^{\frac{2}{4+(N-2)q_2-(N+\alpha)}}.$$
Therefore, if $\alpha>N-4$, then the statement is valid with $\xi_{\varepsilon}=\varepsilon^{-\frac{(N+\alpha)-(N-2)q_2}{2[4+(N-2)q_2-(N+\alpha)]}}$.

If $\alpha\leq N-4$, then for any $\varepsilon_{n} \to \infty$, up to a subsequence, we can assume that $w_{\varepsilon_n}\to W_{\rho}$ in $\mathcal{D}^{1,2}(\mathbb R^N)$
with $\rho \in (0,\rho_0]$. Moreover, $w_{\varepsilon_n}\to W_{\rho}$ in $L^2(\mathbb R^N)$ if and only if $\rho=\rho_0$.
Let $M_{\varepsilon_n}=w_{\varepsilon_n}(0)$ and $z_{\varepsilon_n}=M_{\varepsilon_n}[W_{\rho_0}(0)]^{-1}$,
where $\rho_0$ is given in Lemma \ref{lem 1.1d}, we further perform a scaling
$$\bar{w}_{\varepsilon_n}(x)=z_{\varepsilon_n}^{-1}w_{\varepsilon_n}(z_{\varepsilon_n}^{-\frac{2}{N-2}}x),$$
then $$\bar{w}_{\varepsilon_n}(0)=z_{\varepsilon_n}^{-1}w_{\varepsilon_n}(0)=M_{\varepsilon_n}^{-1}[W_{\rho_0}(0)]w_{\varepsilon_n}(0)=W_{\rho_0}(0),$$
and $\bar{w}_{\varepsilon_n}$ satisfies the rescaled equation
\begin{equation}\label{1.4a}
\begin{aligned}
&-\Delta \bar{w}_{\varepsilon_n}+ \varepsilon_n^{-2\sigma}z_{\varepsilon_n}^{-\frac{4}{N-2}}\bar{w}_{\varepsilon_n}\\
=&\left(I_{\alpha}*(|\bar{w}_{\varepsilon_n}|^{\frac{N+\alpha}{N-2}} +\varepsilon_1^{-1}z_{\varepsilon_n}^{-\frac{N+\alpha}{N-2}}
G(\varepsilon_2z_{\varepsilon_n}\bar{w}_{\varepsilon_n}))\right)\\
\cdot&\left(\frac{N+\alpha}{N-2}|\bar{w}_{\varepsilon_n}|^{\frac{N+\alpha}{N-2}-2}\bar{w}_{\varepsilon_n}
+\varepsilon_1^{-1}\varepsilon_2z_{\varepsilon_n}^{-\frac{N+\alpha}{N-2}+1}g(\varepsilon_2z_{\varepsilon_n}\bar{w}_{\varepsilon_n})\right).\\
\end{aligned}
\end{equation}
Moreover, we have
$$
\|\bar{w}_{\varepsilon_n}\|_2^2=z_{\varepsilon_n}^{\frac{4}{N-2}}\|w_{\varepsilon_n}\|_2^2\ \ \ \mbox{and}\ \ \
\|\nabla \bar{w}_{\varepsilon_n}\|_2^2=\|\nabla w_{\varepsilon_n}\|_2^2.
$$

By Lemma \ref{lem 1.1d}, for any $\varepsilon_n \to \infty$, there exists $\rho\geq\rho_0$ such that
$$M_{\varepsilon_n}=w_{\varepsilon_n}(0)\to W_{\rho}(0)=\rho^{-\frac{N-2}{2}}W_1(0)\leq \rho_0^{-\frac{N-2}{2}}W_1(0)<\infty,$$
which yields that $M_{\varepsilon_n}\leq C$ for some $C>0$ and any large $\varepsilon_n>0$. Suppose that there exists a sequence $\varepsilon_n\to \infty$ such that $M_{\varepsilon_n}\to 0$.
By Lemma \ref{lem 1.1d}, up to a subsequence, $M_{\varepsilon_n}=w_{\varepsilon_n}(0) \to W_{\rho}(0)\neq 0$ for some $\rho \in (0,\infty)$. This leads to a contradiction. Therefore, there exists some $c>0$ such that $M_{\varepsilon_n}\geq c>0$.
Let
$$\xi_{\varepsilon_n}=z_{\varepsilon_n}^{-\frac{2}{N-2}}\varepsilon_n^{-\frac{(N+\alpha)-(N-2)q_2}{2[4+(N-2)q_2-(N+\alpha)]}},$$
then
$$\xi_{\varepsilon_n}\thicksim \varepsilon_n^{-\frac{(N+\alpha)-(N-2)q_2}{2[4+(N-2)q_2-(N+\alpha)]}},$$
and for large $\varepsilon_n>0$, the rescaled family of ground states
$$\bar{w}_{\varepsilon_n}(x)=\xi_{\varepsilon_n}^{\frac{N-2}{2}}v_{\varepsilon_n}(\xi_{\varepsilon_n}x)$$
satisfies
$$\|\nabla \bar{w}_{\varepsilon_n}\|^2\thicksim \int_{\mathbb R^N}
(I_{\alpha}*|\bar{w}_{\varepsilon_n}|^{\frac{N+\alpha}{N-2}})|\bar{w}_{\varepsilon_n}|^{\frac{N+\alpha}{N-2}}
\thicksim \|\bar{w}_{\varepsilon_n}\|_2^2\thicksim1,$$
and as $\varepsilon_n \to \infty$, $\bar{w}_{\varepsilon_n}$ converges to the extremal function $W_{\rho_0}$ in $\mathcal{D}^{1,2}(\mathbb R^N)$.
Then by Lemma \ref{lem 1.1d}, we also have
$\bar{w}_{\varepsilon_n}\to W_{\rho_0}$ in $L^2(\mathbb R^N)$. Thus we conclude that $\bar{w}_{\varepsilon_n}\to W_{\rho_0}$ in $H^1(\mathbb R^N)$.
Finally, the uniqueness of limit implies that $\bar{w}_{\varepsilon}\to W_{\rho_0}$ in $H^1(\mathbb R^N)$ as $\varepsilon\to \infty$.  The proof is complete.
\end{proof}

\begin{proof}[Proof of Theorem \ref{th1.3} in the case $N\ge 5$]
Since $w_{\varepsilon}\in \mathcal{N}_{\varepsilon}$, it follows from \eqref{1.5} and \eqref{2.1} that
$$\begin{aligned}
m_{\varepsilon}&=\frac{1}{2}\int_{\mathbb R^N}|\nabla w_{\varepsilon}|^2+\varepsilon^{-\sigma}|w_{\varepsilon}|^2\\
&\ \ \ -\frac{1}{2}\int_{\mathbb R^N}\int_{\mathbb R^N}
\frac{(|w_{\varepsilon}|^{\frac{N+\alpha}{N-2}}+\varepsilon_1^{-1}G(\varepsilon_2w_{\varepsilon}))(|w_{\varepsilon}|^{\frac{N+\alpha}{N-2}}
+\varepsilon_1^{-1}G(\varepsilon_2w_{\varepsilon}))}{|x-y|^{N-\alpha}}dxdy\\
&=\frac{2+\alpha}{2(N+\alpha)}\int_{\mathbb R^N}|\nabla w_{\varepsilon}|^2+O(\varepsilon^{-\sigma}).\\
\end{aligned}$$
Similarly, we also have
$$m_{\infty}=\frac{2+\alpha}{2(N+\alpha)}\int_{\mathbb R^N}|\nabla W_1|^2
=\frac{2+\alpha}{2(N-2)}\int_{\mathbb R^N}(I_{\alpha}*|W_1|^{\frac{N+\alpha}{N-2}})|W_1|^{\frac{N+\alpha}{N-2}}.$$
According to Lemma \ref{lem 1.1b}, we get that
$$\int_{\mathbb R^N}|\nabla W_1|^2-\int_{\mathbb R^N}|\nabla w_{\varepsilon}|^2=\frac{2(N+\alpha)}{2+\alpha}(m_{\infty}-m_{\varepsilon})+O(\varepsilon^{-\sigma})=O(\varepsilon^{-\sigma}).$$
Also, $\|\nabla W_1\|_2^2=\frac{N+\alpha}{N-2}\int_{\mathbb R^N}(I_{\alpha}*|W_1|^{\frac{N+\alpha}{N-2}})|W_1|^{\frac{N+\alpha}{N-2}}
=\left(\frac{N-2}{N+\alpha}\right)^{\frac{N-2}{2+\alpha}}S_{\alpha}^{\frac{N+\alpha}{2+\alpha}}$, we conclude that
$$\|\nabla w_{\varepsilon}\|_2^2=\left(\frac{N-2}{N+\alpha}\right)^{\frac{N-2}{2+\alpha}}S_{\alpha}^{\frac{N+\alpha}{2+\alpha}}+O(\varepsilon^{-\sigma}).$$
Then by $w_{\varepsilon}\in \mathcal{N}_{\varepsilon}$, we obtain
$$\begin{aligned}
\int_{\mathbb R^N}(I_{\alpha}*|w_{\varepsilon}|^{\frac{N+\alpha}{N-2}})|w_{\varepsilon}|^{\frac{N+\alpha}{N-2}}
=\frac{N-2}{N+\alpha}\|\nabla w_{\varepsilon}\|_2^2+O(\varepsilon^{-\sigma})=\left(\frac{N-2}{N+\alpha}S_{\alpha}\right)^{\frac{N+\alpha}{2+\alpha}}+O(\varepsilon^{-\sigma}).
\end{aligned}$$
Finally, according to Lemma \ref{lem 1.1c} and Lemma \ref{lem 1.1d}, we complete the proof.
\end{proof}

\section{The proof of main result for the cases $N=3,4$}
In this section, we prove that the asymptotic behavior of positive ground state solutions of \eqref{1.1} for the dimension $N=3,4$,
which are more complicated than $N\geq5$. In order to get a more accurate estimation, Kelvin transformation and Moser iteration technique are used during the process of proof.

\begin{lemma}\label{lem 4.1}
Assume that  (H1), (H2), (H3) and (H4) hold, then
$$
\varepsilon^{-\sigma}\gtrsim m_{\infty}-m_{\varepsilon}\gtrsim\left \{\begin{array}{rcl}
\varepsilon^{-\frac{3+\alpha-q_2}{2(q_2-1-\alpha)}}, \ \qquad  &if & N=3,\\
\left(\varepsilon\ln \varepsilon\right)^{-\frac{4+\alpha-2q_2}{2q_2-\alpha}},\  \ &if & N=4.
\end{array} \right.
$$ as $\varepsilon \to \infty$.
\end{lemma}
\begin{proof}
Arguing as in the proof of Lemma \ref{lem 1.1b}, we prove that $m_\infty-m_\varepsilon\lesssim \varepsilon^{-\sigma}$ as $\varepsilon\to \infty$.
To proceed, let $\rho>0$, $R$ is a large parameter, and $\eta_{R}\in C_0^{\infty}(\mathbb R)$ is a cut-off function such that
$\eta_{R}(r)=1$ for $|r|<R$, $0<\eta_{R}(r)<1$ for $R<|r|<2R$, $\eta_{R}(r)=0$ for $|r|>2R$ and $|\eta'_{R}|\leq \frac{2}{R}$. For $l\gg 1$, a straightforward
computation shows that
$$\int_{\mathbb R^N}|\nabla \eta_l W_1|^2=
\left\{
\begin{aligned}
&\frac{2(N+\alpha)m_{\infty}}{2+\alpha}+O(l^{-2}),
\ & if& \ N=4,\\
&\frac{2(N+\alpha)m_{\infty}}{2+\alpha}+O(l^{-1}),
\ & if&\ N=3.\\
\end{aligned}
\right.$$
$$\int_{\mathbb R^N}(I_{\alpha}*|\eta_{l}W_1|^{\frac{N+\alpha}{N-2}})|\eta_{l}W_1|^{\frac{N+\alpha}{N-2}}
=\frac{2(N-2)m_{\infty}}{2+\alpha}+O(l^{-N}).$$
$$\int_{\mathbb R^N}|\eta_l W_1|^2=
\left\{
\begin{aligned}
&\ln l(1+o(1)),
\ & if& \ N=4,\\
&l(1+o(1)),
\ & if&\ N=3.\\
\end{aligned}
\right.$$
By Lemma \ref{lem 2.1}, we find
$$\begin{aligned}
m_{\varepsilon}&\leq \sup_{t\geq 0}J_{\varepsilon}((\eta_{R}W_{\rho})_{t})=J_{\varepsilon}((\eta_{R}W_{\rho})_{t_{\varepsilon}})\\
&\leq \sup_{t\geq0}\left(\frac{t^{N-2}}{2}\int_{\mathbb R^N}|\nabla \left(\eta_{R}W_{\rho}\right)|^2
-\frac{t^{N+\alpha}}{2}\int_{\mathbb R^N}(I_{\alpha}*|\eta_{R}W_{\rho}|^{\frac{N+\alpha}{N-2}})|\eta_{R}W_{\rho}|^{\frac{N+\alpha}{N-2}}\right)\\
&-\Bigg(\int_{\mathbb R^N}(I_{\alpha}*|(\eta_{R}W_{\rho})_{t_{\varepsilon}}|^{\frac{N+\alpha}{N-2}})\varepsilon_1^{-1}G(\varepsilon_2(\eta_{R}W_{\rho})_{t_{\varepsilon}})dxdy
-\frac{1}{2}\varepsilon^{-\sigma}\int_{\mathbb R^N}|(\eta_{R}W_{\rho})_{t_{\varepsilon}}|^2dx\\
&\ \ \ \ \ \ +\frac{1}{2}\int_{\mathbb R^N}
(I_{\alpha}*\varepsilon_1^{-1}G(\varepsilon_2(\eta_{R}W_{\rho})_{t_{\varepsilon}}))\varepsilon_1^{-1}G(\varepsilon_2(\eta_{R}W_{\rho})_{t_{\varepsilon}})dxdy\Bigg)\\
&=(I)-(II),\\
\end{aligned}$$
where $\rho$ is a constant and $t_{\varepsilon}\in(0,\infty)$ is the unique critical point of the function $h(t)$ defined by
$$\begin{aligned}
h(t)&=\frac{t^{N-2}}{2}\int_{\mathbb R^N}|\nabla \left(\eta_{R}W_{\rho}\right)|^2
+\frac{\varepsilon^{-\sigma}t^{N}}{2}\int_{\mathbb R^N}|\eta_{R}W_{\rho}|^2\\
&-\frac{t^{N+\alpha}}{2}\int_{\mathbb R^N}(I_{\alpha}*(|\eta_{R}W_{\rho}|^{\frac{N+\alpha}{N-2}}+\varepsilon_1^{-1}G(\varepsilon_2\eta_{R}W_{\rho})))
(|\eta_{R}W_{\rho}|^{\frac{N+\alpha}{N-2}}+\varepsilon_1^{-1}G(\varepsilon_2\eta_{R}W_{\rho})).\\
\end{aligned}$$
That is, $t=t_{\varepsilon}$ solves the following equation $l_1(t)=l_2(t)$, where
$$l_1(t)=\frac{N-2}{2t^2}\int_{\mathbb R^N}|\nabla \left(\eta_{R}W_{\rho}\right)|^2,$$
and
$$\begin{aligned}
l_2(t)=\frac{(N+\alpha)t^{\alpha}}{2}\int_{\mathbb R^N}
&(I_{\alpha}*(|\eta_{R}W_{\rho}|^{\frac{N+\alpha}{N-2}}+\varepsilon_1^{-1}G(\varepsilon_2\eta_{R}W_{\rho})))\\
\cdot&(|\eta_{R}W_{\rho}|^{\frac{N+\alpha}{N-2}}+\varepsilon_1^{-1}G(\varepsilon_2\eta_{R}W_{\rho}))
-\frac{N\varepsilon^{-\sigma}}{2}\int_{\mathbb R^N}|\eta_{R}W_{\rho}|^2.\\
\end{aligned}$$
If $t_{\varepsilon}\geq1$, then
$$\begin{aligned}
\frac{N-2}{2t^2}\int_{\mathbb R^N}|\nabla \left(\eta_{R}W_{\rho}\right)|^2
&\geq\frac{N+\alpha}{2}\int_{\mathbb R^N}(I_{\alpha}*(|\eta_{R}W_{\rho}|^{\frac{N+\alpha}{N-2}}
+\varepsilon_1^{-1}G(\varepsilon_2\eta_{R}W_{\rho})))\\
&\cdot\left(|\eta_{R}W_{\rho}|^{\frac{N+\alpha}{N-2}}+\varepsilon_1^{-1}G(\varepsilon_2\eta_{R}W_{\rho})\right)
-\frac{N\varepsilon^{-\sigma}}{2}\int_{\mathbb R^N}|\eta_{R}W_{\rho}|^2,\\
\end{aligned}$$
and according to \eqref{1.5}, \eqref{2.1}, we have
$$
t_{\varepsilon}\leq\left(\frac{N-2}{N+\alpha}\frac{\int_{\mathbb R^N}|\nabla \left(\eta_{R}W_{\rho}\right)|^2}
{\int_{\mathbb R^N}(I_{\alpha}*|\eta_{R}W_{\rho}|^{\frac{N+\alpha}{N-2}})|\eta_{R}W_{\rho}|^{\frac{N+\alpha}{N-2}}}\right)^{\frac{1}{2}}.
$$
If $t_{\varepsilon}\leq1$, then
$$\begin{aligned}
\frac{N-2}{2t^2}\int_{\mathbb R^N}|\nabla \left(\eta_{R}W_{\rho}\right)|^2
&\leq \frac{N+\alpha}{2}\int_{\mathbb R^N}(I_{\alpha}*(|\eta_{R}W_{\rho}|^{\frac{N+\alpha}{N-2}}
+\varepsilon_1^{-1}G(\varepsilon_2\eta_{R}W_{\rho})))\\
&\cdot\left(|\eta_{R}W_{\rho}|^{\frac{N+\alpha}{N-2}}+\varepsilon_1^{-1}G(\varepsilon_2\eta_{R}W_{\rho})\right)
-\frac{N\varepsilon^{-\sigma}}{2}\int_{\mathbb R^N}|\eta_{R}W_{\rho}|^2,\\
\end{aligned}$$
and hence
$$\begin{aligned}
t_{\varepsilon}
\geq&\left(\frac{N-2}{N+\alpha}\frac{\int_{\mathbb R^N}|\nabla \left(\eta_{R}W_{\rho}\right)|^2}
{\int_{\mathbb R^N}(I_{\alpha}*|\eta_{R}W_{\rho}|^{\frac{N+\alpha}{N-2}})|\eta_{R}W_{\rho}|^{\frac{N+\alpha}{N-2}}}\right)^{\frac{1}{2}}\\
\cdot&\Bigg(1-\frac{
2\int_{\mathbb R^N}(I_{\alpha}*|\eta_{R}W_{\rho}|^{\frac{N+\alpha}{N-2}})\varepsilon_1^{-1}G(\varepsilon_2\eta_{R}W_{\rho})
-\frac{N\varepsilon^{-\sigma}}{N+\alpha}\int_{\mathbb R^N}|\eta_{R}W_{\rho}|^2}
{\int_{\mathbb R^N}(I_{\alpha}*|\eta_{R}W_{\rho}|^{\frac{N+\alpha}{N-2}})|\eta_{R}W_{\rho}|^{\frac{N+\alpha}{N-2}}}\\
&\ \ \ \ -\frac{\int_{\mathbb R^N}(I_{\alpha}*\varepsilon_1^{-1}G(\varepsilon_2\eta_{R}W_{\rho}))\varepsilon_1^{-1}G(\varepsilon_2\eta_{R}W_{\rho})}
{\int_{\mathbb R^N}(I_{\alpha}*|\eta_{R}W_{\rho}|^{\frac{N+\alpha}{N-2}})|\eta_{R}W_{\rho}|^{\frac{N+\alpha}{N-2}}}\Bigg)^{\frac{1}{2}}.\\
\end{aligned}$$
Therefore, by \eqref{1.5}, we obtain
$$\begin{aligned}
|t_{\varepsilon}-1|&\leq \left|\left(\frac{N-2}{N+\alpha}\frac{\int_{\mathbb R^N}|\nabla \left(\eta_{R}W_{\rho}\right)|^2}
{\int_{\mathbb R^N}(I_{\alpha}*|\eta_{R}W_{\rho}|^{\frac{N+\alpha}{N-2}})|\eta_{R}W_{\rho}|^{\frac{N+\alpha}{N-2}}}\right)^{\frac{1}{2}}-1\right|\\
&\ \ +\varepsilon^{-\sigma}\left(\frac{N-2}{N+\alpha}\frac{\int_{\mathbb R^N}|\nabla \left(\eta_{R}W_{\rho}\right)|^2}
{\int_{\mathbb R^N}(I_{\alpha}*|\eta_{R}W_{\rho}|^{\frac{N+\alpha}{N-2}})|\eta_{R}W_{\rho}|^{\frac{N+\alpha}{N-2}}}\right)^{\frac{1}{2}}\\
&\ \ \ \ \ \ \ \ \cdot\left(\frac{\psi(\rho)}{\int_{\mathbb R^N}
(I_{\alpha}*|\eta_{R}W_{\rho}|^{\frac{N+\alpha}{N-2}})|\eta_{R}W_{\rho}|^{\frac{N+\alpha}{N-2}}}\right),
\end{aligned}$$
where
$$
\psi(\rho):=C\int_{\mathbb R^N}(I_{\alpha}*|\eta_{R}W_{\rho}|^{\frac{N+\alpha}{N-2}})|\eta_{R}W_{\rho}|^{q_2}
-\frac{N}{N+\alpha}\int_{\mathbb R^N}|\eta_{R}W_{\rho}|^2.
$$

Set $l=R/\rho$, then for
$$(I):=\sup_{t\geq0}\left(\frac{t^{N-2}}{2}\int_{\mathbb R^N}|\nabla \left(\eta_{R}W_{\rho}\right)|^2
-\frac{t^{N+\alpha}}{2}\int_{\mathbb R^N}(I_{\alpha}*|\eta_{R}W_{\rho}|^{\frac{N+\alpha}{N-2}})|\eta_{R}W_{\rho}|^{\frac{N+\alpha}{N-2}}\right),$$
we conclude that $(I)$ gets its maximum value in $t_{\max}$ given by
$$t_{\max}=\left(\frac{N-2}{N+\alpha}\cdot\frac{\int_{\mathbb R^N}|\nabla \left(\eta_{R}W_{\rho}\right)|^2}
{\int_{\mathbb R^N}(I_{\alpha}*|\eta_{R}W_{\rho}|^{\frac{N+\alpha}{N-2}})|\eta_{R}W_{\rho}|^{\frac{N+\alpha}{N-2}}}\right)^{\frac{1}{2+\alpha}}.$$

Hence we have
$$\begin{aligned}
(I)
&=\frac{1}{2}\left(\frac{N-2}{N+\alpha}\right)^{\frac{N-2}{2+\alpha}}\left(1-\frac{N-2}{N+\alpha}\right)
\frac{\left(\int_{\mathbb R^N}|\nabla \left(\eta_{R}W_{\rho}\right)|^2\right)^{\frac{N+\alpha}{2+\alpha}}}
{\left(\int_{\mathbb R^N}(I_{\alpha}*|\eta_{R}W_{\rho}|^{\frac{N+\alpha}{N-2}})|\eta_{R}W_{\rho}|^{\frac{N+\alpha}{N-2}}\right)^{\frac{N-2}{2+\alpha}}}\\
&=\frac{1}{2}\left(\frac{N-2}{N+\alpha}\right)^{\frac{N-2}{2+\alpha}}\left(\frac{2+\alpha}{N+\alpha}\right)
\frac{\left(\frac{2(N+\alpha)m_{\infty}}{2+\alpha}+O(l^{-N+2})\right)^{\frac{N+\alpha}{2+\alpha}}}
{\left(\frac{2(N-2)m_{\infty}}{2+\alpha}+O(l^{-N})\right)^{\frac{N-2}{2+\alpha}}}\\
&=
\begin{cases}
m_{\infty}+O(l^{-2})\ \ \ if\ \ \ N=4,\\
m_{\infty}+O(l^{-1})\ \ \ if\ \ \ N=3.
\end{cases}\\
\end{aligned}$$

Since
$$\begin{aligned}
\varphi(\rho):&=
c\int_{\mathbb R^N}(I_{\alpha}*|\eta_{R}W_{\rho}|^{\frac{N+\alpha}{N-2}})|\eta_{R}W_{\rho}|^{q_2}
-\frac{1}{2}\int_{\mathbb R^N}|\eta_{R}W_{\rho}|^2\\
&=\rho^{\frac{N+\alpha}{2}-\frac{(N-2)q_2}{2}}c\int_{\mathbb R^N}(I_{\alpha}*|\eta_{l}W_{1}|^{\frac{N+\alpha}{N-2}})|\eta_{l}W_{1}|^{q_2}
-\frac{1}{2}\rho^2\int_{\mathbb R^N}|\eta_{l}W_{1}|^2\\
\end{aligned}$$
takes its maximum value $\varphi(\rho_l)$ at the unique point
$$\begin{aligned}
\rho_{l}:&=\left(\frac{[(N+\alpha)-(N-2)q_2]c}{2}
\cdot\frac{\int_{\mathbb R^N}(I_{\alpha}*|\eta_{l}W_{1}|^{\frac{N+\alpha}{N-2}})|\eta_{l}W_{1}|^{q_2}}{\int_{\mathbb R^N}|\eta_{l}W_{1}|^2}\right)^
{\frac{2}{4-(N+\alpha)+(N-2)q_2}}\\
&\sim
\begin{cases}
 (\ln l)^{-\frac{2}{4-(N+\alpha)+(N-2)q_2}}\ \ \ if\ \ \ N=4,\\
 l^{-\frac{2}{4-(N+\alpha)+(N-2)q_2}}\ \ \ \ \ \ \ \ if\ \ \ N=3,
\end{cases}\\
\end{aligned}$$
and
$$\begin{aligned}
\varphi(\rho_l)&=\sup_{\rho>0}\varphi(\rho)=\frac{4-[(N+\alpha)-(N-2)q_2]}{2[(N+\alpha)-(N-2)q_2]}\cdot\rho_{l}^2\cdot\int_{\mathbb R^N}|\eta_{l}W_{1}|^2\\
&\leq C\left(\int_{\mathbb R^N}|\eta_{l}W_{1}|^{2^*}dx\right)^{(1+\frac{\tilde{q}_2-2}{2^*-2})\cdot\frac{N+\alpha}{2N}\cdot\frac{4}{4-(N+\alpha)+(N-2)q_2}}\\
&\leq C\left(\int_{\mathbb R^N}|W_{1}|^{2^*}dx\right)^{(1+\frac{\tilde{q}_2-2}{2^*-2})\cdot\frac{N+\alpha}{2N}\cdot\frac{4}{4-(N+\alpha)+(N-2)q_2}}\leq C,\\
\end{aligned}$$
where the interpolation inequality \eqref{2.2} is used.
Besides, as $l \to \infty$, $$\int_{\mathbb R^N}(I_{\alpha}*|\eta_{l}W_{1}|^{\frac{N+\alpha}{N-2}})|\eta_{l}W_{1}|^{q_2}  \to\
\int_{\mathbb R^N}(I_{\alpha}*|W_{1}|^{\frac{N+\alpha}{N-2}})|W_{1}|^{q_2},$$
it follows that
$$\begin{aligned}
\varphi(\rho_l)&=C\left(\frac{\int_{\mathbb R^N}(I_{\alpha}*|\eta_{l}W_{1}|^{\frac{N+\alpha}{N-2}})|\eta_{l}W_{1}|^{q_2}}
{\int_{\mathbb R^N}|\eta_{l}W_{1}|^2}\right)^{\frac{4}{4-(N+\alpha)+(N-2)q_2}}\cdot\int_{\mathbb R^N}|\eta_{l}W_{1}|^2\\
&=
\begin{cases}
C(\ln l(1+o(1)))^{-\frac{(N+\alpha)-(N-2)q_2}{4-(N+\alpha)+(N-2)q_2}}\ \ if\ \ N=4,\\
C(l(1+o(1)))^{-\frac{(N+\alpha)-(N-2)q_2}{4-(N+\alpha)+(N-2)q_2}}\ \ \ \ \ if\ \ N=3.
\end{cases}\\
\end{aligned}$$

Similarly, we can prove that $\psi(\rho)$ is bounded and as $l\to +\infty$,
$$\begin{aligned}
|t_{\varepsilon}-1|
\to \frac{C\varepsilon^{-2\sigma}}{\int_{\mathbb R^N}
(I_{\alpha}*|\eta_{R}W_{\rho}|^{\frac{N+\alpha}{N-2}})|\eta_{R}W_{\rho}|^{\frac{N+\alpha}{N-2}}}.
\end{aligned}$$
Thus, for large $\varepsilon>0$, according to \eqref{1.5}, we have
$$\begin{aligned}
(II)&\geq\varepsilon^{-\sigma}\Bigg(\varphi(\rho_l)+(t_{\varepsilon}^{N}-1)\varphi(\rho_l)\\
&\ \ \ +t_{\varepsilon}^{N}(t_{\varepsilon}^{\alpha}-1)\rho_{l}^{\frac{N+\alpha}{2}-\frac{(N-2)q_2}{2}}\int_{\mathbb R^N}
(I_{\alpha}*|\eta_{l}W_{1}|^{\frac{N+\alpha}{N-2}})|\eta_{l}W_{1}|^{q_2}\Bigg)\\
&\sim \varepsilon^{-\sigma}
\begin{cases}
C(\ln l(1+o(1)))^{-\frac{(N+\alpha)-(N-2)q_2}{4-(N+\alpha)+(N-2)q_2}}\ \ if\ \ N=4,\\
C(l(1+o(1)))^{-\frac{(N+\alpha)-(N-2)q_2}{4-(N+\alpha)+(N-2)q_2}}\ \ \ \ \ if\ \ N=3.
\end{cases}\\
\end{aligned}$$

It follows that if $N=4$, then
$$m_{\varepsilon}\leq (I)-(II)\leq m_{\infty}+O(l^{-2})-C\varepsilon^{-\sigma}(\ln l)^{-\frac{(N+\alpha)-(N-2)q_2}{4-(N+\alpha)+(N-2)q_2}},$$
take $l=\varepsilon^{M}$, then
$$m_{\varepsilon}\leq m_{\infty}+O(\varepsilon^{-2M})-C\varepsilon^{-\sigma}M^{-\frac{(N+\alpha)-(N-2)q_2}{4-(N+\alpha)+(N-2)q_2}}
\left(\ln\varepsilon\right)^{-\frac{(N+\alpha)-(N-2)q_2}{4-(N+\alpha)+(N-2)q_2}},$$
if $M>-\frac{4+\alpha-2q_2}{2[2q_2-\alpha]}$, then $-2M<-\sigma=-\frac{(N+\alpha)-(N-2)q_2}{4+(N-2)q_2-(N+\alpha)}=-\frac{4+\alpha-2q_2}{2q_2-\alpha}$, and hence
$$m_{\varepsilon}\leq m_{\infty}-\varepsilon^{-\sigma}\left(\ln \varepsilon\right)^{-\frac{(N+\alpha)-(N-2)q_2}{4-(N+\alpha)+(N-2)q_2}}\bar{\omega}
=m_{\infty}-\left(\varepsilon\ln \varepsilon\right)^{-\frac{4+\alpha-2q_2}{2q_2-\alpha}}\bar{\omega},$$
for large $\varepsilon>0$, where
$$\bar{\omega}=\frac{1}{2}CM^{-\frac{(N+\alpha)-(N-2)q_2}{4-(N+\alpha)+(N-2)q_2}}=\frac{1}{2}CM^{-\frac{4+\alpha-2q_2}{2q_2-\alpha}}.$$

If $N=3$, then
$$m_{\varepsilon}\leq (I)-(II)\leq m_{\infty}+O(l^{-1})-C\varepsilon^{-\sigma}l^{-\frac{(N+\alpha)-(N-2)q_2}{4-(N+\alpha)+(N-2)q_2}},$$
take $l=\delta^{-1}\varepsilon^{-\tau}$, we get
$$m_{\varepsilon}\leq m_{\infty}+\delta O(\varepsilon^{\tau})-C\varepsilon^{-\sigma}\delta^{\frac{(N+\alpha)-(N-2)q_2}{4-(N+\alpha)+(N-2)q_2}}
\varepsilon^{\frac{\tau[(N+\alpha)-(N-2)q_2]}{4-(N+\alpha)+(N-2)q_2}},$$
if $\tau=-\frac{(N+\alpha)-(N-2)q_2}{2[2-(N+\alpha)+(N-2)q_2]}$,
then
$$m_{\varepsilon}\leq m_{\infty}+\left(\delta O(1)-C\delta^{\frac{(N+\alpha)-(N-2)q_2}{4-(N+\alpha)+(N-2)q_2}}\right)\varepsilon^{\tau}.$$
Since
$$1>\frac{(N+\alpha)-(N-2)q_2}{4-(N+\alpha)+(N-2)q_2},$$
it follows that for some $\delta>0$, there exists $\bar{\omega}>0$ such that
$$m_{\varepsilon}\leq m_{\infty}-\varepsilon^{-\frac{(N+\alpha)-(N-2)q_2}{2[2-(N+\alpha)+(N-2)q_2]}}\bar{\omega}=m_{\infty}-\varepsilon^{-\frac{3+\alpha-q_2}{2(q_2-1-\alpha)}}\bar{\omega}.$$
This completes the proof.
\end{proof}

\smallskip

In the lower dimension cases $N=3,4$, by Lemma 2.6, arguing as in the proof of \cite[Lemma 5.7]{MM}, we show that there exists $\xi_\varepsilon\in (0,\infty)$ with $\xi_\varepsilon\to 0$ such that
\begin{equation}\label{4.1}
\tilde{w}_{\varepsilon}(x)=\xi_{\varepsilon}^{\frac{N-2}{2}}w_{\varepsilon}(\xi_{\varepsilon}x),
\end{equation}
converges to $W_1$ in $D^{1,2}(\mathbb R^N)$ and $L^{2^*}(\mathbb R^N)$ as $\varepsilon\to \infty$. Clearly, $\tilde w=\tilde w_\varepsilon$ satisfies
\begin{equation}\label{4.2}
\begin{aligned}
-\Delta \tilde{w}+ \varepsilon^{-\sigma}\xi_{\varepsilon}^2\tilde{w}=&\left(I_{\alpha}*(|\tilde{w}|^{\frac{N+\alpha}{N-2}}
+\varepsilon_1^{-1}\xi_{\varepsilon}^{\frac{N+\alpha}{2}}G(\varepsilon_2\xi_{\varepsilon}^{-\frac{N-2}{2}}\tilde{w}))\right)\\
\cdot &\left(\frac{N+\alpha}{N-2}|\tilde{w}|^{\frac{N+\alpha}{N-2}-2}+\varepsilon_1^{-1}\varepsilon_2\xi_{\varepsilon}^{\frac{2+\alpha}{2}}
g(\varepsilon_2\xi_{\varepsilon}^{-\frac{N-2}{2}}\tilde{w})
\tilde{w}^{-1}\right)\tilde{w}.\\
\end{aligned}
\end{equation}
And the corresponding energy functional of the above equation is given by
$$\begin{aligned}
\tilde{J}_{\varepsilon}(\tilde{w})=&\frac{1}{2}\int_{\mathbb R^N}|\nabla \tilde{w}|^2+\varepsilon^{-\sigma}\xi_{\varepsilon}^2|\tilde{w}|^2\\
&\quad -\frac{1}{2}\int_{\mathbb R^N}(I_{\alpha}*(|\tilde{w}|^{\frac{N+\alpha}{N-2}}+\varepsilon_1^{-1}\xi_{\varepsilon}^{\frac{N+\alpha}{2}}G(\varepsilon_2\xi_{\varepsilon}^{-\frac{N-2}{2}}\tilde{w})))\\ &\ \ \ \ \ \ \ \ \ \ \ \ \ \cdot(|\tilde{w}|^{\frac{N+\alpha}{N-2}}+\varepsilon_1^{-1}\xi_{\varepsilon}^{\frac{N+\alpha}{2}}G(\varepsilon_2\xi_{\varepsilon}^{-\frac{N-2}{2}}\tilde{w})).\\
\end{aligned}$$
Clearly, we have $J_{\varepsilon}(w_\varepsilon)=\tilde{J}_{\varepsilon}(\tilde{w_\varepsilon})$.
Furthermore, we have the following Lemma.
\begin{lemma}\label{lem 4.2}
Let $w_\varepsilon(x)=\varepsilon^{-\frac{(N-2)\sigma}{4}}v_\varepsilon(\varepsilon^{-\frac{\sigma}{2}}x),$
 and let $\tilde w_\varepsilon$ be  given by \eqref{4.1}, then the following statements hold true
\begin{itemize}
\item[(i)] $\|\nabla v_\varepsilon\|_2^2=\|\nabla w_\varepsilon\|_2^2=\|\nabla \tilde{w_\varepsilon}\|_2^2$, \   $\varepsilon^{\sigma}\|v_\varepsilon\|_2^2=\|w_\varepsilon\|_2^2=\xi_{\varepsilon}^2\|\tilde{w_\varepsilon}\|_2^2$,
\smallskip
\item[(ii)] $\int_{\mathbb R^N}(I_{\alpha}*|\tilde{w_\varepsilon}|^{\frac{N+\alpha}{N-2}})|\tilde{w_\varepsilon}|^{\frac{N+\alpha}{N-2}}
            =\int_{\mathbb R^N}(I_{\alpha}*|w_\varepsilon|^{\frac{N+\alpha}{N-2}})|w_\varepsilon|^{\frac{N+\alpha}{N-2}}
            =\int_{\mathbb R^N}(I_{\alpha}*|v_\varepsilon|^{\frac{N+\alpha}{N-2}})|v_\varepsilon|^{\frac{N+\alpha}{N-2}}$.
\end{itemize}
\end{lemma}
The corresponding Nehari and Pohoz\v aev's identities of \eqref{4.2} are as follows
$$
\begin{aligned}
\int_{\mathbb R^N}|\nabla \tilde{w}_{\varepsilon}|^2+\varepsilon^{-\sigma}\xi_{\varepsilon}^2\int_{\mathbb R^N}|\tilde{w}_{\varepsilon}|^2
=&\int_{\mathbb R^N}(I_{\alpha}*(|\tilde{w}_{\varepsilon}|^{\frac{N+\alpha}{N-2}}
+\varepsilon_1^{-1}\xi_{\varepsilon}^{\frac{N+\alpha}{2}}G(\varepsilon_2\xi_{\varepsilon}^{-\frac{N-2}{2}}\tilde{w}_{\varepsilon})))\\
\cdot&(\frac{N+\alpha}{N-2}|\tilde{w}_{\varepsilon}|^{\frac{N+\alpha}{N-2}}
+\varepsilon_1^{-1}\varepsilon_2\xi_{\varepsilon}^{\frac{2+\alpha}{2}}g(\varepsilon_2\xi_{\varepsilon}^{-\frac{N-2}{2}}\tilde{w}_{\varepsilon})\tilde{w}_{\varepsilon}),\\
\end{aligned}
$$
and
$$\begin{aligned}
&\frac{N-2}{2}\int_{\mathbb R^N}|\nabla \tilde{w}_{\varepsilon}|^2+\frac{N}{2}\varepsilon^{-\sigma}\xi_{\varepsilon}^2\int_{\mathbb R^N}|\tilde{w}_{\varepsilon}|^2\\
=&\frac{N+\alpha}{2}\int_{\mathbb R^N}
(I_{\alpha}*(|\tilde{w}_{\varepsilon}|^{\frac{N+\alpha}{N-2}}+\varepsilon_1^{-1}\xi_{\varepsilon}^{\frac{N+\alpha}{2}}G(\varepsilon_2\xi_{\varepsilon}^{-\frac{N-2}{2}}\tilde{w}_{\varepsilon})))\\
&\ \ \ \ \ \ \cdot(|\tilde{w}_{\varepsilon}|^{\frac{N+\alpha}{N-2}}+\varepsilon_1^{-1}\xi_{\varepsilon}^{\frac{N+\alpha}{2}}G(\varepsilon_2\xi_{\varepsilon}^{-\frac{N-2}{2}}\tilde{w}_{\varepsilon})),\\
\end{aligned}$$
so, we obtain
\begin{equation}\label{4.3}
\begin{aligned}
\varepsilon^{-\sigma}\xi_{\varepsilon}^2\|\tilde{w}_{\varepsilon}\|_2^2
=&\int_{\mathbb R^N}(I_{\alpha}*(|\tilde{w}_{\varepsilon}|^{\frac{N+\alpha}{N-2}}
+\varepsilon_1^{-1}\xi_{\varepsilon}^{\frac{N+\alpha}{2}}G(\varepsilon_2\xi_{\varepsilon}^{-\frac{N-2}{2}}\tilde{w}_{\varepsilon})))\\
\cdot&(\frac{N+\alpha}{2}\varepsilon_1^{-1}\xi_{\varepsilon}^{\frac{N+\alpha}{2}}G(\varepsilon_2\xi_{\varepsilon}^{-\frac{N-2}{2}}\tilde{w}_{\varepsilon})
-\frac{N-2}{2}\varepsilon_1^{-1}\varepsilon_2\xi_{\varepsilon}^{\frac{2+\alpha}{2}}g(\varepsilon_2\xi_{\varepsilon}^{-\frac{N-2}{2}}\tilde{w}_{\varepsilon})\tilde{w}_{\varepsilon}).\\
\end{aligned}
\end{equation}

To control norm $\|\tilde{w}_{\varepsilon}\|_2^2$, we introduce the following lemma:
\begin{lemma}\label{lem 4.2a}
There exists a constant $c>0$ such that
\begin{equation}\label{4.4}
\tilde{w}_{\varepsilon}(x)\geq c|x|^{-(N-2)}\mathrm{exp} (-\varepsilon^{-\sigma/2}\xi_{\varepsilon}|x|),\ \ |x|\geq1.
\end{equation}
\end{lemma}
The proof  of Lemma \ref{lem 4.2a} is similar to that in \cite[Lemma 4.8] {MM2014}. Consequently, we have the following two lemmas.

\begin{lemma}\label{lem 4.2b}
If $N=3$, then $\|\tilde{w}_{\varepsilon}\|_2^2\gtrsim \varepsilon^{\sigma/2}\xi_{\varepsilon}^{-1}$.
\end{lemma}

\begin{lemma}\label{lem 4.2c}
If $N=4$, then $\|\tilde{w}_{\varepsilon}\|_2^2\gtrsim -\ln (\varepsilon^{-\sigma}\xi_{\varepsilon}^{2})$.
\end{lemma}

Now, we are in the position to state our result concerning the best decay estimation of $\tilde{w}_{\varepsilon}$.
\begin{proposition}\label{Pro 4.1}
Assume that  (H1), (H2), (H3) and (H4)  hold, then
there exists a constant $C>0$ such that for large $\varepsilon>0$,
$$\tilde{w}_{\varepsilon}(x)\leq C(1+|x|)^{-(N-2)}, \ \ \ x \in \mathbb R^N.$$
\end{proposition}
The proof of Proposition \ref{Pro 4.1} will be given in the Appendix.
From Proposition \ref{Pro 4.1}, we obtain the following

\begin{lemma}\label{lem 4.3}
Assume that (H1), (H2), (H3) and (H4) hold,  then
$\int_{\mathbb R^N}(I_{\alpha}*|\tilde{w}_{\varepsilon}|^{\frac{N+\alpha}{N-2}})|\tilde{w}_{\varepsilon}|^{q_2}\sim1$ as $\varepsilon \to \infty$. Furthermore, for any $s>\frac{N}{N-2}$, $\tilde{w}_{\varepsilon}\to W_1$
in $L^{s}(\mathbb R^N)$ as $\varepsilon \to \infty$.
\end{lemma}
\begin{proof}
By virtue of Proposition \ref{Pro 4.1}, there exists a constant $C>0$ such that for all large $\varepsilon>0$,
$$\tilde{w}_{\varepsilon}(x)\leq \frac{C}{(1+|x|)^{N-2}},\ \ \forall x \in \mathbb R^N.$$
Therefore, for any  $s>\frac{N}{N-2}$, $\tilde{w}_{\varepsilon}$ is bounded uniformly in $L^{s}(\mathbb R^N)$
for large $\varepsilon>0$, and by dominated convergence theorem, we get that $\tilde{w}_{\varepsilon} \to W_1$ in $L^{s}(\mathbb R^N)$ as $\varepsilon \to \infty$.

Since $\frac{2Nq_2}{N+\alpha}>\frac{N}{N-2}$,
by the Hardy-Littlewood-Sobolev inequality, we find that as $\varepsilon \to \infty$,
$$\int_{\mathbb R^N}(I_{\alpha}*|\tilde{w}_{\varepsilon}-W_1|^{\frac{N+\alpha}{N-2}})
|\tilde{w}_{\varepsilon}-W_1|^{q_2}\leq C\|\tilde{w}_{\varepsilon}-W_1\|_{\frac{2Nq_2}{N+\alpha}}^{q_2}\to 0.$$
Arguing as in \cite[Lemma 2.1]{LW-1}, we show that
$$\begin{aligned}
&\quad\lim_{\varepsilon \to \infty}\int_{\mathbb R^N}(I_{\alpha}*|\tilde{w}_{\varepsilon}|^{\frac{N+\alpha}{N-2}})|\tilde{w}_{\varepsilon}|^{q_2}
-\int_{\mathbb R^N}(I_{\alpha}*|\tilde{w}_{\varepsilon}-W_1|^{\frac{N+\alpha}{N-2}})|\tilde{w}_{\varepsilon}-W_1|^{q_2}\\
&\qquad =\int_{\mathbb R^N}(I_{\alpha}*|W_1|^{\frac{N+\alpha}{N-2}})|W_1|^{q_2},\\
\end{aligned}$$
therefore, it follows that
$$\lim_{\varepsilon \to \infty}\int_{\mathbb R^N}(I_{\alpha}*|\tilde{w}_{\varepsilon}|^{\frac{N+\alpha}{N-2}})|\tilde{w}_{\varepsilon}|^{q_2}
=\int_{\mathbb R^N}(I_{\alpha}*|W_1|^{\frac{N+\alpha}{N-2}})|W_1|^{q_2}.$$
Particularly, we have that $\int_{\mathbb R^N}(I_{\alpha}*|\tilde{w}_{\varepsilon}|^{\frac{N+\alpha}{N-2}})|\tilde{w}_{\varepsilon}|^{q_2}\sim 1$ as $\varepsilon\to \infty$.
\end{proof}

\begin{proof}[Proof of Theorem \ref{th1.3} in the case $N=4,3$]
By \eqref{3.1} and Lemma \ref{lem 4.2}, we have $\tau(\tilde{w}_{\varepsilon})=\tau(w_\varepsilon)$, so by \eqref{1.5} and \eqref{4.3}, we get
\begin{equation}\label{4.28}
\begin{aligned}
m_{\infty}&\leq \sup_{t\geq 0}\tilde{J}_{\infty}((\tilde{w}_{\varepsilon})_t)\leq \sup_{t\geq 0}\tilde{J}_{\varepsilon}((\tilde{w}_{\varepsilon})_t)
-\frac{t_{\varepsilon}^N}{2}\int_{\mathbb R^N}\varepsilon^{-\sigma} \xi_{\varepsilon}^2 |\tilde{w}_{\varepsilon}|^2\\
&\quad+\frac{1}{2}\int_{\mathbb R^N}
(I_{\alpha}*(2|(\tilde{w}_{\varepsilon})_{t_{\varepsilon}}|^{\frac{N+\alpha}{N-2}}
+\varepsilon_1^{-1}\xi_{\varepsilon}^{\frac{N+\alpha}{2}}G(\varepsilon_2\xi_{\varepsilon}^{-\frac{N-2}{2}}(\tilde{w}_{\varepsilon})_{t_{\varepsilon}})))\\
&\ \ \ \ \ \ \cdot\varepsilon_1^{-1}\xi_{\varepsilon}^{\frac{N+\alpha}{2}}G(\varepsilon_2\xi_{\varepsilon}^{-\frac{N-2}{2}}(\tilde{w}_{\varepsilon})_{t_{\varepsilon}})\\
&\leq m_{\varepsilon}+C\varepsilon^{-\sigma}\xi_{\varepsilon}^{\frac{N+\alpha}{2}-\frac{(N-2)q_2}{2}}
\left(\int_{\mathbb R^N}(I_{\alpha}*|\tilde{w}_{\varepsilon}|^{\frac{N+\alpha}{N-2}})|\tilde{w}_{\varepsilon}|^{q_2}+o_\varepsilon(1)\right),\\
\end{aligned}
\end{equation}
which, together with Lemma \ref{lem 4.1},  implies that
$$\begin{aligned}
&\quad\xi_{\varepsilon}^{\frac{N+\alpha}{2}-\frac{(N-2)q_2}{2}}\left(\int_{\mathbb R^N}(I_{\alpha}*|\tilde{w}_{\varepsilon}|^{\frac{N+\alpha}{N-2}})|\tilde{w}_{\varepsilon}|^{q_2}+o_\varepsilon(1)\right)\\
&\gtrsim \varepsilon^{\sigma}\delta_{\varepsilon}
\gtrsim
\begin{cases}
\begin{split}
&(\ln\varepsilon)^{-\frac{4+\alpha-2q_2}{2q_2-\alpha}}\ \ \ &if\ \ \ N=4,\\
&\varepsilon^{-\frac{(3+\alpha-q_2)^2}{2(q_2-1-\alpha)(q_2+1-\alpha)}}\ \ \ &if\ \ \ N=3,
\end{split}
\end{cases}
\end{aligned}$$
where $\delta_{\varepsilon}:=m_{\infty}-m_{\varepsilon}$.
Therefore, from Lemma \ref{lem 4.3}, we have
$$\xi_{\varepsilon}^{\frac{N+\alpha}{2}-\frac{(N-2)q_2}{2}}
\gtrsim
\begin{cases}
\begin{split}
&(\ln\varepsilon)^{-\frac{4+\alpha-2q_2}{2q_2-\alpha}}\ \ \ &if\ \ \ N=4,\\
&\varepsilon^{-\frac{(3+\alpha-q_2)^2}{2(q_2-1-\alpha)(q_2+1-\alpha)}}\ \ \ &if\ \ \ N=3.
\end{split}
\end{cases}
$$

On the other hand, if $N=3$, then by \eqref{1.5}, \eqref{4.3}, Lemma \ref{lem 4.3} and Lemma \ref{lem 4.2b}, we have
$$\xi_{\varepsilon}^{\frac{q_2+1-\alpha}{2}}\lesssim \frac{1}{\|\tilde{w}_{\varepsilon}\|_2^2}\lesssim \varepsilon^{-\sigma/2}\xi_{\varepsilon},$$
then $\xi_{\varepsilon}^{\frac{q_2-1-\alpha}{2}}\lesssim \varepsilon^{-\sigma/2}$.
Hence, we have
\begin{equation}\label{4.29}
\xi_{\varepsilon}\lesssim \varepsilon^{-\frac{3+\alpha-q_2}{(q_2-1-\alpha)(q_2+1-\alpha)}}.
\end{equation}

If $N=4$, then by \eqref{1.5}, \eqref{4.3}, Lemma \ref{lem 4.3} and Lemma \ref{lem 4.2c}, we have
$$\xi_{\varepsilon}^{q_2-\frac{\alpha}{2}}\lesssim \frac{1}{\|\tilde{w}_{\varepsilon}\|_2^2}\lesssim \frac{1}{-\ln (\varepsilon^{-\sigma}\xi_{\varepsilon}^2)}.$$
Note that
$$-\ln (\varepsilon^{-\sigma}\xi_{\varepsilon}^2)=\sigma \ln\varepsilon+2\ln\frac{1}{\xi_{\varepsilon}}\geq \sigma \ln\varepsilon,$$
it follows that
$$\xi_{\varepsilon}^{q_2-\frac{\alpha}{2}}\lesssim \frac{1}{\|\tilde{w}_{\varepsilon}\|_2^2}\lesssim \left(\ln\varepsilon\right)^{-1},$$
so we obtain that
\begin{equation}\label{4.30}
\xi_{\varepsilon}\lesssim \left(\ln\varepsilon\right)^{-\frac{2}{2q_2-\alpha}}.
\end{equation}
Thus, it follows from \eqref{4.28}-\eqref{4.30} and Lemma \ref{lem 4.3} that
$$\delta_{\varepsilon}=m_{\infty}-m_{\varepsilon}\lesssim \varepsilon^{-\sigma}\xi_{\varepsilon}^{\frac{N+\alpha}{2}-\frac{(N-2)q_2}{2}}
\lesssim
\begin{cases}
\begin{split}
&\left(\varepsilon\ln\varepsilon\right)^{-\frac{4+\alpha-2q_2}{2q_2-\alpha}},\ \ \ &if\ \ \ N=4,\\
&\varepsilon^{-\frac{3+\alpha-q_2}{2(q_2-1-\alpha)}},\ \ \ &if\ \ \ N=3,
\end{split}
\end{cases}
$$
which together with Lemma \ref{lem 4.1} imply that
\begin{equation}\label{4.31}
\delta_{\varepsilon}\thicksim \varepsilon^{-\sigma}\xi_{\varepsilon}^{\frac{N+\alpha}{2}-\frac{(N-2)q_2}{2}}
\thicksim
\begin{cases}
\begin{split}
&\left(\varepsilon\ln\varepsilon\right)^{-\frac{4+\alpha-2q_2}{2q_2-\alpha}},\ \ \ &if\ \ \ N=4,\\
&\varepsilon^{-\frac{3+\alpha-q_2}{2(q_2-1-\alpha)}},\ \ \ &if\ \ \ N=3.
\end{split}
\end{cases}
\end{equation}
Arguing as in  the proof of \cite[Theorem 2.2]{MM}, we also have
\begin{equation}\label{4.32}
\|\nabla \tilde{w}_{\varepsilon}\|_2^2=\left(\frac{N-2}{N+\alpha}\right)^{\frac{N-2}{2+\alpha}}S_{\alpha}^{\frac{N+\alpha}{2+\alpha}}
+
\begin{cases}
\begin{split}
&O\left(\left(\varepsilon\ln\varepsilon\right)^{-\frac{4+\alpha-2q_2}{2q_2-\alpha}}\right),\  &if\   N=4,\\
&O\left(\varepsilon^{-\frac{3+\alpha-q_2}{2(q_2-1-\alpha)}}\right),\  &if\  N=3,
\end{split}
\end{cases}
\end{equation}
$$\begin{aligned}
\int_{\mathbb R^N}(I_{\alpha}*|\tilde{w}_{\varepsilon}|^{\frac{N+\alpha}{N-2}})|\tilde{w}_{\varepsilon}|^{\frac{N+\alpha}{N-2}}
&=\left(\frac{N-2}{N+\alpha}S_{\alpha}\right)^{\frac{N+\alpha}{2+\alpha}}\\
&\ \ \ +
\begin{cases}
\begin{split}
&O\left(\left(\varepsilon\ln\varepsilon\right)^{-\frac{4+\alpha-2q_2}{2q_2-\alpha}}\right),\ \ \ &if\ \ \ N=4,\\
&O\left(\varepsilon^{-\frac{3+\alpha-q_2}{2(q_2-1-\alpha)}}\right),\ \ \ &if\ \ \ N=3.\\
\end{split}
\end{cases}
\end{aligned}$$

Finally, by \eqref{4.3}, Lemma \ref{lem 4.2b} and Lemma \ref{lem 4.2c}, we obtain
$$
\|\tilde{w}_{\varepsilon}\|_2^2
\thicksim
\begin{cases}
\begin{split}
&\ln\varepsilon,\ \ \ &if\ \ \ N=4,\\
&\varepsilon^{\frac{3+\alpha-q_2}{2(q_2-1-\alpha)}},\ \ \ &if\ \ \ N=3.
\end{split}
\end{cases}
$$
By \eqref{1.6}, \eqref{1.8} and \eqref{4.1}, we complete the proof of Theorem \ref{th1.3} in the case $N=4,3$.
\end{proof}

\section{Appendix}

In this appendix, we give a proof of Proposition \ref{Pro 4.1}, which concerns the optimal decay estimate of rescaled ground state, and establish the existence of ground state of our problem \eqref{1.1}.

\subsection{The Proof of Proposition \ref{Pro 4.1}}

In this subsection, we give the proof of Proposition \ref{Pro 4.1}.  The proof is similar to that of \cite[Proposition 6.17]{MM}, for the readers' convenience, we outline  the proof below.

In order to prove Proposition \ref{Pro 4.1}, we also introduce the Kelvin transform of $\tilde{w}_{\varepsilon}$ defined in \eqref{4.1} given as follows
\begin{equation}\label{4.5}
K[\tilde{w}_{\varepsilon}](x):=|x|^{-(N-2)}\tilde{w}_{\varepsilon}\left(\frac{x}{|x|^2}\right).
\end{equation}
It is easy to see that $\|K[\tilde{w}_{\varepsilon}]\|_{L^{\infty}(B_1)}\lesssim 1$ implies that
$$\sup_{|x|\geq 1}\tilde{w}_{\varepsilon}(x)\lesssim |x|^{-(N-2)}.$$
At first,   we  show that there exists $\varepsilon_0>0$ such that
\begin{equation}\label{4.6}
\sup_{\varepsilon \in (\varepsilon_0, \infty)}\parallel K[\tilde{w}_{\varepsilon}]\parallel_{L^{\infty}(B_1)}< \infty.
\end{equation}
It is easy to verify that $K[\tilde{w}_{\varepsilon}]$ satisfies
\begin{equation}\label{4.7}
\begin{aligned}
-\Delta K[\tilde{w}_{\varepsilon}]+&\frac{\varepsilon^{-\sigma}\xi_{\varepsilon}^2}{|x|^4}K[\tilde{w}_{\varepsilon}]
=\frac{1}{|x|^4}\left(I_{\alpha}*(|\tilde{w}_{\varepsilon}|^{\frac{N+\alpha}{N-2}}
+\varepsilon_1^{-1}\xi_{\varepsilon}^{\frac{N+\alpha}{2}}G(\varepsilon_2\xi_{\varepsilon}^{-\frac{N-2}{2}}\tilde{w}_{\varepsilon}))\right)\\
&\cdot\left(\frac{N+\alpha}{N-2}|\tilde{w}_{\varepsilon}|^{\frac{N+\alpha}{N-2}-2}
+\varepsilon_1^{-1}\varepsilon_2\xi_{\varepsilon}^{\frac{2+\alpha}{2}}g(\varepsilon_2\xi_{\varepsilon}^{-\frac{N-2}{2}}\tilde{w}_{\varepsilon})\tilde{w}^{-1}_{\varepsilon}\right)K[\tilde{w}_{\varepsilon}].\\
\end{aligned}
\end{equation}

\begin{lemma}\label{lem 1}
Assume that (H1), (H2), (H3) and (H4) hold, then
there exist constants $L_0>0$ and $C_0>0$ such that for any large
$\varepsilon>0$ and $|x|\geq L_0 \varepsilon^{\sigma/2}\xi_{\varepsilon}^{-1}$,
$$\tilde{w}_{\varepsilon}(x)\leq C_0 \varepsilon^{-\frac{(N-2)\sigma}{4}}\xi_{\varepsilon}^{\frac{(N-2)}{2}} \mathrm{exp} (-\frac{1}{2}\varepsilon^{-\sigma/2}\xi_{\varepsilon}|x|).$$
\end{lemma}
\begin{proof}
By \eqref{3.6}, if $|x|\geq L_0 \varepsilon^{\sigma/2}\xi_{\varepsilon}^{-1}$ with $L_0>0$ large enough, then we have
$$\begin{aligned}
&\frac{N+\alpha}{N-2}\left(I_{\alpha}*|\tilde{w}_{\varepsilon}|^{\frac{N+\alpha}{N-2}}\right)(x)|\tilde{w}_{\varepsilon}|^{\frac{N+\alpha}{N-2}-2}(x)\\
=&\frac{N+\alpha}{N-2}\xi_{\varepsilon}^2\left(I_{\alpha}*|w_{\varepsilon}|^{\frac{N+\alpha}{N-2}}\right)(\xi_{\varepsilon}x)|w_{\varepsilon}|^{\frac{N+\alpha}{N-2}-2}(\xi_{\varepsilon}x)\\
\leq&C\xi_{\varepsilon}^2|\xi_{\varepsilon}x|^{-\frac{N^2-N\alpha+4\alpha}{2(N-2)}}
\leq C\xi_{\varepsilon}^2|\xi_{\varepsilon}L_0 \varepsilon^{\sigma/2}\xi_{\varepsilon}^{-1}|^{-\frac{N^2-N\alpha+4\alpha}{2(N-2)}}\\
=&C\xi_{\varepsilon}^2L_0^{-\frac{N^2-N\alpha+4\alpha}{2(N-2)}}\varepsilon^{\frac{-\sigma(N^2-N\alpha+4\alpha)}{4(N-2)}}
\leq\frac{\varepsilon^{-\sigma}\xi_{\varepsilon}^2}{8},\\
\end{aligned}$$
where we use the fact that
$\frac{N^2-N\alpha+4\alpha}{4(N-2)}>1$.
For $|x|\geq L_0 \varepsilon^{\sigma/2}\xi_{\varepsilon}^{-1}$, combining \eqref{1.5}, \eqref{3.4} with \eqref{3.5}, we get
$$\begin{aligned}
&\left(I_{\alpha}*|\tilde{w}_{\varepsilon}|^{\frac{N+\alpha}{N-2}}\right)(x)
\varepsilon_1^{-1}\varepsilon_2\xi_{\varepsilon}^{\frac{2+\alpha}{2}}g(\varepsilon_2\xi_{\varepsilon}^{-\frac{N-2}{2}}\tilde{w}_{\varepsilon}(x))\tilde{w}^{-1}_{\varepsilon}(x)\\
\leq &C\xi_{\varepsilon}^{N}\left(I_{\alpha}*|w_{\varepsilon}|^{\frac{N+\alpha}{N-2}}\right)(\xi_{\varepsilon}x)\varepsilon_1^{-1}
\left(|\varepsilon_2\xi_{\varepsilon}^{-\frac{N-2}{2}}\tilde{w}_{\varepsilon}(x)|^{q_1}+|\varepsilon_2\xi_{\varepsilon}^{-\frac{N-2}{2}}\tilde{w}_{\varepsilon}(x)|^{q_2}\right)\tilde{w}^{-2}_{\varepsilon}(x)\\
\leq &C\varepsilon^{-\sigma}\xi_{\varepsilon}^2|\xi_{\varepsilon}x|^{-N+\alpha}|\xi_{\varepsilon}x|^{-\frac{N(q_2-2)}{2}}\leq \frac{\varepsilon^{-\sigma}\xi_{\varepsilon}^2}{8},\\
\end{aligned}$$
and
$$\begin{aligned}
&\left(I_{\alpha}*\varepsilon_1^{-1}\xi_{\varepsilon}^{\frac{N+\alpha}{2}}G(\varepsilon_2\xi_{\varepsilon}^{-\frac{N-2}{2}}\tilde{w}_{\varepsilon})\right)(x)
\varepsilon_1^{-1}\varepsilon_2\xi_{\varepsilon}^{\frac{2+\alpha}{2}}g(\varepsilon_2\xi_{\varepsilon}^{-\frac{N-2}{2}}\tilde{w}_{\varepsilon}(x))\tilde{w}^{-1}_{\varepsilon}(x)\\
\leq &C\xi_{\varepsilon}^{N}\varepsilon_1^{-1}
\left(I_{\alpha}*\varepsilon_1^{-1}G(\varepsilon_2w_{\varepsilon})\right)(\xi_{\varepsilon} x)\\ &\cdot\left(|\varepsilon_2\xi_{\varepsilon}^{-\frac{N-2}{2}}\tilde{w}_{\varepsilon}(x)|^{q_1}+|\varepsilon_2\xi_{\varepsilon}^{-\frac{N-2}{2}}\tilde{w}_{\varepsilon}(x)|^{q_2}\right)\tilde{w}^{-2}_{\varepsilon}(x)\\
\leq &C\xi_{\varepsilon}^{2}\varepsilon_1^{-2}
\left(I_{\alpha}*(|\varepsilon_2w_{\varepsilon}|^{q_1}+|\varepsilon_2w_{\varepsilon}|^{q_2})\right)(\xi_{\varepsilon} x)\\ &\cdot\left(|\varepsilon_2w_{\varepsilon}(\xi_{\varepsilon} x)|^{q_1}
+|\varepsilon_2w_{\varepsilon}(\xi_{\varepsilon} x)|^{q_2}\right)w^{-2}_{\varepsilon}(\xi_{\varepsilon} x)\\
\leq &C\varepsilon^{-2\sigma}\xi_{\varepsilon}^2|\xi_{\varepsilon}x|^{-N+\alpha}|\xi_{\varepsilon}x|^{-\frac{N(q_2-2)}{2}}\leq \frac{\varepsilon^{-\sigma}\xi_{\varepsilon}^2}{8},
\end{aligned}$$
and
$$\begin{aligned}
&\left(I_{\alpha}*\varepsilon_1^{-1}\xi_{\varepsilon}^{\frac{N+\alpha}{2}}G(\varepsilon_2\xi_{\varepsilon}^{-\frac{N-2}{2}}\tilde{w}_{\varepsilon})\right)(x)
\frac{N+\alpha}{N-2}|\tilde{w}_{\varepsilon}|^{\frac{N+\alpha}{N-2}-2}\\
\leq &C\xi_{\varepsilon}^{\frac{N-\alpha}{2}}\left(I_{\alpha}*\varepsilon_1^{-1}G(\varepsilon_2w_{\varepsilon})\right)(\xi_{\varepsilon} x)
|\xi_{\varepsilon}^{\frac{N-2}{2}}w_{\varepsilon}(\xi_{\varepsilon}x)|^{\frac{N+\alpha}{N-2}-2}\\
\leq &C\xi_{\varepsilon}^{2}\varepsilon_1^{-1}
\left(I_{\alpha}*(|\varepsilon_2w_{\varepsilon}|^{q_1}+|\varepsilon_2w_{\varepsilon}|^{q_2})\right)(\xi_{\varepsilon}x)
|w_{\varepsilon}(\xi_{\varepsilon}x)|^{\frac{N+\alpha}{N-2}-2}\\
\leq &C\varepsilon^{-\sigma}\xi_{\varepsilon}^{2}|\xi_{\varepsilon}x|^{-N+\alpha}|\xi_{\varepsilon}x|^{-\frac{N}{2}\left(\frac{N+\alpha}{N-2}-2\right)}
\leq \frac{\varepsilon^{-\sigma}\xi_{\varepsilon}^2}{8},
\end{aligned}$$
where $(I_{\alpha}*|w_{\varepsilon}|^{q_i})(x)\leq |x|^{-N+\alpha}$ is similar to \eqref{3.5} as $2<q_i<\frac{N+\alpha}{N-2}$.
Therefore, we obtain
$$-\Delta \tilde{w}_{\varepsilon}(x)+\frac{1}{2}\varepsilon^{-\sigma}\xi_{\varepsilon}^2\tilde{w}_{\varepsilon}(x)\leq 0,\ \ \mbox{for\ \  all}\ \ |x|\geq L_0 \varepsilon^{\sigma/2}\xi_{\varepsilon}^{-1}.$$

Let $R>L_0 \varepsilon^{\sigma/2}\xi_{\varepsilon}^{-1}$ and we introduce a positive function
$$\Psi_{R}(r):=\mathrm{exp}\left(-\frac{1}{2}\varepsilon^{-\sigma/2}\xi_{\varepsilon}\left(r-L_0\varepsilon^{\sigma/2}\xi_{\varepsilon}^{-1}\right)\right)
+\mathrm{exp}\left(\frac{1}{2}\varepsilon^{-\sigma/2}\xi_{\varepsilon}\left(r-R\right)\right).
$$
Arguing as in \cite[Lemma 3.2]{AIIKN2019}, the comparison principle implies that if $L_0\varepsilon^{\sigma/2}\xi_{\varepsilon}^{-1}\leq|x|\leq R$, then
$$\tilde{w}_{\varepsilon}(x)\leq CL_0^{-\frac{N-2}{2}}\varepsilon^{-\frac{(N-2)\sigma}{4}}\xi_{\varepsilon}^{\frac{N-2}{2}}\Psi_{R}(|x|).$$
Since $R>L_0\varepsilon^{\sigma/2}\xi_{\varepsilon}^{-1}$ is arbitrary, taking $R\to \infty$, we find that
$$\tilde{w}_{\varepsilon}(x)\leq C L_0^{-\frac{N-2}{2}}\mathrm{exp}\left(\frac{L_0}{2}\right)\varepsilon^{-\frac{(N-2)\sigma}{4}}\xi_{\varepsilon}^{\frac{N-2}{2}} \mathrm{exp} (-\frac{1}{2}\varepsilon^{-\sigma/2}\xi_{\varepsilon}|x|),$$
for all $|x|\geq L_0 \varepsilon^{\sigma/2}\xi_{\varepsilon}^{-1}$. The proof is complete.
\end{proof}

From Lemma \ref{lem 1}, we see that if $|x|\leq L_0^{-1} \varepsilon^{-\sigma/2}\xi_{\varepsilon}$, then
\begin{equation}\label{4.8}
K[\tilde{w}_{\varepsilon}](x)\lesssim \frac{1}{|x|^{N-2}} \varepsilon^{-\frac{(N-2)\sigma}{4}}\xi_{\varepsilon}^{\frac{N-2}{2}} \mathrm{exp} (-\frac{1}{2}\varepsilon^{-\sigma/2}\xi_{\varepsilon}|x|^{-1}).
\end{equation}
Let $a(x)=\frac{\varepsilon^{-\sigma}\xi_{\varepsilon}^2}{|x|^4},$
and
$$\begin{aligned}
b(x)=\frac{1}{|x|^4}&\left(I_{\alpha}*(|\tilde{w}_{\varepsilon}|^{\frac{N+\alpha}{N-2}}
+\varepsilon_1^{-1}\xi_{\varepsilon}^{\frac{N+\alpha}{2}}G(\varepsilon_2\xi_{\varepsilon}^{-\frac{N-2}{2}}\tilde{w}_{\varepsilon}))\right)\\
\cdot&\left(\frac{N+\alpha}{N-2}|\tilde{w}_{\varepsilon}|^{\frac{N+\alpha}{N-2}-2}
+\varepsilon_1^{-1}\varepsilon_2\xi_{\varepsilon}^{\frac{2+\alpha}{2}}g(\varepsilon_2\xi_{\varepsilon}^{-\frac{N-2}{2}}\tilde{w}_{\varepsilon})\tilde{w}^{-1}_{\varepsilon}\right).\\
\end{aligned}$$
Then \eqref{4.7} reads as
$$-\Delta K[\tilde{w}_{\varepsilon}]+a(x)K[\tilde{w}_{\varepsilon}]=b(x)K[\tilde{w}_{\varepsilon}].$$
We shall apply the Moser iteration to prove \eqref{4.6}. For any $v \in H^1(B_4)$, \eqref{4.8} implies that
\begin{equation}\label{4.9}
\int_{B_4}\frac{\varepsilon^{-\sigma}\xi_{\varepsilon}^2}{|x|^4}K[\tilde{w}_{\varepsilon}](x)|v(x)|dx< \infty.
\end{equation}

In the following,  we denote
$\eta_i=\frac{N+\alpha-(N-2)q_i}{2}$, since $q_1\le q_2$,  we then obtain
$$\begin{aligned}
|b(x)|
&\lesssim \frac{1}{|x|^4}\left\{\left(I_{\alpha}*|\tilde{w}_{\varepsilon}|^{\frac{N+\alpha}{N-2}}\right)(\frac{x}{|x|^2})|w_\varepsilon|^{\frac{N+\alpha}{N-2}-2}(\frac{x}{|x|^2})\right.\\
&+\sum_{i=1}^2\varepsilon^{-\sigma}\xi_\varepsilon^{\eta_i}\left(I_{\alpha}*|\tilde{w}_{\varepsilon}|^{\frac{N+\alpha}{N-2}}\right)(\frac{x}{|x|^2})|\tilde w_\varepsilon|^{q_i-2}(\frac{x}{|x|^2})\\
&+\sum_{i=1}^2\varepsilon^{-\sigma}\xi_\varepsilon^{\eta_i}\left(I_{\alpha}*|\tilde{w}_{\varepsilon}|^{q_i}\right)(\frac{x}{|x|^2})|\tilde w_\varepsilon|^{\frac{N+\alpha}{N-2}-2}(\frac{x}{|x|^2})\\
&+\left.\sum_{i,j=1}^2\varepsilon^{-2\sigma}\xi_\varepsilon^{\eta_i+\eta_j}\left(I_{\alpha}*|\tilde{w}_{\varepsilon}|^{q_i}\right)(\frac{x}{|x|^2})|\tilde w_\varepsilon|^{q_j-2}(\frac{x}{|x|^2})\right\}.\\
\end{aligned}$$

Since $\tilde{w}_{\varepsilon} \to W_1$ in $L^{2^*}(\mathbb R^N)$ as $\varepsilon \to \infty$, we have
$\lim_{\varepsilon\to\infty}\|K[\tilde{w}_{\varepsilon}]-K[W_1]\|_{L^{2^*}(\mathbb R^N)}=0$.
Arguing as in \cite[Proposition 5.14]{MM}, we show that
\begin{equation}\label{4.12}
\begin{aligned}
&\lim_{\varepsilon\to \infty}\Bigg\|\frac{1}{|x|^4}\left(I_{\alpha}*|\tilde{w}_{\varepsilon}|^{\frac{N+\alpha}{N-2}}\right)(\frac{x}{|x|^2})\tilde{w}_{\varepsilon}^{\frac{N+\alpha}{N-2}-2}(\frac{x}{|x|^2})\Bigg\|_{L^{\frac{N}{2}}}\\
&=\Bigg\|\frac{1}{|x|^4}\left(I_{\alpha}*|W_1|^{\frac{N+\alpha}{N-2}}\right)(\frac{x}{|x|^2})W_1^{\frac{N+\alpha}{N-2}-2}(\frac{x}{|x|^2})\Bigg\|_{L^{\frac{N}{2}}}.
\end{aligned}
\end{equation}
In the following,  we denote $\eta_{\varepsilon}=L_0^{-1}\varepsilon^{-\sigma/2}\xi_{\varepsilon}$.
\begin{lemma}\label{lem 2}
Assume that (H1), (H2), (H3) and (H4) hold, then for any $i=1,2$,
\begin{equation}\label{4.13}
\lim_{\varepsilon\to\infty}\int_{|x|\leq4}\left|\frac{\varepsilon^{-\sigma}
\xi_{\varepsilon}^{\eta_i}}{|x|^4}
\left(I_{\alpha}*|\tilde{w}_{\varepsilon}|^{\frac{N+\alpha}{N-2}}\right)(\frac{x}{|x|^2})|\tilde w_\varepsilon|^{q_i-2}(\frac{x}{|x|^2})\right|^{\frac{N}{2}}dx=0.
\end{equation}
\end{lemma}
\begin{proof}
We divide the integral into two parts:
$$
I_{\varepsilon}^{(1)}\left(\frac{N}{2}\right):=\varepsilon^{-\frac{N\sigma}{2}}
\xi_{\varepsilon}^{\frac{\eta_iN}{2}}\int_{|z|\geq \eta_\varepsilon^{-1}}\left|
\left(I_{\alpha}*|\tilde{w}_{\varepsilon}|^{\frac{N+\alpha}{N-2}}\right)(z)|\tilde{w}_{\varepsilon}|^{q_i-2}(z)\right|^{\frac{N}{2}}dz,
$$
$$
I_{\varepsilon}^{(2)}\left(\frac{N}{2}\right):=\varepsilon^{-\frac{N\sigma}{2}}
\xi_{\varepsilon}^{\frac{\eta_iN}{2}}\int_{1/4\le |z|\le\eta_\varepsilon^{-1}}\left|
\left(I_{\alpha}*|\tilde{w}_{\varepsilon}|^{\frac{N+\alpha}{N-2}}\right)(z)|\tilde{w}_{\varepsilon}|^{q_i-2}(z)\right|^{\frac{N}{2}}dz.
$$
By \eqref{3.5}, we have
$$
(I_{\alpha}*|\tilde{w}_{\varepsilon}|^{\frac{N+\alpha}{N-2}})(z)
=\xi_{\varepsilon}^{\frac{N-\alpha}{2}}(I_{\alpha}*|w_{\varepsilon}|^{\frac{N+\alpha}{N-2}})(\xi_{\varepsilon}z)
\lesssim \xi_{\varepsilon}^{-\frac{N-\alpha}{2}}|z|^{-N+\alpha}.
$$
Therefore, by  Lemma \ref{lem 1}, we get
$$
\begin{aligned}
I_{\varepsilon}^{(1)}\left(\frac{N}{2}\right)&\lesssim\varepsilon^{-\frac{N\sigma}{2}}
\xi_{\varepsilon}^{\frac{\eta_iN}{2}}\int_{|z|\geq\eta_{\varepsilon}^{-1}}\left(\xi_{\varepsilon}^{-\frac{N-\alpha}{2}}|z|^{-N+\alpha}\right)^{\frac{N}{2}}\\
& \ \ \ \cdot\left(\varepsilon^{-\frac{(N-2)\sigma}{4}}\xi_{\varepsilon}^{\frac{N-2}{2}}
\mathrm{exp}\left(-\frac{1}{2}\varepsilon^{-\sigma/2}\xi_{\varepsilon}|z|\right)\right)^{\frac{N(q_i-2)}{2}}dz\\
&\lesssim\varepsilon^{-\frac{N\sigma(N-2)(q_i-2)}{8}-\frac{N\sigma(N-\alpha)}{4}}
\xi_{\varepsilon}^{\frac{\eta_iN}{2}-\frac{N-\alpha}{2}\frac{N}{2}+\frac{(N-2)(q_i-2)N}{4}+\frac{N(N-\alpha)}{2}-N}\\
& \ \ \ \cdot\int_{L_0}^{\infty}|r|^{-\frac{N(N-\alpha)}{2}+N-1}\mathrm{exp}\left(-\frac{N(q_i-2)}{4}|r|\right)dr\\
&\lesssim \varepsilon^{\frac{-N\sigma}{8}[(N-2)(q_i-2)+2(N-\alpha)]}\to \ 0,
\end{aligned}$$
as $\varepsilon \to \infty$, here we use the fact that
$$\frac{\eta_iN}{2}-\frac{N-\alpha}{2}\cdot\frac{N}{2}+\frac{(N-2)(q_i-2)N}{4}+\frac{N(N-\alpha)}{2}-N=0.$$
And by using the H\"{o}lder inequality, the Hardy-Littlewood-Sobolev inequality and the facts that
$$2\leq \frac{N+\alpha}{N-2}\cdot\frac{2N}{2N-(N-2)q_i+2\alpha}\leq \frac{2N}{N-2},$$
and
$$\frac{N\eta_i}{2}+\frac{N(N+\alpha)}{4}-\frac{N[2N-(N-2)q_i+2\alpha]}{4}=0,$$
we obtain that as $\varepsilon \to \infty$,
$$
\begin{aligned}
I_{\varepsilon}^{(2)}\left(\frac{N}{2}\right)&\lesssim \varepsilon^{-\frac{N\sigma}{2}}
\xi_{\varepsilon}^{\frac{N\eta_i}{2}}\left(\int_{ \frac{1}{4}\leq|z|\leq \eta^{-1}_{\varepsilon}}
|\tilde{w}_{\varepsilon}|^{2^*}dz\right)^{\frac{(N-2)(q_i-2)}{4}}\\
&\cdot \left(\int_{ \frac{1}{4}\leq|z|\leq \eta_{\varepsilon}^{-1}}
\left|(I_{\alpha}*|\tilde{w}_{\varepsilon}|^{\frac{N+\alpha}{N-2}})(z)\right|^{\frac{2N}{2N-(N-2)q_i}}dz\right)^{\frac{2N-(N-2)q_i}{4}}\\
&\lesssim\varepsilon^{-\frac{N\sigma}{2}}
\xi_{\varepsilon}^{\frac{N\eta_i}{2}}\left(\int_{ \frac{1}{4}\leq|z|\leq \eta_{\varepsilon}^{-1}}
|\tilde{w}_{\varepsilon}|^{\frac{N+\alpha}{N-2}\cdot\frac{2N}{2N-(N-2)q_i+2\alpha}}dz\right)^{\frac{2N-(N-2)q_i+2\alpha}{4}}\\
&\lesssim\varepsilon^{-\frac{N\sigma}{2}}
\xi_{\varepsilon}^{\frac{N\eta_i}{2}+\frac{N(N+\alpha)}{4}-\frac{N[2N-(N-2)q_i+2\alpha]}{4}}\\
&\cdot \left(\int_{\mathbb R^N}|w_{\varepsilon}|^{\frac{N+\alpha}{N-2}\cdot\frac{2N}{2N-(N-2)q_i+2\alpha}}dz\right)^{\frac{2N-(N-2)q_i+2\alpha}{4}}\\
&\lesssim\varepsilon^{-\frac{N\sigma}{2}}
 \ \to\ 0,
\end{aligned}
$$
Thus the conclusion follows since $\lim_{\varepsilon\to\infty}I_{\varepsilon}^{(1)}\left(\frac{N}{2}\right)+I_{\varepsilon}^{(2)}\left(\frac{N}{2}\right)=0$. The proof is complete.
\end{proof}

Similarly, we also have the following lemma.

\begin{lemma}\label{lem 2b}
Assume that (H1), (H2), (H3) and (H4) hold, then for any $i=1,2$,
\begin{equation}\label{4.13b}
\lim_{\varepsilon\to\infty}\int_{|x|\leq4}\left|\frac{\varepsilon^{-\sigma}
\xi_{\varepsilon}^{\eta_i}}{|x|^4}
\left(I_{\alpha}*|\tilde{w}_{\varepsilon}|^{q_i}\right)(\frac{x}{|x|^2})|\tilde w_\varepsilon|^{\frac{N+\alpha}{N-2}-2}(\frac{x}{|x|^2})\right|^{\frac{N}{2}}dx=0.
\end{equation}
\end{lemma}

\begin{lemma}\label{lem 2a}
Assume that (H1), (H2), (H3) and (H4) hold, then for any $i,j=1,2$,
\begin{equation}\label{4.14}
\lim_{\varepsilon\to\infty}\int_{|x|\leq4}\left|\frac{\varepsilon^{-2\sigma}\xi_\varepsilon^{\eta_i+\eta_j}}{|x|^4}
\left(I_{\alpha}*|\tilde{w}_{\varepsilon}|^{q_i}\right)(\frac{x}{|x|^2})|\tilde w_\varepsilon|^{q_j-2}(\frac{x}{|x|^2})\right|^{\frac{N}{2}}dx=0.
\end{equation}
\end{lemma}
\begin{proof}
Similar to Lemma \ref{lem 2}, we divide the integral into two parts: $I_{\varepsilon}^{(1)}\left(\frac{N}{2}\right)$ denotes the integral on $\{x| \ |x|\leq\eta_{\varepsilon}\}$, $I_{\varepsilon}^{(2)}\left(\frac{N}{2}\right)$ denotes the integral on $\{x| \ \eta_{\varepsilon}\leq|x|\leq4\}$.

For the first integral part, by virtue of Lemma \ref{lem 2.7}, for all $|z|\geq L_0\varepsilon^{\sigma/2}\xi_{\varepsilon}^{-1}$,
$$\begin{aligned}
(I_{\alpha}*|\tilde{w}_{\varepsilon}|^{q_i})(z)&=\xi_{\varepsilon}^{\frac{(N-2)q_i}{2}-\alpha}(I_{\alpha}*|w_{\varepsilon}|^{q_i})(\xi_{\varepsilon}z)
\lesssim\xi_{\varepsilon}^{\frac{(N-2)q_i}{2}-N}|z|^{-N+\alpha},
\end{aligned}$$
for $q_i>2$. Therefore, by Lemma \ref{lem 1}, it follows that
$$\begin{aligned}
I_{\varepsilon}^{(1)}\left(\frac{N}{2}\right)&\lesssim \varepsilon^{-\frac{N\sigma}{2}-\frac{N\sigma(N-2)(q_j-2)}{8}-\frac{N\sigma(N-\alpha)}{4}}\\
&\ \ \ \cdot\xi_{\varepsilon}^{\frac{N(\eta_i+\eta_j)}{2}+\frac{N}{2}(\frac{(N-2)q_i}{2}-N)+\frac{N(N-2)(q_j-2)}{4}+\frac{N(N-\alpha)}{2}-N}\\
&\ \ \ \cdot\int_{|z|\geq L_0}|z|^{-\frac{N(N-\alpha)}{2}}
\mathrm{exp}\left(-\frac{N(q_j-2)}{4}|z|\right)dz\\
&\lesssim \varepsilon^{-\frac{N\sigma}{2}-\frac{N\sigma(N-2)(q_j-2)}{8}-\frac{N\sigma(N-\alpha)}{4}} \to 0,\\
\end{aligned}$$
where
$$\frac{N(\eta_i+\eta_j)}{2}+\frac{N}{2}\left(\frac{(N-2)q_i}{2}-N\right)+\frac{N(N-2)(q_j-2)}{4}+\frac{N(N-\alpha)}{2}-N=0.$$

For the second integral part, by the H\"{o}lder inequality and the Hardy-Littlewood-Sobolev inequality, we get that
$$\begin{aligned}
I_{\varepsilon}^{(2)}\left(\frac{N}{2}\right)
&\lesssim \varepsilon^{-N\sigma}\xi_{\varepsilon}^{\frac{(\eta_i+\eta_j)N}{2}}\left(\int_{\frac{1}{4}\leq|z|\leq \eta_{\varepsilon}^{-1}}
|\tilde{w}_{\varepsilon}|^{\frac{N(q_j-2)}{2}\cdot\frac{4}{(N-2)(q_j-2)}}dz\right)^{\frac{(N-2)(q_j-2)}{4}}\\
&\ \ \ \cdot
\left(\int_{\frac{1}{4}\leq|z|\leq \eta_{\varepsilon}^{-1}}
\left|(I_{\alpha}*|\tilde{w}_{\varepsilon}|^{q_i})\right|^{\frac{N}{2}\cdot\frac{4}{2N-(N-2)q_j}}dz\right)^{\frac{2N-(N-2)q_j}{4}}.\\
\end{aligned}$$

In what follows, , we split into two cases:

{\bf Case 1}: When $\alpha\leq N-2$, we have $\frac{2\alpha}{N-2}\leq2<q_i, q_j<\frac{N+\alpha}{N-2}$,
hence $$2\leq \frac{2Nq_i}{2N-(N-2)q_j+2\alpha}\leq\frac{2N}{N-2}.$$
Also, by using the fact
$$\frac{(\eta_i+\eta_j)N}{2}+\frac{N(N-2)q_i}{4}-\frac{N(2N-(N-2)q_j+2\alpha)}{4}=0,$$
thus by the boundedness of $w_\varepsilon$ in $H^1(\mathbb R^N)$, we obtain
$$\begin{aligned}
I_{\varepsilon}^{(2)}\left(\frac{N}{2}\right)&\lesssim\varepsilon^{-N\sigma}
\xi_{\varepsilon}^{\frac{(\eta_i+\eta_j)N}{2}+\frac{N(N-2)q_i}{4}-\frac{N(2N-(N-2)q_j+2\alpha)}{4}}\\
&\ \ \ \cdot\left(\int_{\frac{1}{4}\leq|z|\leq \eta_{\varepsilon}^{-1}}
|w_{\varepsilon}|^{\frac{2Nq_i}{2N-(N-2)q_j+2\alpha}}dz\right)^{\frac{2N-(N-2)q_j+2\alpha}{4}}\\
&\lesssim\varepsilon^{-N\sigma}\to0\ \ \mbox{as}\ \ \varepsilon\to\infty.
\end{aligned}$$
{\bf Case 2}: When $\alpha>N-2$, we have $2<\frac{2\alpha}{N-2}<\frac{N+\alpha}{N-2}$. Let
$$\begin{aligned}
Q_{\varepsilon}:&=\int_{\frac{1}{4}\leq|z|\leq \eta_{\varepsilon}^{-1}}
\left|(I_{\alpha}*|\tilde{w}_{\varepsilon}|^{q_i})\right|^{\frac{2N}{2N-(N-2)q_j}}dz\\
&\leq \left(\int_{\frac{1}{4}\leq|z|\leq \eta_{\varepsilon}^{-1}}
\left|(I_{\alpha}*|\tilde{w}_{\varepsilon}|^{q_i})\right|^{\frac{2Nq}{2N-(N-2)q_j}}dz\right)^{\frac{1}{q}}
\left(\int_{\frac{1}{4}\leq|z|\leq \eta_{\varepsilon}^{-1}}dz\right)^{\frac{q-1}{q}}\\
&\lesssim\varepsilon^{\frac{N\sigma(q-1)}{2q}}\xi_{\varepsilon}^{-\frac{N(q-1)}{q}}
\left(\int_{\frac{1}{4}\leq|z|\leq \eta_{\varepsilon}^{-1}}|\tilde{w}_{\varepsilon}|^{q_is}dz\right)^{\frac{1}{qs}\cdot\frac{2Nq}{2N-(N-2)q_j}}\\
&\lesssim\varepsilon^{\frac{N\sigma(q-1)}{2q}}
\xi_{\varepsilon}^{-\frac{N(q-1)}{q}+\frac{N-2}{2}\cdot\frac{2Nq_i}{2N-(N-2)q_j}-\frac{N}{qs}\cdot\frac{2Nq}{2N-(N-2)q_j}}
|w_{\varepsilon}|_{q_is}^{\frac{2Nq_i}{2N-(N-2)q_j}},\\
\end{aligned}$$
where $s=\frac{2Nq}{2N-(N-2)q_j+2\alpha q}$.

When $2<q_i\leq \frac{2\alpha}{N-2}$ and $2<q_j<\frac{N+\alpha}{N-2}$, we can choose $q>1$ satisfying $$\frac{2N-(N-2)q_j}{2(N-\alpha)}<q<\infty,$$
then
\begin{equation}\label{4.16}
1<s:=\frac{2Nq}{2N-(N-2)q_j+2\alpha q}<\frac{N}{\alpha},
\end{equation}
so $2<q_i<q_is<\frac{Nq_i}{\alpha}\leq\frac{2N}{N-2}$. Since $w_{\varepsilon}$ is bounded in $H^1(\mathbb R^N)$, we have
$$Q_{\varepsilon}\lesssim\varepsilon^{\frac{N\sigma(q-1)}{2q}}
\xi_{\varepsilon}^{-\frac{N(q-1)}{q}+\frac{N-2}{2}\cdot\frac{2Nq_i}{2N-(N-2)q_j}-\frac{N}{qs}\cdot\frac{2Nq}{2N-(N-2)q_j}}.$$
When $\frac{2\alpha}{N-2}<q_i<\frac{N+\alpha}{N-2}$ and $2<q_j<\frac{N+\alpha}{N-2}$,
we get that
$$\quad\quad\quad\quad \frac{2N-(N-2)q_j}{q_i(N-2)-2\alpha}>1\ \ \ \mbox{and}\ \ \ \frac{2N-(N-2)q_j}{2(N-\alpha)}<q<\frac{2N-(N-2)q_j}{q_i(N-2)-2\alpha},$$
which together with \eqref{4.16} imply that
 $2<q_is<\frac{2N}{N-2}$, then it follows from the boundedness of $w_{\varepsilon}$ in $H^1(\mathbb R^N)$ that
$$Q_{\varepsilon}\lesssim\varepsilon^{\frac{N\sigma(q-1)}{2q}}
\xi_{\varepsilon}^{-\frac{N(q-1)}{q}+\frac{N-2}{2}\cdot\frac{2Nq_i}{2N-(N-2)q_j}-\frac{N}{qs}\cdot\frac{2Nq}{2N-(N-2)q_j}}.$$
Therefore, from case 2, we get
$$\begin{aligned}
I_{\varepsilon}^{(2)}\left(\frac{N}{2}\right)&\lesssim \varepsilon^{-N\sigma}
\varepsilon^{\frac{N\sigma(q-1)}{2q}\cdot\frac{2N-(N-2)q_j}{4}}\\
&\ \ \cdot\xi_{\varepsilon}^{\frac{(\eta_i+\eta_j)N}{2}}
\xi_{\varepsilon}^{\left(-\frac{N(q-1)}{q}+\frac{N-2}{2}\cdot\frac{2Nq_i}{2N-(N-2)q_j}-\frac{N}{qs}\cdot\frac{2Nq}{2N-(N-2)q_j}\right)\frac{2N-(N-2)q_j}{4}}\\
&\lesssim \varepsilon^{-N\sigma (1-\frac{q-1}{2q}\cdot\frac{2N-(N-2)q_j}{4})}\to\ 0,\ \ \ \mbox{as}\ \ \ \varepsilon\to \infty.
\end{aligned}$$
Besides, $q>1$ implies that $1-\frac{q-1}{2q}\cdot\frac{2N-(N-2)q_j}{4}>0$ and $\Gamma(\alpha,q_i,q_j)=0$, where
$$\begin{aligned}
\Gamma(\alpha,q_i,q_j)
&=\frac{(\eta_i+\eta_j)N}{2}-\frac{N(q-1)(2N-(N-2)q_j)}{4q}+\frac{N(N-2)q_i}{4}-\frac{N^2}{2s}\\
&=\frac{(\eta_i+\eta_j)N}{2}-\frac{N(2N-(N-2)q_j)}{4}+\frac{N(2N-(N-2)q_j)}{4q}\\
&\ \ \ +\frac{N(N-2)q_i}{4}-\frac{N^2}{2}\cdot\frac{2N-(N-2)q_j+2\alpha q}{2Nq}=0.\\
\end{aligned}$$
Thus, the conclusion follows since $\lim_{\varepsilon\to\infty}I_{\varepsilon}^{(1)}\left(\frac{N}{2}\right)+I_{\varepsilon}^{(2)}\left(\frac{N}{2}\right)=0$.
\end{proof}

From \eqref{4.8}-\eqref{4.9}, \eqref{4.13}-\eqref{4.14} and Lemma \ref{lem 2.8} (i), we deduce that for any $r>1$, there exists
$\varepsilon_r>0$ such that
\begin{equation}\label{4.17}
\sup_{\varepsilon \in (\varepsilon_r,\infty)}\parallel K[\tilde{w}_{\varepsilon}]^r\parallel_{H^1(B_1)}\leq C_r.
\end{equation}

To verify the conditions in Lemma \ref{lem 2.8} (ii), we introduce the following three lemmas:

\begin{lemma}\label{lem 4}
Assume that (H1), (H2), (H3) and (H4) hold, then there exists $r_0>\frac{N}{2}$ such that
\begin{equation}\label{4.18}
\lim_{\varepsilon\to\infty}\int_{|x|\leq4}\left|\frac{\varepsilon^{-\sigma}
\xi_{\varepsilon}^{\eta_i}}{|x|^4}\left(I_{\alpha}*|\tilde{w}_{\varepsilon}|^{\frac{N+\alpha}{N-2}}\right)(\frac{x}{|x|^2})
|\tilde{w}_{\varepsilon}|^{q_i-2}(\frac{x}{|x|^2})\right|^{r_0}dx=0.
\end{equation}
\end{lemma}
\begin{proof}
Similar to Lemma \ref{lem 2}, we divide the integral into two parts: $I_{\varepsilon}^{(1)}\left(r_0\right)$ denotes the integral on $\{x| \ |x|\leq\eta_{\varepsilon}\}$, $I_{\varepsilon}^{(2)}\left(r_0\right)$ denotes the integral on $\{x|\ \eta_{\varepsilon}\leq|x|\leq4\}$.

For the first integral part, by H\"{o}lder inequality, we have
\begin{equation}\label{4.19}
\begin{aligned}
I_{\varepsilon}^{(1)}(r_0)&\lesssim \varepsilon^{-\sigma r_0}
\xi_{\varepsilon}^{\eta_i r_0}
\int_{|x|\leq \eta_{\varepsilon}}\left|\frac{1}{|x|^{\gamma_i}}\left(I_{\alpha}*|\tilde{w}_{\varepsilon}|^{\frac{N+\alpha}{N-2}}\right)(\frac{x}{|x|^2})
K[\tilde{w}_{\varepsilon}]^{q_i-2}(x)\right|^{r_0}dx\\
&\lesssim \varepsilon^{-\sigma r_0}
\xi_{\varepsilon}^{\eta_i r_0}
\left(\int_{|x|\leq \eta_{\varepsilon}}\frac{1}{|x|^{\iota_1r_0s_1}}K[\tilde{w}_{\varepsilon}]^{(q_i-2)r_0\theta s_1}dx\right)^{\frac{1}{s_1}}\\
&\ \ \ \cdot\left(\int_{|x|\leq \eta_{\varepsilon}}\frac{1}{|x|^{\iota_2r_0s_2}}
\left|\left(I_{\alpha}*|\tilde{w}_{\varepsilon}|^{\frac{N+\alpha}{N-2}}\right)(\frac{x}{|x|^2})\right|^{r_0s_2}dx\right)^{\frac{1}{s_2}}\\
&\ \ \ \cdot\left(\int_{|x|\leq \eta_{\varepsilon}}\left|K[\tilde{w}_{\varepsilon}]\right|^{(q_i-2)(1-\theta)r_0s_3}dx\right)^{\frac{1}{s_3}},\\
\end{aligned}
\end{equation}
where
$\gamma_i=\iota_1+\iota_2
=2N-(N-2)q_i$ with $\iota_1=(N+\alpha)-(N-2)q_i$ and $\iota_2=N-\alpha$.

Let
$$s_1=\frac{2N}{2N-(N-\alpha)(r_0+\kappa)},\ s_2=\frac{2N}{(N-\alpha)r_0}\ \ \mbox{and}\ \ \ s_3=\frac{2N}{(N-\alpha)\kappa},$$
it is easy to see that $\frac{1}{s_1}+\frac{1}{s_2}+\frac{1}{s_3}=1$.

Using the fact that $\iota_2r_0s_2=2N$ and the Hardy-Littlewood-Sobolev inequality, we get
\begin{equation}\label{4.20}
\begin{aligned}
&\ \ \ \int_{|x|\leq \eta_{\varepsilon}}\frac{1}{|x|^{\iota_2r_0s_2}}
\left|\left(I_{\alpha}*|\tilde{w}_{\varepsilon}|^{\frac{N+\alpha}{N-2}}\right)(\frac{x}{|x|^2})\right|^{r_0s_2}dx\\
&=\int_{|z|\geq \eta_{\varepsilon}^{-1}}\left|(I_{\alpha}*|\tilde{w}_{\varepsilon}|^{\frac{N+\alpha}{N-2}})(z)\right|^{r_0s_2}dz\\
&\lesssim \left(\int_{\mathbb R^N}|\tilde{w}_{\varepsilon}|^{\frac{N+\alpha}{N-2}\cdot\frac{Nr_0s_2}{r_0s_2\alpha+N}}dz\right)
^{\frac{r_0s_2\alpha+N}{Nr_0s_2}\cdot r_0s_2}=\left(\int_{\mathbb R^N}|\tilde{w}_{\varepsilon}|^{2^*}dz\right)
^{\frac{r_0s_2\alpha+N}{N}}\leq C.\\
\end{aligned}
\end{equation}
Also, for small $\kappa>0$, $(q_i-2)(1-\theta)r_0s_3>2^*$, hence from \eqref{4.17}, we have
\begin{equation}\label{4.21}
\int_{|x|\leq \eta_{\varepsilon}}|K[\tilde{w}_{\varepsilon}]|^{(q_i-2)(1-\theta)r_0s_3}dx\leq C<\infty.
\end{equation}
Therefore, by \eqref{4.19}-\eqref{4.21} and \eqref{4.8}, we have
$$\begin{aligned}
I_{\varepsilon}^{(1)}(r_0)&\lesssim \varepsilon^{-\sigma r_0} 
\xi_{\varepsilon}^{\eta_i r_0}
\left(\int_{|x|\leq \eta_{\varepsilon}}\frac{1}{|x|^{\iota_1r_0s_1}}K[\tilde{w}_{\varepsilon}]^{(q_i-2)r_0\theta s_1}dx\right)^{\frac{1}{s_1}}\\
&\lesssim \varepsilon^{-\sigma r_0-\frac{\sigma(N-2)(q_i-2)r_0\theta}{4}} 
\xi_{\varepsilon}^{\eta_i r_0+\frac{(N-2)(q_i-2)r_0\theta}{2}}\\
&\ \ \ \cdot\left(\int_{|x|\leq \eta_{\varepsilon}}\frac{\mathrm{exp}\left(-\frac{1}{2}(q_i-2)r_0\theta s_1\varepsilon^{-\sigma/2}\xi_{\varepsilon}|x|^{-1}\right)}
{|x|^{\iota_1r_0s_1+(N-2)(q_i-2)r_0\theta s_1}}dx\right)^{\frac{1}{s_1}}\\
&= \varepsilon^{-\Gamma_1(r_0,s_1)}
\cdot\xi_{\varepsilon}^{\Gamma_2(r_0,s_1)} \cdot\left(\int_{L_0}^{\infty}\frac{\mathrm{exp}\left(-\frac{1}{2}(q_i-2)r_0\theta s_1|R|\right)}
{|R|^{-\iota_1r_0s_1-(N-2)(q_i-2)r_0\theta s_1+N+1}}dR\right)^{\frac{1}{s_1}},\\
\end{aligned}$$
where
$$\Gamma_1(r_0,s_1)=\sigma r_0-\frac{\sigma(N-2)(q_i-2)r_0\theta}{4}-\frac{\sigma}{2}(\iota_1r_0-\frac{N}{s_1}),$$
and $$\Gamma_2(r_0,s_1)=\eta_i r_0-\frac{(N-2)(q_i-2)r_0\theta}{2}-\iota_1r_0+\frac{N}{s_1}.$$
Since $q_i>2$, it is easy to check that
$$\begin{aligned}
\Gamma_1(\frac{N}{2},\frac{4}{4-(N-\alpha)})
&=\frac{N\sigma}{2} \left(1+\frac{(N-2)(2-\theta)q_i}{4}+\frac{(N-2)\theta}{2}-\frac{3N+\alpha}{4}+1\right)\\
&>\frac{N\sigma}{2}\cdot \frac{N-\alpha}{4}>0,
\end{aligned}$$
and
$$\begin{aligned}
\Gamma_2(\frac{N}{2},\frac{4}{4-(N-\alpha)})
=&\frac{N(N-2)(1-\theta)q_i}{4}-\frac{N^2}{2}+\frac{N(N-2)\theta}{2}+N\\
>&\frac{N(N-2)(1-\theta)}{2}-\frac{N^2}{2}+\frac{N(N-2)\theta}{2}+N=0.\\
\end{aligned}$$
Therefore, we can choose $r_0>\frac{N}{2}$ with $r_0-\frac{N}{2}$ being small, and $\kappa>0$ small such that
$\Gamma_1(r_0,s_1)>0$ and $\Gamma_2(r_0,s_1)>0$.
Thus, we conclude that
$$I_{\varepsilon}^{(1)}(r_0)\lesssim \varepsilon^{-\Gamma_1(r_0,s_1)}\xi_{\varepsilon}^{\Gamma_2(r_0,s_1)} \ \to\ 0,\ \ \mbox{as}\ \ \varepsilon\to\infty.$$

For the second integral part, by H\"{o}lder inequality, we know
\begin{equation}\label{4.22}
\begin{aligned}
I_{\varepsilon}^{(2)}(r_0)&\lesssim \varepsilon^{-\sigma r_0}
\xi_{\varepsilon}^{\eta_i r_0}\int_{\eta_{\varepsilon}\leq |x|\leq4}\left|\frac{1}{|x|^{\gamma_i}}
\left(I_{\alpha}*|\tilde{w}_{\varepsilon}|^{\frac{N+\alpha}{N-2}}\right)(\frac{x}{|x|^2})
K[\tilde{w}_{\varepsilon}]^{q_i-2}(x)\right|^{r_0}dx\\
&\lesssim \varepsilon^{-\sigma r_0}
\xi_{\varepsilon}^{\eta_i r_0}\left(\int_{\eta_{\varepsilon}\leq |x|\leq4}\frac{1}{|x|^{\iota_1 r_0 \tilde{s}_1}}dx\right)^{\frac{1}{\tilde{s}_1}}\\
&\ \ \ \cdot\left(\int_{\eta_{\varepsilon}\leq |x|\leq4}\frac{1}{|x|^{\iota_2 r_0 s_2}}
\left|\left(I_{\alpha}*|\tilde{w}_{\varepsilon}|^{\frac{N+\alpha}{N-2}}\right)(\frac{x}{|x|^2})\right|^{r_0 s_2}dx\right)^{\frac{1}{s_2}}\\
&\ \ \ \cdot\left(\int_{\eta_{\varepsilon}\leq |x|\leq4}K[\tilde{w}_{\varepsilon}]^{(q_i-2)r_0 \tilde{s}_3}dx\right)^{\frac{1}{\tilde{s}_3}}.\\
\end{aligned}
\end{equation}
Let $\kappa>0$ be a small number, and $\tau(r_0)$ be given as follows
$$\tau(r_0)=\frac{2N-[(N-2)(q_i-2)+(N-\alpha)]r_0}{(N-2)(q_i-2)r_0}>0,$$
for
$r_0 \in \left(\frac{N}{2},\frac{2N}{(N-2)(q_i-2)+(N-\alpha)}\right)$.
Put $$\tilde{s}_1=\frac{2N(1+\kappa)^{-1}}{2N-[(N-2)(q_i-2)+(N-\alpha)]r_0},$$
$$s_2=\frac{2N}{(N-\alpha)r_0}\ \ \mbox{and} \ \ \tilde{s}_3=\frac{2N}{(N-2)(q_i-2)r_0(1-\kappa\tau(r_0))},$$
then it is easy to see that
$\frac{1}{\tilde{s}_1}+\frac{1}{s_2}+\frac{1}{\tilde{s}_3}=1$, $(q_i-2)r_0\tilde{s}_3=\frac{2N}{N-2}(1-\kappa\tau(r_0))^{-1}>2^*$, and hence
it follows from \eqref{4.17} that
$$\int_{\eta_{\varepsilon}\leq |x|\leq4}K[\tilde{w}_{\varepsilon}]^{(q_i-2)r_0 \tilde{s}_3}dx\leq C<\infty.$$
Since $2<q_i<\frac{N+\alpha}{N-2}$ and $r_0>\frac{N}{2}$, we have
$$
\iota_1r_0\tilde{s}_1-N=\frac{2Nr_0(1+\kappa)^{-1}[(N+\alpha)-(N-2)q_i]}{2N-[(N-2)(q_i-2)+(N-\alpha)]r_0}-N>\frac{N}{2},
$$
which implies that
$$
\begin{aligned}
&\quad\int_{\eta_{\varepsilon}\leq |x|\leq4}\frac{1}{|x|^{\iota_1 r_0 \tilde{s}_1}}dx=\int_{\eta_{\varepsilon}}^{4}r^{-\iota_1 r_0 \tilde{s}_1+N-1}dr\\
&=\frac{1}{\iota_1 r_0 \tilde{s}_1-N}
\left(L_0^{\iota_1 r_0 \tilde{s}_1-N}\varepsilon^{\frac{\sigma}{2}(\iota_1 r_0 \tilde{s}_1-N)}\xi_{\varepsilon}^{-\iota_1 r_0 \tilde{s}_1+N}
-4^{-\iota_1 r_0 \tilde{s}_1+N}\right)\\
&\lesssim \varepsilon^{\frac{\sigma}{2}(\iota_1 r_0 \tilde{s}_1-N)}\xi_{\varepsilon}^{-\iota_1 r_0 \tilde{s}_1+N}.\\
\end{aligned}
$$
Combining the above inequality with \eqref{4.20}, we have
$$I_{\varepsilon}^{(2)}(r_0)\lesssim \varepsilon^{-\sigma r_0+\frac{\sigma}{2}(\iota_1 r_0-\frac{N}{\tilde{s}_1})}\xi_{\varepsilon}^{\eta_ir_0-\iota_1 r_0+\frac{N}{\tilde{s}_1}}.$$

Let$$\Gamma_3(r_0,\tilde{s}_1):=\sigma r_0-\frac{\sigma}{2}(\iota_1 r_0-\frac{N}{\tilde{s}_1}),
\ \ \ \Gamma_4(r_0,\tilde{s}_1):=\eta_ir_0-\iota_1 r_0+\frac{N}{\tilde{s}_1},$$
then it is obvious that
$$
\begin{aligned}
\Gamma_3(\frac{N}{2},\tilde{s}_1)&>\Gamma_3\left(\frac{N}{2},\frac{4}{4-[(N-2)(q_i-2)+(N-\alpha)]}\right)\\
&=\frac{N\sigma}{4}\left(2+2-\frac{(N-2)(q_i-2)+(N-\alpha)}{2}-(N+\alpha)+(N-2)q_i\right)\\
&=\frac{N\sigma}{4}\left(\frac{(N-2)q_i}{2}-\frac{N+\alpha}{2}+2\right)>\frac{N\sigma}{4}\cdot\frac{N-\alpha}{2}>0,\\
\end{aligned}
$$
and
$$
\begin{aligned}
\Gamma_4(\frac{N}{2},\tilde{s}_1)&>\Gamma_4\left(\frac{N}{2},\frac{4}{4-[(N-2)(q_i-2)+(N-\alpha)]}\right)
=N-2\cdot\frac{N}{2}=0.
\end{aligned}
$$
Therefore, we can choose $r_0>\frac{N}{2}$ and $\kappa>0$ small such that $\Gamma_3(r_0,\tilde{s}_1)>0$, $\Gamma_4(r_0,\tilde{s}_1)>0$, then we obtain
$$I_{\varepsilon}^{(2)}(r_0)\lesssim \varepsilon^{-\Gamma_3(r_0,\tilde{s}_1)}\xi_{\varepsilon}^{\Gamma_4(r_0,\tilde{s}_1)}\ \to 0,\ \ \ \mbox{as}\ \ \ \varepsilon\to\infty.$$
Hence, $\lim_{\varepsilon\to\infty}I_{\varepsilon}(r_0)=\lim_{\varepsilon\to\infty}(I_{\varepsilon}^{(1)}(r_0)+I_{\varepsilon}^{(2)}(r_0))=0$.
The proof is complete.
\end{proof}
Similarly, we also have the following result.
\begin{lemma}\label{lem 2c}
Assume that (H1), (H2), (H3) and (H4) hold, then for any $i=1,2$,
\begin{equation}\label{4.13c}
\lim_{\varepsilon\to\infty}\int_{|x|\leq4}\left|\frac{\varepsilon^{-\sigma}
\xi_{\varepsilon}^{\eta_i}}{|x|^4}
\left(I_{\alpha}*|\tilde{w}_{\varepsilon}|^{q_i}\right)(\frac{x}{|x|^2})|\tilde w_\varepsilon|^{\frac{N+\alpha}{N-2}-2}(\frac{x}{|x|^2})\right|^{r_0}dx=0.
\end{equation}
\end{lemma}
\begin{lemma}\label{lem 5}
Assume that (H1), (H2), (H3) and (H4) hold, then there exists $r_0>\frac{N}{2}$ such that
\begin{equation}\label{4.24}
\begin{aligned}
\lim_{\varepsilon\to\infty}\int_{|x|\leq4}\left|\frac{\varepsilon^{-2\sigma}
\xi_{\varepsilon}^{\eta_i+\eta_j}}{|x|^4}
\left(I_{\alpha}*|\tilde{w}_{\varepsilon}|^{q_i}\right)(\frac{x}{|x|^2})|\tilde{w}_{\varepsilon}|^{q_j-2}(\frac{x}{|x|^2})\right|^{r_0}dx=0.\\
\end{aligned}
\end{equation}
\end{lemma}
\begin{proof}
Similar to Lemma \ref{lem 2}, we divide the integral into two parts: $I_{\varepsilon}^{(1)}\left(r_0\right)$ denotes the integral on $\{x|
 \ |x|\leq\eta_{\varepsilon}\}$, $I_{\varepsilon}^{(2)}\left(r_0\right)$ denotes the integral on $\{x|\ \eta_{\varepsilon}\leq|x|\leq4\}$.

For the first integral part, by H\"{o}lder inequality, we have
\begin{equation}\label{4.25}
\begin{aligned}
I_{\varepsilon}^{(1)}(r_0)&\lesssim \varepsilon^{-2\sigma r_0}
\xi_{\varepsilon}^{(\eta_i+\eta_j) r_0}
\int_{|x|\leq \eta_{\varepsilon}}\left|\frac{1}{|x|^{\gamma_j}}\left(I_{\alpha}*|\tilde{w}_{\varepsilon}|^{q_i}\right)(\frac{x}{|x|^2})
K[\tilde{w}_{\varepsilon}]^{q_j-2}(x)\right|^{r_0}dx\\
&\lesssim \varepsilon^{-2\sigma r_0}
\xi_{\varepsilon}^{(\eta_i+\eta_j) r_0}
\left(\int_{|x|\leq \eta_{\varepsilon}}\frac{1}{|x|^{\iota_1r_0s_1}}K[\tilde{w}_{\varepsilon}]^{(q_j-2)r_0\theta s_1}dx\right)^{\frac{1}{s_1}}\\
&\ \ \ \cdot\left(\int_{|x|\leq \eta_{\varepsilon}}\frac{1}{|x|^{\iota_2r_0s_2}}
\left|\left(I_{\alpha}*|\tilde{w}_{\varepsilon}|^{q_i}\right)(\frac{x}{|x|^2})\right|^{r_0s_2}dx\right)^{\frac{1}{s_2}}\\
&\ \ \ \cdot\left(\int_{|x|\leq \eta_{\varepsilon}}\left|K[\tilde{w}_{\varepsilon}]\right|^{(q_j-2)(1-\theta)r_0s_3}dx\right)^{\frac{1}{s_3}},\\
\end{aligned}
\end{equation}
where
$\gamma_j=\iota_1+\iota_2
=2N-(N-2)q_j$ with $\iota_1=(N+\alpha)-(N-2)q_j$ and $\iota_2=N-\alpha$.

Let
$$s_1=\frac{2N}{2N-(N-\alpha)(r_0+\kappa)},\ s_2=\frac{2N}{(N-\alpha)r_0}\ \ \mbox{and}\ \ \ s_3=\frac{2N}{(N-\alpha)\kappa},$$
it is easy to see that $\frac{1}{s_1}+\frac{1}{s_2}+\frac{1}{s_3}=1$.

Using the fact that $\iota_2r_0s_2=2N$ and the Hardy-Littlewood-Sobolev inequality, we get
\begin{equation}\label{4.20}
\begin{aligned}
&\ \ \ \int_{|x|\leq \eta_{\varepsilon}}\frac{1}{|x|^{\iota_2r_0s_2}}
\left|\left(I_{\alpha}*|\tilde{w}_{\varepsilon}|^{q_i}\right)(\frac{x}{|x|^2})\right|^{r_0s_2}dx\\
&=\int_{|z|\geq \eta_{\varepsilon}^{-1}}\left|(I_{\alpha}*|\tilde{w}_{\varepsilon}|^{q_i})(z)\right|^{r_0s_2}dz
\lesssim \left(\int_{\mathbb R^N}|\tilde{w}_{\varepsilon}|^{q_i\cdot\frac{Nr_0s_2}{r_0s_2\alpha+N}}dz\right)^{\frac{r_0s_2\alpha+N}{Nr_0s_2}\cdot r_0s_2}\\
&=\left(\int_{\mathbb R^N}|\tilde{w}_{\varepsilon}|^{\frac{2Nq_i}{N+\alpha}}dz\right)^{\frac{N+\alpha}{N-\alpha}}
 =\xi_{\varepsilon}^{\frac{N(N-2)q_i}{N-\alpha}-\frac{N(N+\alpha)}{N-\alpha}}
 \left(\int_{\mathbb R^N}|w_{\varepsilon}|^{\frac{2Nq_i}{N+\alpha}}dz\right)^{\frac{N+\alpha}{N-\alpha}}.\\
\end{aligned}
\end{equation}
Also, for small $\kappa>0$, $(q_j-2)(1-\theta)r_0s_3>2^*$, hence from \eqref{4.17}, we have
\begin{equation}\label{4.21a}
\int_{|x|\leq \eta_{\varepsilon}}|K[\tilde{w}_{\varepsilon}]|^{(q_j-2)(1-\theta)r_0s_3}dx\leq C<\infty.
\end{equation}
Therefore, by \eqref{4.25}-\eqref{4.21a} and \eqref{4.8}, we have
$$\begin{aligned}
I_{\varepsilon}^{(1)}(r_0)&\lesssim \varepsilon^{-2\sigma r_0}
\xi_{\varepsilon}^{(\eta_i+\eta_j)r_0+\frac{(N-2)q_ir_0}{2}-\frac{(N+\alpha)r_0}{2}}
\left(\int_{|x|\leq \eta_{\varepsilon}}\frac{1}{|x|^{\iota_1r_0s_1}}K[\tilde{w}_{\varepsilon}]^{(q_j-2)r_0\theta s_1}dx\right)^{\frac{1}{s_1}}\\
&\lesssim \varepsilon^{-2\sigma r_0-\frac{\sigma(N-2)(q_j-2)r_0\theta}{4}}
\xi_{\varepsilon}^{(\eta_i+\eta_j) r_0+\frac{(N-2)(q_j-2)r_0\theta}{2}+\frac{(N-2)q_ir_0}{2}-\frac{(N+\alpha)r_0}{2}}\\
&\ \ \ \cdot\left(\int_{|x|\leq \eta_{\varepsilon}}\frac{\mathrm{exp}\left(-\frac{1}{2}(q_j-2)r_0\theta s_1\varepsilon^{-\sigma/2}\xi_{\varepsilon}|x|^{-1}\right)}
{|x|^{\iota_1r_0s_1+(N-2)(q_j-2)r_0\theta s_1}}dx\right)^{\frac{1}{s_1}}\\
&= \varepsilon^{-\Gamma_1(r_0,s_1)}
\cdot\xi_{\varepsilon}^{\Gamma_2(r_0,s_1)}
\cdot\left(\int_{L_0}^{\infty}\frac{\mathrm{exp}\left(-\frac{1}{2}(q_j-2)r_0\theta s_1|R|\right)}
{|R|^{-\iota_1r_0s_1-(N-2)(q_j-2)r_0\theta s_1+N+1}}dR\right)^{\frac{1}{s_1}},\\
\end{aligned}$$
where
 $$\Gamma_1(r_0,s_1)=2\sigma r_0-\frac{\sigma(N-2)(q_j-2)r_0\theta}{4}-\frac{\sigma}{2}(\iota_1r_0-\frac{N}{s_1}),$$
and $$\Gamma_2(r_0,s_1)=(\eta_i+\eta_j)r_0-\frac{(N-2)(q_j-2)r_0\theta}{2}-\iota_1r_0+\frac{N}{s_1}+\frac{(N-2)q_ir_0}{2}-\frac{(N+\alpha)r_0}{2},$$
Since $q_j>2$, it is easy to check that
$$\begin{aligned}
\Gamma_1(\frac{N}{2},\frac{4}{4-(N-\alpha)})
&=\frac{N\sigma}{2} \left(2+\frac{(N-2)(2-\theta)q_j}{4}+\frac{(N-2)\theta}{2}-\frac{3N+\alpha}{4}+1\right)\\
&>\frac{N\sigma}{2}\left(\frac{N-\alpha}{4}+1\right)>0,
\end{aligned}$$
and
$$\begin{aligned}
\Gamma_2(\frac{N}{2},\frac{4}{4-(N-\alpha)})
=&\frac{N(N-2)(1-\theta)q_j}{4}-\frac{N^2}{2}+\frac{N(N-2)\theta}{2}+N\\
>&\frac{N(N-2)(1-\theta)}{2}-\frac{N^2}{2}+\frac{N(N-2)\theta}{2}+N=0.\\
\end{aligned}$$
Therefore, we can choose $r_0>\frac{N}{2}$ with $r_0-\frac{N}{2}$ being small, and $\kappa>0$ small such that
$\Gamma_1(r_0,s_1)>0$ and $\Gamma_2(r_0,s_1)>0$.
Thus, we conclude that
$$I_{\varepsilon}^{(1)}(r_0)\lesssim \varepsilon^{-\Gamma_1(r_0,s_1)}\xi_{\varepsilon}^{\Gamma_2(r_0,s_1)} \ \to\ 0,\ \ \mbox{as}\ \ \varepsilon\to\infty.$$


For the second the integral part, for $r_0\in(\frac{N}{2},\frac{2N}{2N-(N-2)q_j})$, put
$$s_1=\frac{2N}{[2N-(N-2)q_j]r_0},\ \ \ s_2=\frac{2N}{2N-[2N-(N-2)q_j]r_0}.$$
Clearly, $\frac{1}{s_1}+\frac{1}{s_2}=1$.
Denote $\gamma_j=2N-(N-2)q_j$, then by the H\"{o}lder inequality, we have
\begin{equation}\label{4.27}
\begin{aligned}
I_{\varepsilon}^{(2)}(r_0)&\lesssim\varepsilon^{-2\sigma r_0}\xi_{\varepsilon}^{(\eta_i+\eta_j)r_0}\int_{\eta_{\varepsilon}\leq|x|\leq4}
\left|\frac{1}{|x|^{\gamma_j}}\left(I_{\alpha}*|\tilde{w}_{\varepsilon}|^{q_i}\right)(\frac{x}{|x|^2})|K[\tilde{w}_{\varepsilon}]|^{q_j-2}(x)\right|^{r_0}dx\\
&\lesssim\varepsilon^{-2\sigma r_0}\xi_{\varepsilon}^{(\eta_i+\eta_j)r_0}
\left(\int_{\eta_{\varepsilon}\leq |x|\leq 4}
\left|\frac{1}{|x|^{\gamma_j}}\left(I_{\alpha}*|\tilde{w}_{\varepsilon}|^{q_i}\right)(\frac{x}{|x|^2})\right|^{r_0s_1}dx\right)^{\frac{1}{s_1}}\\
&\quad\quad\quad\quad\quad\quad\quad\quad\cdot
\left(\int_{\eta_{\varepsilon}\leq |x|\leq 4}\left|K[\tilde{w}_{\varepsilon}]\right|^{(q_j-2)r_0s_2}dx\right)^{\frac{1}{s_2}},\\
&\lesssim\varepsilon^{-2\sigma r_0}\xi_{\varepsilon}^{(\eta_i+\eta_j)r_0}
\left(\int_{\frac{1}{4}\leq |z|\leq \eta_{\varepsilon}^{-1}}
\left|\left(I_{\alpha}*|\tilde{w}_{\varepsilon}|^{q_i}\right)(z)\right|^{\frac{2N}{2N-(N-2)q_j}}dx\right)^{\frac{[2N-(N-2)q_j]r_0}{2N}},\\
\end{aligned}
\end{equation}
in the last inequality, we have used the fact that
$$\int_{\eta_{\varepsilon}\leq |x|\leq 4}\left|K[\tilde{w}_{\varepsilon}]\right|^{(q_j-2)r_0s_2}dx\leq C<\infty,$$
which follows from \eqref{4.17} and $ (q_j-2)r_0s_2=\frac{2Nr_0(q_j-2)}{2N-[2N-(N-2)q_j]r_0}>2^*$.

To proceed,  we split into two cases:

{\bf Case 1}: When $\alpha\leq N-2$, we have $\frac{2\alpha}{N-2}\leq2<q_i, q_j<\frac{N+\alpha}{N-2}$,
hence $$2\leq \frac{2Nq_i}{2N-(N-2)q_j+2\alpha}\leq\frac{2N}{N-2}.$$
By the boundedness of $w_\varepsilon$ in $H^1(\mathbb R^N)$ and the fact
$$(\eta_i+\eta_j)r_0+\frac{(N-2)q_ir_0}{2}-\frac{[2N-(N-2)q_j+2\alpha]r_0}{2}=0,$$
thus, we obtain
$$\begin{aligned}
I_{\varepsilon}^{(2)}(r_0)&\lesssim\varepsilon^{-2\sigma r_0}
\xi_{\varepsilon}^{(\eta_i+\eta_j)r_0+\frac{(N-2)q_ir_0}{2}-\frac{[2N-(N-2)q_j+2\alpha]r_0}{2}}\\
&\ \ \ \cdot\left(\int_{\frac{1}{4}\leq|z|\leq \eta_{\varepsilon}^{-1}}
|w_{\varepsilon}|^{\frac{2Nq_i}{2N-(N-2)q_j+2\alpha}}dz\right)^{\frac{[2N-(N-2)q_j+2\alpha]r_0}{2N}}\\
&\lesssim\varepsilon^{-2\sigma r_0} \to0\ \ \mbox{as}\ \ \varepsilon\to\infty.
\end{aligned}$$

{\bf Case 2}: When $\alpha>N-2$, we have $2<\frac{2\alpha}{N-2}<\frac{N+\alpha}{N-2}$. Same as Case2 of the second part of Lemma \ref{lem 2a}, we have
$$\begin{aligned}
Q_{\varepsilon}:&=\int_{\frac{1}{4}\leq|z|\leq \eta_{\varepsilon}^{-1}}
\left|(I_{\alpha}*|\tilde{w}_{\varepsilon}|^{q_i})\right|^{\frac{2N}{2N-(N-2)q_j}}dz\\
&\lesssim\varepsilon^{\frac{N\sigma(q-1)}{2q}}
\xi_{\varepsilon}^{-\frac{N(q-1)}{q}+\frac{N-2}{2}\cdot\frac{2Nq_i}{2N-(N-2)q_j}-\frac{N}{qs}\cdot\frac{2Nq}{2N-(N-2)q_j}}.\\
\end{aligned}$$
Therefore, from case 2, we get
$$\begin{aligned}
I_{\varepsilon}^{(2)}(r_0)&\lesssim \varepsilon^{-2\sigma r_0}
\varepsilon^{\frac{N\sigma(q-1)}{2q}\cdot\frac{[2N-(N-2)q_j]r_0}{2N}}
\cdot\xi_{\varepsilon}^{\Gamma(\alpha,q_i,q_j)}\\
&\lesssim \varepsilon^{-2\sigma r_0 (1-\frac{q-1}{2q}\cdot\frac{2N-(N-2)q_j}{4})}\to\ 0,\ \ \ \mbox{as}\ \ \ \varepsilon\to \infty.
\end{aligned}$$
where we have used the fact that $1-\frac{q-1}{2q}\cdot\frac{2N-(N-2)q_j}{4}>0$ for $q>1$ and
$$\begin{aligned}
\Gamma(\alpha,q_i,q_j)&=(\eta_i+\eta_j)r_0+\frac{[2N-(N-2)q_j]r_0}{2N}\\
\cdot&\left(-\frac{N(q-1)}{q}+\frac{N-2}{2}\cdot\frac{2Nq_i}{2N-(N-2)q_j}-\frac{N}{qs}\cdot\frac{2Nq}{2N-(N-2)q_j}\right)=0.\\
\end{aligned}$$
Thus, the conclusion follows since $\lim_{\varepsilon\to\infty}I_{\varepsilon}^{(1)}(r_0)+I_{\varepsilon}^{(2)}(r_0)=0$.
\end{proof}

Finally, let $\bar{a}(x)=\varepsilon^{-\sigma}\xi_{\varepsilon}^2,$
and
$$\begin{aligned}
\bar{b}(x)=&\left(I_{\alpha}*(|\tilde{w}_{\varepsilon}|^{\frac{N+\alpha}{N-2}}
+\varepsilon_1^{-1}\xi_{\varepsilon}^{\frac{N+\alpha}{2}}G(\varepsilon_2\xi_{\varepsilon}^{-\frac{N-2}{2}}\tilde{w}_{\varepsilon}))\right)\\
\cdot&\left(\frac{N+\alpha}{N-2}|\tilde{w}_{\varepsilon}|^{\frac{N+\alpha}{N-2}-2}
+\varepsilon_1^{-1}\varepsilon_2\xi_{\varepsilon}^{\frac{2+\alpha}{2}}g(\varepsilon_2\xi_{\varepsilon}^{-\frac{N-2}{2}}\tilde{w}_{\varepsilon})\tilde{w}^{-1}_{\varepsilon}\right).\\
\end{aligned}$$
Then equation \eqref{4.2} reads as
$$-\Delta \tilde{w}_{\varepsilon}+ \bar{a}(x)\tilde{w}_{\varepsilon}=\bar{b}(x)\tilde{w}_{\varepsilon}.$$
Applying again the Moser iteration, similar to the above argument,  we prove that
\begin{equation}\label{4.6a}
\sup_{\varepsilon \in (\varepsilon_0, \infty)}\parallel\tilde{w}_{\varepsilon}\parallel_{L^{\infty}(B_1)}< \infty.
\end{equation}
Combining \eqref{4.6} and \eqref{4.6a}, we conclude the proof of  Proposition\ref{Pro 4.1}.

\smallskip

\subsection{The Existence and Nonexistence of Ground State}
In this subsection,   for the readers'  convenience, we sketch of the proof for the existence of ground state solution to  the equation \eqref{1.1}.  Our existence result is the following theorem, which  proof is similar to that in \cite{LW-1}.

\begin{theorem}\label{th6.1}
Assume that (H1) and (H2) hold. If  $N=3$, we further assume $q_2>1+\alpha$, then the problem \eqref{1.1}
admits a positive ground state $u_{\varepsilon} \in H^1(\mathbb R^N)$, which is radially symmetric and radially nonincreasing.
\end{theorem}
Firstly, it is known \cite{LM2019, MS2015} that any weak solution of \eqref{1.1} in $H^1(\mathbb R^N)$ has additional regularity properties, which allows us to establish the Poho\v{z}aev identity for all finite energy solutions.
\begin{lemma}\label{lem 6.1}
Assume that (H1) and (H2) hold. If $u\in H^1(\mathbb R^N)$ is a solution of \eqref{1.1}, then $u\in W_{loc}^{2,q}(\mathbb R^N)$ for every $q>1$. Moreover, $u$  satisfies the Poho\v{z}aev identity
\begin{equation}\label{6.2}
\begin{aligned}
P_{\varepsilon}(u)&=
\frac{N-2}{2}\int_{\mathbb R^N}|\nabla u|^2+\frac{N\varepsilon}{2}\int_{\mathbb R^N}|u|^2\\
&-\frac{N+\alpha}{2}\int_{\mathbb R^N}(I_{\alpha}*(|u|^{\frac{N+\alpha}{N-2}}
+G(u)))(|u|^{\frac{N+\alpha}{N-2}}+G(u))=0.\\
\end{aligned}
\end{equation}
\end{lemma}

For any function $u\in H^1(\mathbb R^N)$ and $t\geq0$, define $u_{t} : \mathbb R^N \to \mathbb R$ as in Lemma \ref{lem 2.1}.
Then we have the following result.
\begin{lemma}\label{lem 6.2}
Assume   (H1) and  (H2) hold. Then for every $u\in H^1(\mathbb R^N)\setminus\{0\}$, there exists a unique $t_0>0$ such that $P_{\varepsilon}(u_{t_0})=0$. Moreover, $I_{\varepsilon}(u_{t_0})=\max_{t\geq0}I_{\varepsilon}(u_t)$.
\end{lemma}
\begin{proof}
By a direct calculation, we know that
$$
\begin{aligned}
f(t):=I_{\varepsilon}(u_{t})&=\frac{t^{N-2}}{2}\int_{\mathbb R^N}|\nabla u|^2+\frac{t^{N}\varepsilon}{2}\int_{\mathbb R^N}|u|^2\\
&-\frac{t^{N+\alpha}}{2}\int_{\mathbb R^N}(I_{\alpha}*(|u|^{\frac{N+\alpha}{N-2}}+G(u)))
(|u|^{\frac{N+\alpha}{N-2}}+G(u)).
\end{aligned}
$$
It's clear that $f(t)$ has a unique critical point $t_0$ which corresponding to its maximum. Hence, $I_{\varepsilon}(u_{t_0})=\max_{t\geq0}I_{\varepsilon}(u_t)$ and
$$
\begin{aligned}
0=f'(t_0)&=\frac{N-2}{2}t_0^{N-3}\int_{\mathbb R^N}|\nabla u|^2+\frac{N\varepsilon}{2}t_0^{N-1}\int_{\mathbb R^N}|u|^2\\
&-\frac{N+\alpha}{2}t_0^{N+\alpha-1}\int_{\mathbb R^N}(I_{\alpha}*(|u|^{\frac{N+\alpha}{N-2}}+G(u)))
(|u|^{\frac{N+\alpha}{N-2}}+G(u)).
\end{aligned}
$$
That is, $P_{\varepsilon}(u_{t_0})=0$. The proof is complete.
\end{proof}
Define $$m_{\varepsilon}=\inf_{u\in H^{1}(\mathbb R^N)\backslash\{0\}}\sup_{t\geq 0}J_{\varepsilon}(u_t).$$
Lemma \ref{lem 2.1} and Lemma \ref{lem 6.2} imply that there exists a positive and radially nonincreasing function sequence
$\{u_n\}\in H_r^1(\mathbb R^N)\setminus \{0\}$ such that
\begin{equation}\label{6.3}
I_{\varepsilon}(u_n)\to m_{\varepsilon},\ \ I'_{\varepsilon}(u_n)\to0 \ \ and\ \  P_{\varepsilon}(u_n)=0.
\end{equation}
For such a function sequence, we have the following result.

\begin{lemma}\label{lem 6.3}
Assume (H1) and (H2) hold, then the above  sequence $\{u_{n}\}$ is bounded in $H_r^1(\mathbb R^N)$.
\end{lemma}
\begin{proof}
For any fixed $\varepsilon>0$  and the sequence satisfying \eqref{6.3}, by \eqref{1.5}, we get
\begin{equation}\label{6.4}
\begin{aligned}
m_{\varepsilon}+1&\geq I_{\varepsilon}(u_n)-\frac{1}{2q_1}I'_{\varepsilon}(u_n)u_n\\
&=\left(\frac{1}{2}-\frac{1}{2q_1}\right)\int_{\mathbb R^N}|\nabla u_n|^2+\varepsilon|u_n|^2\\
&+\left(\frac{1}{2q_1}\frac{N+\alpha}{N-2}-\frac{1}{2}\right)\int_{\mathbb R^N}(I_{\alpha}*(|u_n|^{\frac{N+\alpha}{N-2}}+G(u_n)))|u_n|^{\frac{N+\alpha}{N-2}}\\
&+\frac{1}{2q_1}\int_{\mathbb R^N}(I_{\alpha}*(|u_n|^{\frac{N+\alpha}{N-2}}+G(u_n)))(g(u_n)u_n-q_1G(u_n))\\
&\geq\left(\frac{1}{2}-\frac{1}{2q_1}\right)\frac{1}{2}\int_{\mathbb R^N}|\nabla u_n|^2+\varepsilon|u_n|^2.\\
\end{aligned}
\end{equation}
Thus, the sequence $\{u_n\}$ is bounded in $H_r^1(\mathbb R^N)$ as $n \to \infty$.
\end{proof}

In the following, we will give the lower and upper estimates of $m_{\varepsilon}$. For any $\kappa>0$, we define
\begin{equation}\label{6.5}
u_{\kappa}(x)=\varphi(x)U_{\kappa}(x),
\end{equation}
where $\varphi(x)\in C_{c}^{\infty}(\mathbb R^N)$ is a cut-off function satisfying: $0\leq \varphi(x)\leq1$ in $\mathbb R^N$; $\varphi(x)=1$ in $B_1(0)$; $\varphi(x)=0$ in $\mathbb R^N\setminus B_2$. Here, $B_{r}$ denotes the ball in $\mathbb R^N$ of center at origin and radius $r$.
$$U_{\kappa}(x)=\kappa^{-\frac{N-2}{2}}U_1(\frac{x}{\kappa})=[N(N-2)]^{\frac{N-2}{4}}\left(\frac{\kappa^2}{\kappa^2+|x|^2}\right)^{\frac{N-2}{2}}.$$

By \cite{BN1983} (see also \cite{W1996}), we have the following estimates.
\begin{equation}\label{6.6}
\int_{\mathbb R^N}|\nabla u_{\kappa}|^2=S^{\frac{N}{2}}+O(\kappa^{N-2}), \ \ N\geq3,
\end{equation}
\begin{equation}\label{6.7}
\int_{\mathbb R^N}|u_{\kappa}|^2=\left\{
\begin{aligned}
&K_1\kappa^2+O(\kappa^{N}),
\ & if& \ N\geq5,\\
&K_1\kappa^2|\ln \kappa|+O(\kappa^{2}),
\ & if&\ N=4,\\
&K_1\kappa+O(\kappa^{2}),
\ & if&\ N=3,\\
\end{aligned}
\right.
\end{equation}
where $K_1>0$. By a direct calculation, there exists $K_2$, $K_3>0$ such that
\begin{equation}\label{6.8}
\int_{\mathbb R^N}(I_{\alpha}*|u_{\kappa}|^{q})|u_{\kappa}|^{q}\geq K_2\kappa^{-(N-2)q+N+\alpha},
\end{equation}
and
\begin{equation}\label{6.9}
\int_{\mathbb R^N}(I_{\alpha}*|u_{\kappa}|^{\frac{N+\alpha}{N-2}})|u_{\kappa}|^{q}\geq K_3\kappa^{\frac{N+\alpha}{2}-\frac{(N-2)q}{2}}.
\end{equation}
Moreover, similar as in \cite{GY} and \cite{GY2017}, by direct computation,
\begin{equation}\label{6.10}
\int_{\mathbb R^N}(I_{\alpha}*|u_{\kappa}|^{\frac{N+\alpha}{N-2}})|u_{\kappa}|^{\frac{N+\alpha}{N-2}}\geq
(C_{*}(N,\alpha))^{\frac{N}{2}}S_{\alpha}^{\frac{N+\alpha}{2}}+O(\kappa^{\frac{N+\alpha}{2}}).
\end{equation}

\begin{lemma}\label{lem 6.4}
Assume (H1) and (H2) hold. If $N=3$, we further assume $q_2>1+\alpha$, then
\begin{equation}\label{6.11}
0<m_{\varepsilon}<\frac{2+\alpha}{2(N-2)}\left(\frac{N-2}{N+\alpha}S_{\alpha}\right)^{\frac{N+\alpha}{2+\alpha}}.
\end{equation}
\end{lemma}
\begin{proof}
We use $u_{\kappa}$ to estimate $m_{\varepsilon}$, where $u_{\kappa}$ is defined in \eqref{6.5}.
By \eqref{6.3}, Hardy-Littlewood-Sobolev inequality and Sobolev imbedding theorem, there exists $C_1$, $C_2>0$ such that
$$
\begin{aligned}
0=P_{\varepsilon}(u_n)&=
\frac{N-2}{2}\int_{\mathbb R^N}|\nabla u_n|^2+\frac{N\varepsilon}{2}\int_{\mathbb R^N}|u_n|^2\\
&-\frac{N+\alpha}{2}\int_{\mathbb R^N}(I_{\alpha}*(|u_n|^{\frac{N+\alpha}{N-2}}
+G(u_n)))(|u_n|^{\frac{N+\alpha}{N-2}}+G(u_n))\\
&\geq C_1\|u_n\|_{H^1(\mathbb R^N)}^{2}-C_2\left(\|u_n\|_{H^1(\mathbb R^N)}^{\frac{2(N+\alpha)}{N-2}}+\|u_n\|_{H^1(\mathbb R^N)}^{2q_i}
+\|u_n\|_{H^1(\mathbb R^N)}^{\frac{(N+\alpha)}{N-2}+q_i}\right),\\
\end{aligned}
$$
where $i=1$, $2$, and which implies that there exists $C_3>0$ such that
\begin{equation}\label{6.12}
\|u_n\|_{H^1(\mathbb R^N)}\geq C_3.
\end{equation}
Combining \eqref{6.4} and \eqref{6.12}, we obtain that $m_{\varepsilon}>0$.
By Lemma \ref{lem 6.2}, there exists a unique $t_{\kappa}$
such that $P_{\varepsilon}((u_{\kappa})_{t_{\kappa}})=0$ and $I_{\varepsilon}((u_{\kappa})_{t_{\kappa}})=\sup_{t\geq0}I_{\varepsilon}((u_{\kappa})_{t})$.
Thus, $m_{\varepsilon}\leq \sup_{t\geq0}I_{\varepsilon}((u_{\kappa})_{t})$. By a direct calculation, we  have 
\begin{equation}\label{6.13}
\begin{aligned}
I_{\varepsilon}((u_{\kappa})_{t})&=\frac{t^{N-2}}{2}\int_{\mathbb R^N}|\nabla u_{\kappa}|^2+\frac{t^{N}\varepsilon}{2}\int_{\mathbb R^N}|u_{\kappa}|^2\\
&-\frac{t^{N+\alpha}}{2}\int_{\mathbb R^N}(I_{\alpha}*(|u_{\kappa}|^{\frac{N+\alpha}{N-2}}+G(u_{\kappa})))
(|u_{\kappa}|^{\frac{N+\alpha}{N-2}}+G(u_{\kappa}))\\
&\leq\frac{t^{N-2}}{2}\left(S^{\frac{N}{2}}+O(\kappa^{N-2})\right)
-\frac{t^{N+\alpha}}{2}\left((C_{*}(N,\alpha))^{\frac{N}{2}}S_{\alpha}^{\frac{N+\alpha}{2}}+O(\kappa^{\frac{N+\alpha}{2}})\right)\\
&\ \ \ -c_1t^{N+\alpha}\kappa^{\frac{N+\alpha}{2}-\frac{(N-2)q_i}{2}}-c_2t^{N+\alpha}\kappa^{-(N-2)q_i+N+\alpha}\\
&\ \ \ +\frac{t^{N}\varepsilon}{2}
\begin{cases}
\begin{split}
&K_1\kappa^2+O(\kappa^{N}),
\ & \ N\geq5,\\
&K_1\kappa^2|\ln \kappa|+O(\kappa^{2}),
\ &\ N=4,\\
&K_1\kappa+O(\kappa^{2}),
\ &\ N=3,\\
\end{split}
\end{cases}
\end{aligned}
\end{equation}
where $c_1$, $c_2>0$.
We claim that for every fixed $\varepsilon>0$ there exists $t_0$, $t_{1}>0$ independent of $\kappa$ such that $t_{\kappa}\in[t_0,t_{1}]$ for $\kappa>0$ small. Suppose by contradiction that $t_{\kappa}\to0$ or $t_{\kappa}\to\infty$ as $\kappa\to0$. Then \eqref{6.13} implies that $\sup_{t\geq0}I_{\varepsilon}((u_{\kappa})_{t})\leq0$ as $\kappa\to0$ and then $m_{\varepsilon}\leq 0$, which contradicts $m_{\varepsilon}>0$. Hence, the claim holds.

According $S_{\alpha}=\frac{S}{C_*(N,\alpha)^{\frac{N-2}{N+\alpha}}}$ (see \cite{GY} Lemma 1.2) and by direct calculation, we obtain
$$
\frac{t^{N-2}}{2}S^{\frac{N}{2}}
-\frac{t^{N+\alpha}}{2}(C_{*}(N,\alpha))^{\frac{N}{2}}S_{\alpha}^{\frac{N+\alpha}{2}}\leq \frac{2+\alpha}{2(N-2)}\left(\frac{N-2}{N+\alpha}S_{\alpha}\right)^{\frac{N+\alpha}{2+\alpha}}.$$
From (H4), we have
$$\sup_{t\geq0}I_{\varepsilon}((u_{\kappa})_{t})<\frac{2+\alpha}{2(N-2)}\left(\frac{N-2}{N+\alpha}S_{\alpha}\right)^{\frac{N+\alpha}{2+\alpha}}.$$
The proof is complete.
\end{proof}

\begin{proof}[Proof of Theorem \ref{th6.1}]
Let $\{u_n\}\in H_r^1(\mathbb R^N)\setminus \{0\}$ satisfying \eqref{6.3} is a positive and radially nonincreasing function sequence. By Lemma \ref{lem 6.2}, $\{u_n\}$ is bounded in $H_r^1(\mathbb R^N)$.
Hence, there exists a function $u \in H_r^1(\mathbb R^N)$ such that up to a subsequence,
$$u_{n}\rightharpoonup u\ \ \mbox{in}\ H^1(\mathbb R^N),\ \ \ u_{n}\rightarrow u\ \ \mbox{in}\ \ L^{q}(\mathbb R^N)\ \  \mbox{for any} \ q\in(2,2^*),$$
and
$$u_{n}\rightarrow u\ \ \mbox{in}\ \ L_{loc}^2(\mathbb R^N),\ \ \ u_{n}(x)\rightarrow u(x)\ \ \mbox{a.e.  on} \ \ \mathbb R^N.$$

First, we claim that $u\not\equiv 0$. Suppose by contradiction that $u\equiv 0$. By \eqref{6.3}, we get that
\begin{equation}\label{6.14}
\begin{aligned}
&\frac{N}{2}I'_{\varepsilon}(u_n)u_n-P_{\varepsilon}(u_n)\\
=&\int_{\mathbb R^N}|\nabla u_n|^2
+\left(\frac{N+\alpha}{2}-\frac{N}{2}\frac{N+\alpha}{N-2}\right)\int_{\mathbb R^N}(I_{\alpha}*(|u_n|^{\frac{N+\alpha}{N-2}}+G(u_n)))|u_n|^{\frac{N+\alpha}{N-2}}\\
+&\int_{\mathbb R^N}(I_{\alpha}*(|u_n|^{\frac{N+\alpha}{N-2}}+G(u_n)))\left(\frac{N+\alpha}{2}G(u_n)-\frac{N}{2}g(u_n)u_n\right)\rightarrow 0.
\end{aligned}
\end{equation}
According  to $u_{n}\to 0$ in $L^q(\mathbb R^N)$, we have
$$\int_{\mathbb R^N}|\nabla u_n|^2=\frac{N+\alpha}{N-2}\int_{\mathbb R^N}(I_{\alpha}*(|u_n|^{\frac{N+\alpha}{N-2}}))|u_n|^{\frac{N+\alpha}{N-2}}+o(1)\to l.$$
Therefore, by \eqref{1}, we have $l\geq S_{\alpha} \left(\frac{N-2}{N+\alpha}l\right)^{\frac{N-2}{N+\alpha}}$, and we deduce that
$$either\ \  l=0,\ \  or\ \ l\geq \left(\frac{N-2}{N+\alpha}\right)^{\frac{N-2}{2+\alpha}}S_{\alpha}^{\frac{N+\alpha}{2+\alpha}}.$$
Let us suppose at first that $l\geq \left(\frac{N-2}{N+\alpha}\right)^{\frac{N-2}{2+\alpha}}S_{\alpha}^{\frac{N+\alpha}{2+\alpha}}$.
Since $I_{\varepsilon}(u_n)\to m_{\varepsilon}$ and $P_{\varepsilon}(u_n)=0$, we have that
$$
\begin{aligned}
m_{\varepsilon}+o(1)&=I_{\varepsilon}(u_n)=I_{\varepsilon}(u_n)-\frac{1}{N}P_{\varepsilon}(u_n)\\
&=\left(\frac{1}{2}-\frac{N-2}{2N}\right)\int_{\mathbb R^N}|\nabla u_n|^2\\
&+\left(\frac{N+\alpha}{2N}-\frac{1}{2}\right)
\int_{\mathbb R^N}(I_{\alpha}*(|u_n|^{\frac{N+\alpha}{N-2}}
+G(u_n)))(|u_n|^{\frac{N+\alpha}{N-2}}+G(u_n))\\
&=\frac{1}{N}\int_{\mathbb R^N}|\nabla u_n|^2+\frac{\alpha}{2N}\int_{\mathbb R^N}(I_{\alpha}*(|u_n|^{\frac{N+\alpha}{N-2}}))|u_n|^{\frac{N+\alpha}{N-2}}+o(1)\\
&=\frac{1}{N}l+\frac{\alpha}{2N}\frac{N-2}{N+\alpha}l+o(1) \geq \frac{2+\alpha}{2(N-2)}\left(\frac{N-2}{N+\alpha}S_{\alpha}\right)^{\frac{N+\alpha}{2+\alpha}}.\\
\end{aligned}
$$
However, this is not possible since in this case $m_{\varepsilon}\geq \frac{2+\alpha}{2(N-2)}\left(\frac{N-2}{N+\alpha}S_{\alpha}\right)^{\frac{N+\alpha}{2+\alpha}}$, which contradicts Lemma \ref{lem 6.4}.
If instead $l=0$, we have $\|\nabla u_n\|_2\to 0$ and $\int_{\mathbb R^N}(I_{\alpha}*(|u_n|^{\frac{N+\alpha}{N-2}}))|u_n|^{\frac{N+\alpha}{N-2}}\to0$.
Then $I_{\varepsilon}(u_n)\to 0$, which contradicts Lemma \ref{lem 6.4}. Thus, $u \not\equiv 0$.

Second, similar to the proof of Theorem 1.1 in \cite{LM2019}, for any $\varphi \in C_{c}^{\infty}(\mathbb R^N)$, we have
$$
\begin{aligned}
0=&\int_{\mathbb R^N}\nabla u_n\nabla\varphi+\varepsilon u_n\varphi
-\int_{\mathbb R^N}(I_{\alpha}*(|u_n|^{\frac{N+\alpha}{N-2}}+G(u_n)))(|u_n|^{\frac{4+\alpha-N}{N-2}}u_n\varphi+g(u_n)\varphi)\\
\to&\int_{\mathbb R^N}\nabla u\nabla\varphi+\varepsilon u\varphi
-\int_{\mathbb R^N}(I_{\alpha}*(|u|^{\frac{N+\alpha}{N-2}}+G(u)))(|u|^{\frac{4+\alpha-N}{N-2}}u\varphi+g(u)\varphi)
\end{aligned}
$$
as $n\to \infty$. Passing to the limit by weak convergence, we obtain
$$-\Delta u+ \varepsilon u=\left(I_{\alpha}*(|u|^{\frac{N+\alpha}{N-2}}+G(u))\right)\left(\frac{N+\alpha}{N-2}|u|^{\frac{N+\alpha}{N-2}-2}u+g(u)\right),$$
and hence by Poho\v{z}aev identity $P_{\varepsilon}(u)=0$.
Recalling that $\bar{u}_n=u_n-u\to0$ in $H^1(\mathbb R^N)$, we also have
\begin{equation}\label{6.15}
\|\nabla u_n\|_2^2=\|\nabla u\|_2^2+\|\nabla \bar{u}_n\|_2^2+o(1),\ \ \|u_n\|_2^2=\|u\|_2^2+\|\bar{u}_n\|_2^2+o(1),
\end{equation}
and by the Brezis-Lieb Lemma
\begin{equation}\label{6.16}
\begin{aligned}
&\int_{\mathbb R^N}(I_{\alpha}*(|u_n|^{\frac{N+\alpha}{N-2}}))|u_n|^{\frac{N+\alpha}{N-2}}\\
=&\int_{\mathbb R^N}(I_{\alpha}*(|u|^{\frac{N+\alpha}{N-2}}))|u|^{\frac{N+\alpha}{N-2}}
+\int_{\mathbb R^N}(I_{\alpha}*(|\bar{u}_n|^{\frac{N+\alpha}{N-2}}))|\bar{u}_n|^{\frac{N+\alpha}{N-2}}+o(1).\\
\end{aligned}
\end{equation}
Therefore, since $P(u_n)=0$, and $u_n\to u$ in $L^q(\mathbb R^N)$, we deduce that
$$
\begin{aligned}
&\frac{N}{2}I'_{\varepsilon}(u_n)u_n-P_{\varepsilon}(u_n)\\
=&\|\nabla u_n\|_2^2-\left(\frac{N}{2}\frac{N+\alpha}{N-2}-\frac{N+\alpha}{2}\right)
\int_{\mathbb R^N}(I_{\alpha}*(|u_n|^{\frac{N+\alpha}{N-2}}+G(u_n)))|u_n|^{\frac{N+\alpha}{N-2}}\\
-&\int_{\mathbb R^N}(I_{\alpha}*(|u_n|^{\frac{N+\alpha}{N-2}}+G(u_n)))\left(\frac{N+\alpha}{2}G(u_n)-\frac{N}{2}g(u_n)u_n\right)=o(1),
\end{aligned}
$$
hence
$$
\begin{aligned}
\|\nabla u\|_2^2+&\|\nabla \bar{u}_n\|_2^2=\frac{N+\alpha}{N-2}\int_{\mathbb R^N}(I_{\alpha}*G(u))|u|^{\frac{N+\alpha}{N-2}}\\
+&\frac{N+\alpha}{N-2}\Big(\int_{\mathbb R^N}(I_{\alpha}*(|u|^{\frac{N+\alpha}{N-2}}))|u|^{\frac{N+\alpha}{N-2}}
+\int_{\mathbb R^N}(I_{\alpha}*(|\bar{u}_n|^{\frac{N+\alpha}{N-2}}))|\bar{u}_n|^{\frac{N+\alpha}{N-2}}\Big)\\
+&\int_{\mathbb R^N}(I_{\alpha}*(|u|^{\frac{N+\alpha}{N-2}}+G(u)))\left(\frac{N+\alpha}{2}G(u)-\frac{N}{2}g(u)u\right)+o(1).\\
\end{aligned}
$$
But $P_{\varepsilon}(u)=0$, $I'(u)u=o(1)$ and hence $\|\nabla \bar{u}_n\|_2^2
=\frac{N+\alpha}{N-2}\int_{\mathbb R^N}(I_{\alpha}*(|\bar{u}|^{\frac{N+\alpha}{N-2}}))|\bar{u}|^{\frac{N+\alpha}{N-2}}+o(1)$.
We infer that up to a subsequence
$$\lim_{n\to\infty}\|\nabla \bar{u}_n\|_2^2
=\frac{N+\alpha}{N-2}\lim_{n\to\infty}\int_{\mathbb R^N}(I_{\alpha}*(|\bar{u}|^{\frac{N+\alpha}{N-2}}))|\bar{u}|^{\frac{N+\alpha}{N-2}}=l\geq0,$$
by \eqref{1}, we have $l\geq S_{\alpha} \left(\frac{N-2}{N+\alpha}l\right)^{\frac{N-2}{N+\alpha}}$.
Thus, either $l=0$, or $l\geq \left(\frac{N-2}{N+\alpha}\right)^{\frac{N-2}{2+\alpha}}S_{\alpha}^{\frac{N+\alpha}{2+\alpha}}$.
If $l\geq \left(\frac{N-2}{N+\alpha}\right)^{\frac{N-2}{2+\alpha}}S_{\alpha}^{\frac{N+\alpha}{2+\alpha}}$, then by \eqref{6.15} and \eqref{6.16}
$$
\begin{aligned}
m_{\varepsilon}&=\lim_{n\to\infty}I_{\varepsilon}(u_n)\\
&=\lim_{n\to\infty}\left(I_{\varepsilon}(u)+\frac{1}{2}\|\bar{u}_n\|_2^2+\frac{1}{2}\|\nabla \bar{u}_n\|_2^2
-\frac{1}{2}\int_{\mathbb R^N}(I_{\alpha}*(|\bar{u}_n|^{\frac{N+\alpha}{N-2}}))|\bar{u}_n|^{\frac{N+\alpha}{N-2}}\right)\\
&=I_{\varepsilon}(u)+\frac{1}{2}\|\bar{u}_n\|_2^2+\frac{2+\alpha}{2(N+\alpha)}l\geq \frac{2+\alpha}{2(N-2)}\left(\frac{N-2}{N+\alpha}S_{\alpha}\right)^{\frac{N+\alpha}{2+\alpha}},\\
\end{aligned}
$$
which contradicts Lemma \ref{lem 6.4}. Hence $l=0$, $\|\bar{u}\|_2^2=0$ and $I_{\varepsilon}(u)=m_{\varepsilon}$. We show that $u_n\rightarrow u$
in $H^1(\mathbb R^N)$. That is, $u$ ia a nonnegative and radially nonincreasing groundstate solution of \eqref{1.1}. The strongly maximum principle implies that
$u$ is positive.
The proof is complete. 
\end{proof}

\medskip
\medskip

\end{document}